\documentclass[11pt]{article}
\usepackage{amssymb,amsmath,amsthm,verbatim,IEEEtrantools,macros,bm}
\usepackage[margin=1in]{geometry}
\newcommand{\Rho}{P}
\newcommand{\Tau}{T}
\DeclareMathOperator{\realpart}{Re}
\DeclareMathOperator{\supp}{supp}
\DeclareMathOperator{\tr}{tr}
\DeclareMathOperator{\disc}{disc}
\DeclareMathOperator{\diam}{diam}
\title{Oscillatory integrals with phases arising from algebraic number fields}
\author{Robert Fraser}
\begin{document}
\maketitle
\begin{abstract}
We develop a theory of oscillatory integrals whose phase is given by the trace of a polynomial over an algebraic number field. We present an application to the singular integral for a version of Tarry's problem for algebraic integers.
\end{abstract}
\section{Introduction and Background}
\subsection{Oscillatory Integrals associated to real and complex polynomials}
In harmonic analysis, oscillatory integrals of the form
\begin{equation}\label{eq:Ipsiphi}
I_{\psi}(\phi) = \int_{\vec x \in \mathbb{R}^n} e^{2 \pi i \phi(\vec x)} \psi(\vec x) \, d \vec x
\end{equation}
are ubiquitous in harmonic analysis. Here, the function $\phi$ is a real-valued function called the \textit{phase} and the function $\psi$ is called a \textit{cutoff}. Integrals of this type are of great importance in Fourier analysis, partial differential equations, and number theory. An excellent resource on integrals of this type is Chapter VIII of Stein \cite{Stein93}.

When the phase $\phi$ is polynomial, it is possible to obtain bounds depending on the algebraic properties of the function $\phi$.

A rich theory has emerged in estimating oscillatory integrals whose phase $\phi$ has an additional complex structure. Problems related to restriction and convolution estimates for complex curves are considered by Bak and Ham \cite{BakHam14}, Chung and Ham \cite{ChungHam19}, Meade \cite{Meade23}, de Dios Pont \cite{deDiosPont20}, and Biggs, Brandes, and Hughes \cite{BiggsBrandesHughes22}. A general theory of oscillatory integrals whose phase has complex structure is presented by Wright \cite{Wright20} for complex curves. We will describe this theory in detail in Subsection \ref{sub:cplx}.

We generalize this theory of complex oscillatory integrals to oscillatory integrals whose phase has an additional structure corresponding to an algebraic number field. 

We will present an application of our main result to the oscillatory integral arising in Tarry's problem. Tarry's problem concerns the existence of nontrivial integer solutions to the Vinogradov system
\begin{IEEEeqnarray*}{rCl}
x_1^1 + \cdots + x_s^1 & = & y_1^1 + \cdots + y_s^1 \\
x_1^2 + \cdots + x_s^2 & = & y_1^2 + \cdots + y_s^2 \\
\vdots & \vdots & \vdots \\
x_1^k + \cdots + x_s^d & = & y_1^d + \cdots + y_s^d.
\end{IEEEeqnarray*}
The singular integral in Tarry's problem is 
\[\int_{\vec \eta \in \mathbb{R}^d}  \left| \int_{t \in [0,1]}  e^{2 \pi i (\eta_1 t + \eta_2 t^2 + \cdots + \eta_d t^d)}\right|^{2s} \, dt  \, d\vec \eta.\]
Estimating this singular integral is important in determining the minimal value of $d$ for which the Vinogradov system has nontrivial integer solutions. It therefore is important to obtain estimates on the extension operator associated to the moment curve:
\[E f(\vec \eta)= \int_{t \in \mathbb{R}} e^{2 \pi i(\eta_1 t + \cdots + \eta_d t^d)} \psi(t) f(t) \, dt\]
where $\psi$ is a smooth cutoff. The $L^p$ mapping properties of this operator were determined by Drury \cite{Drury85}.

In this work, we will consider the singular integral arising in a version of Tarry's problem that counts the number of solutions to the Vinogradov system lying in an algebraic number field. The proof is similar to Wright's estimate \cite{Wright20} of the singular integral for Tarry's problem over the Gaussian integers.
\subsection{Oscillatory integrals and the $H$-functional}
Let $\psi : \mathbb{R}^n \to \mathbb{R}$ be a smooth cutoff function and let $\phi : \mathbb{R}^n \to \mathbb{R}$ be a smooth phase.

Let $I_{\psi}(\phi)$ be the oscillatory integral defined by
\[I_{\psi}(\phi) := \int e^{2 \pi i \phi(\vec x)} \psi(\vec x) \, d\vec x.\]
For a fixed $\psi$, it is interesting to try to obtain bounds on the integral $I_{\psi}(\phi)$ in terms of the derivatives of $\phi$.

In one dimension, a useful bound is given by the van der Corput lemma. This lemma can be found in Stein \cite{Stein93}.

\begin{mylem}\label{lem:vdc}
Let $r \geq 2$ and suppose $\phi$ is a phase function such that $|\phi^{(r)}(x)| \geq 1$ on an interval $[a,b]$. Then
\begin{equation}\label{eq:vdc}
\int_a^b e^{i \lambda \phi(x)} \lesssim \lambda^{-r}
\end{equation}
where the implicit constant does not depend on the phase $\phi$. If $\phi'$ is monotone, then \eqref{eq:vdc} holds for $r = 1$ as well.
\end{mylem}

In higher dimensions, there is a related but inferior bound. The version of the bound presented below can be found in Stein \cite{Stein93}.
\begin{mylem}\label{lem:vdc2}
Let $\psi$ be a smooth function supported in the unit ball. Let $r \in \mathbb{N}$ and let $\phi$ be a real-valued $C^{r+1}$ function satisfying the condition that for some multi-index $\alpha$ with $|\alpha| = r> 0$, we have the lower bound 
\[\left| \partial^{\alpha} \phi(x) \right| \geq 1\]
for all $x$ in the support of $\psi$. Then
\begin{equation}\label{eq:vdc2}
\int e^{i \lambda \phi(x)} \psi(x) \, dx \leq c_r(\phi) \lambda^{-n/r} (\norm{\psi}_{\infty} + \norm{\nabla \psi}_{1})
\end{equation}
where $c_r(\phi)$ is a constant that depends on $\phi$, but remains bounded provided that the $C^{r+1}$ norm of $\phi$ remains bounded.
\end{mylem}
Although Lemma \ref{lem:vdc2} is useful, the dependence of the constant on $\phi$ limits its usefulness. In practice, we will often want to obtain a bound on the oscillatory integral in \eqref{eq:vdc2} that is uniform over a family of functions $\phi$ whose $C^{r+1}$ norms do not remain bounded. 

Although it is not possible to obtain a bound that is uniform over all $C^{r + 1}$ functions, a classical result of Arhipov, Čubarikov, and Karacuba \cite{ArhipovKaracubaCubarikov79} establishes a bound that is uniform over polynomials of bounded degree. In order to state this bound, we need to introduce a quantity called the $H$-functional.

\begin{mydef}\label{def:hrdef}
Let $f : \mathbb{R}^n \to \mathbb{R}$ be a polynomial function. For points $x \in \mathbb{R}^n$, define the \textbf{real pointwise $H$-functional} $H_{\mathbb{R}} f (x)$ to be the quantity
\begin{equation}\label{eq:hrxdef}
H_{\mathbb{R}} f(x) := \max_{|\alpha| \geq 1} \left|\frac{1}{\alpha!} \partial^{\alpha} f(x) \right|^{1/|\alpha|} .
\end{equation}
Let $\psi : \mathbb{R}^n \to \mathbb{R}$ be a compactly supported smooth function. Define the \textbf{real uniform $H$-functional} $H_{\mathbb{R}, \psi} f$ to be
\begin{equation}\label{eq:hrdef}
H_{\mathbb{R}, \psi} f := \inf_{x \in \supp \psi} H_{\mathbb{R}} f(x).
\end{equation}
\end{mydef}
The uniform $H$-functional $H_{\mathbb{R}, \psi} f$ will be large if, for every $x \in \supp \psi$, at least one of the derivatives $\partial^{\alpha} f$ is large. Observe that since each derivative $\partial^{\alpha}$ is linear, we have for $\lambda > 0$:
\[\left| \partial^{\alpha} \lambda f (x) \right|^{1 /|\alpha|} = \lambda^{1/\alpha} \left| \partial^{\alpha} f \right|^{1/|\alpha|}.\]
Hence, if $|\partial^{\alpha} f| \gtrsim 1$ on $\supp \psi$, we have the estimate
\begin{equation}\label{eq:hgrowth}
H_{\mathbb{R}, \psi} \lambda f \gtrsim \lambda^{1/\alpha},
\end{equation}
where the implicit constant is independent of the function $f$. Arhipov, Čubarikov and Karacuba \cite{ArhipovKaracubaCubarikov79} establish the following bound. 
\begin{mythm}\label{thm:hrbound}[Arhipov, Karacuba, Čubarikov]
Let $d \geq 0$ and let $\phi$ be a polynomial of degree at most $d$. Then the following estimate holds:
\begin{equation}\label{eq:hrbound}
\int_{\mathbb{R}^n} e^{i \phi(x)} \psi(x) \, dx \lesssim_{\psi, d} (H_{\mathbb{R}, \psi} \phi )^{-1}.
\end{equation}
In particular, applying \eqref{eq:hgrowth} we recover the following bound. If $\phi$ is such that $|\partial^{\alpha} \phi| \geq 1$ on $\supp \psi$, then 
\begin{equation}\label{eq:sibound}
\int_{\mathbb{R}^n} e^{i \lambda \phi(x)} \psi(x) \, dx  \lesssim_{\psi, d} \lambda^{-1/|\alpha|}.
\end{equation}
\end{mythm}
\subsection{Complex Oscillatory Integrals}\label{sub:cplx}
Wright \cite{Wright20} discusses oscillatory integrals of the following type. Suppose $f: \mathbb{C} \to \mathbb{C}$ is a complex--analytic function. Then we can define the real and imaginary parts $u_f$ and $v_f$ of $f$ by
\begin{equation}\label{eq:ufvfdef}
f(x + iy) = u_f(x,y) + i v_f(x,y).
\end{equation}
More generally, if $f$ is an analytic function from $\mathbb{C}^n$ into $\mathbb{C}$, we dan define $u_f$ and $v_f$ by
\[f(x_1 + iy_1, \ldots, x_n + i y_n) = u_f(x_1, y_1, \ldots, x_n, y_n) + i v_f(x_1, y_1, \ldots, x_n, y_n).\]
Fix a smooth, compactly supported cutoff function $\psi : \mathbb{R}^{2n} \to \mathbb{R}$. Wright considers oscillatory integrals of the form
\[\int e^{i \phi_{\text{W}, f}(x_1, y_1, \ldots, x_n, y_n)} \psi(x_1, y_1, \ldots, x_n, y_n) \, dx_1 d y_1 \cdots d x_n d y_n \]
where the phase $\phi_{\text{W},f}$ is of the form 
\begin{equation}\label{eq:phiwfdef}
\phi_{\text{W}, f} := u_f + v_f
\end{equation}
for a complex polynomial $f$. Note that the functions that can be written as $\phi_{W,f}$ for a complex polynomial $f$ are precisely the harmonic polynomials.

Wright \cite{Wright20} defines a complex version of the $H$-functional.
\begin{mydef}\label{def:hcdef}
Let $f : \mathbb{C}^n \to \mathbb{C}$ be a polynomial function. For points $\bm{\vec x} = (x_1, y_1, \ldots, x_n, y_n) \in \mathbb{R}^{2n}$, define the \textbf{complex pointwise $H$-functional} $H_{\mathbb{C}} f (\bm{\vec x})$ to be the quantity \footnote{Wright defines the complex $H$-functional with exponent $1/|\alpha|$ instead of $2/|\alpha|$, but to conform to the convention in this paper, it is more convenient to use $2/|\alpha|$ as the exponent.}
\begin{equation}\label{eq:hcxdef}
\max_{|\alpha| \geq 1} \left|\frac{1}{\alpha!} \partial^{\alpha} f(x_1 + i y_1, \ldots, x_n + i y_n) \right|^{2/|\alpha|} .
\end{equation}
where the maximum is taken over the complex derivatives of $f$.

Let $\psi : \mathbb{R}^n \to \mathbb{R}$ be a compactly supported smooth function. Define the \textbf{complex uniform $H$-functional} $H_{\mathbb{C}, \psi} f$ to be
\begin{equation}\label{eq:hcdef}
\inf_{\bm{\vec x} \in \supp \psi} H_{\mathbb{C}} f(\bm{\vec x})
\end{equation}
\end{mydef}
Wright \cite{Wright20} establishes the following oscillatory integral bound.
\begin{mythm}[Wright]\label{thm:hcbound}
Let $f$ be a polynomial of degree at most $d$ with complex coefficients. Define $\phi_{\text{W},f}$ as in \eqref{eq:phiwfdef}. Let $\psi : \mathbb{R}^{2n} \to \mathbb{R}$ be a compactly supported smooth cutoff. Then 
\begin{equation}\label{eq:hcbound}
\int e^{i \phi_{W, f}(\bm{\vec x})} \psi(\bm{\vec x}) \, d \bm{\vec x} \lesssim_{\psi, d} (H_{\mathbb{C}, \psi}f)^{-1}
\end{equation}
\end{mythm}
This result is interesting for two reasons. First, the behavior of the oscillatory integral with \textit{real} phase $\phi_{W,f}$ of $2n$ real variables can be entirely understood in terms of the \textit{complex} function $f$ of $n$ complex variables. Second, the exponent $1/|\alpha|$ for real phases from the definition of the $H$-functional used in Theorem \ref{thm:hrbound} is improved to $2/|\alpha|$ for phases $\phi_{W,f}$ in the definition of the $H$-functional used in Theorem \ref{thm:hcbound}. This improved exponent is useful in applications.

In a talk given on this result and a later work, Dall'ara and Wright \cite{DallaraWright21} mention that Theorem \ref{thm:hcbound} holds if $\phi_{W,f}$ is replaced by any nonzero linear combination of $u_f$ and $v_f$. In order to conform to a convention we will use later in this manuscript, we will define $\phi_{\mathbb{C}, f}$ by the equation
\begin{equation}\label{eq:phicdef}
\phi_{\mathbb{C},f}(\bm{\vec x}) := 2 u_f(\bm{\vec x})
\end{equation}
and emphasize that Theorem \ref{thm:hcbound} holds when $\phi_{W, f}$ is replaced by $\phi_{\mathbb{C}, f}$.

Dall'ara and Wright \cite{DallaraWright21} make another interesting observation regarding Theorem \ref{thm:hcbound}. The definitions of $\phi_{\text{W},f}$ and $\phi_{\mathbb{C}, f}$ rely on the decomposition \eqref{eq:ufvfdef} of the function $f$ with respect to the basis $\{1, i\}$. However, the proof of Theorem \ref{thm:hcbound} does not depend on the choice of basis for $\mathbb{C}$. Given a vector space basis $B := \{\omega_1, \omega_2\}$ of $\mathbb{C}$ as an $\mathbb{R}$-vector space, we can write a complex number $z$ in a unique way in the form $z = x_1 \omega_1 + x_2 \omega_2$ where $x_1, x_2 \in \mathbb{R}$. Hence, given a complex polynomial $f: \mathbb{C}^n \to \mathbb{C}$, we can define $u_B$ and $v_B$ by
\[f(x_1 \omega_1 + y_1 \omega_2, \ldots, x_n \omega_1 + y_n \omega_2) := u_B(\vec x, \vec y) \omega_1 + v_B(\vec x, \vec y) \omega_2.\]
If we define the real-valued function $\psi_{W, B, f}$ by 
\[\phi_{\text{W}, B, f}(\vec x, \vec y) := u_B(\vec x, \vec y) + v_B(\vec x, \vec y),\]
then Theorem \ref{thm:hcbound} holds with $\phi_{W,f}$ replaced by $\phi_{W, B, f}$.

The primary goal of this paper is to locate a family of polynomial phases in $\mathbb{R}^{kn}$ and an appropriate $H$-functional for which the bound \eqref{eq:hcbound} can be improved.
\subsection{Algebraic number fields}
In order to introduce the class of oscillatory integrals to be discussed in this article, we must first discuss algebraic number fields. An \textbf{algebraic number field} is a finite extension of $\mathbb{Q}$. If $K$ is an algebraic number field, then the \textbf{degree} of $K$ is the dimension of $K$ as a $\mathbb{Q}$-vector space. Throughout our discussion, we will fix a number field $K$ of degree $k$ and a vector space basis $B$ for $K$.

Observe that if $q \in K$, then multiplication by $q$ is a $\mathbb{Q}$-linear map on $K$. Hence we can associate to $q$ a matrix $A_{K, B}(q)$ corresponding to this linear map with respect to the basis $B$. For example, if $K = \mathbb{Q}(\sqrt{2})$ and $B = B_0 := \{1, \sqrt{2}\}$, then we can express a typical element $q$ of $K$ with respect to the basis $B$ as
\[q = q_1 + q_2 \sqrt{2}.\]
If we take the product with another number $r = r_1 + r_2 \sqrt{2}$, the product is
\[(q_1 r_1 + 2 q_2 r_2) + (q_1 r_2 + q_2 r_1) \sqrt{2},\]
and so the matrix $A_{\mathbb{Q}(\sqrt{2}), B_0}(q)$ is 
\begin{equation}\label{qsqrt2matrix}
A_{\mathbb{Q}(\sqrt{2}), B_0}(q) = \left(\begin{array}{cc}
q_1 & 2 q_2 \\
q_2 & q_1
\end{array}\right); \qquad q_1, q_2 \in \mathbb{Q}.
\end{equation}
We will write $K_B$ for the image of $K$ under the map $A_{K,B}$. Observe that $A_{K,B}$ is a field homomorphism, so $K_B$ is isomorphic to $K$ as a field.

Given a vector $\mathbf{q} = (q_1, \ldots, q_k) \in \mathbb{Q}^k$, we define the scalar $\mathbf{q}_B \in K$ by $\mathbf{q}_B := q_1 \omega_1 + \cdots + q_k \omega_k$. Then the map $\textbf{q} \mapsto A_{K,B}(\mathbf{q}_B)$ can be extended by continuity to a map $A_{K,B}^* : \mathbb{R}^k \to M_{k \times k}(\mathbb{R})$. For the example $K = \mathbb{Q}(\sqrt{2})$ and $B = B_0$ as above, we have 
\begin{equation}\label{rqsqrt2matrix}
A_{\mathbb{Q}(\sqrt{2}), B_0}^*(x_1, x_2) := \left(\begin{array}{cc}
x_1 & 2 x_2 \\
x_2 & x_1
\end{array}\right); \qquad x_1, x_2 \in \mathbb{R}.
\end{equation}
The image of $\mathbb{R}^k$ under the map $A_{K,B}^*$ will be denoted $(\mathbb{R} \otimes_{\mathbb{Q}} K)_B$. The algebra $(\mathbb{R} \otimes_{\mathbb{Q}} K)_B$ is isomorphic as a $\mathbb{Q}$-algebra to the tensor product $\mathbb{R} \otimes_{\mathbb{Q}} K$. We will also use $A_{K,B}^*$ to denote the natural extension of $A_{K,B}^*$ to $\mathbb{C}^k$.

There are $k$ embeddings of the field $K$ into $\mathbb{C}$; we will use $\Sigma$ to denote the set of embeddings. Among these embeddings $\sigma \in \Sigma$, some $k_1 \leq k$ of them will be \textbf{real embeddings}; that is, the image of $\sigma$ will be contained in $\mathbb{R}$; the remaining \textbf{strictly complex embeddings} will come in $k_2$ conjugate pairs, where $k = k_1 + 2k_2$. We will write $\Rho$ for the set of real embeddings and $\Tau$ for the set of strictly complex embeddings. Generally, we will use the letter $\sigma$ to refer to an arbitrary element of $\Sigma$, we will use the letter $\rho$ to refer to an arbitrary element of $\Rho$. We will write $\tilde k$ for $k_1 + k_2$.

 We can lift these embeddings to embeddings of $K_B$ into $\mathbb{C}$; we will use $\sigma$ to denote both an embedding of $K$ into $\mathbb{C}$ and for the lift to $K_B$. For example, for the case $K = \mathbb{Q}(\sqrt{2})$, $B = B_0$, we have the pair of embeddings
\[\sigma_1 \left( \begin{array}{cc}
q_1 & 2 q_2 \\
q_2 & q_1
\end{array} \right) = q_1 + q_2 \sqrt{2}; \quad \sigma_2 \left( \begin{array}{cc}
q_1 & 2 q_2 \\
q_2 & q_1 
\end{array}\right) = q_1 - q_2 \sqrt{2}.\]
Note that there are no strictly complex embeddings; a field such as $\mathbb{Q}(\sqrt{2})$ that admits only real embeddings is called a \textbf{totally real field}. These embeddings can be extended by continuity to ring morphisms on all of $(\mathbb{R} \otimes_{\mathbb{Q}} K)_B$. We will use $\sigma^*$ to denote these extensions. For the example $K = \mathbb{Q}(\sqrt{2})$, $B = B_0$, we have
\[\sigma_1^* \left(
\begin{array}{cc}
x_1 & 2 x_2  \\
x_2 & x_1 \end{array} \right) = x_1 + x_2 \sqrt{2}; \quad \sigma_2^* \left( \begin{array}{cc}
x_1 & 2 x_2 \\
x_2 & x_1 \end{array} \right) = x_1 - x_2 \sqrt{2}.\]
We will use $\sigma^{**}$ to denote the map $\sigma^* \circ A_{K,B}^* : \mathbb{R}^k \to \mathbb{C}$, suppressing the dependence of this composition on $K$ and $B$. 

Note that since $x_1$ and $x_2$ can be any real numbers, the ring homomorphisms $\sigma_1^*$ and $\sigma_2^*$ are not embeddings; that is, they are not injective. In our example, $\sigma_1^{**}(x_1, x_2) = 0$ whenever $x_1 = - \sqrt{2} x_2$. In fact, the maps $\sigma^*$ will only be embeddings if $K$ is an imaginary quadratic field, in which case $\mathbb{R} \otimes_{\mathbb{Q}} K$ is isomorphic to $\mathbb{C}$.

The \textit{field trace} of an element $q \in K$ is the sum
\[\text{tr}_{K / \mathbb{Q}}(q) = \sum_{j=1}^k \sigma_j(q).\]
This agrees with the trace of the matrix $A_{K,B}(q)$. Observe that this trace is independent of the choice of basis $B$. We extend the field trace by continuity to all of $(\mathbb{R} \otimes_{\mathbb{Q}} K)_B$, using $\tr^*$ to denote this extension. Of course, the map $\tr^*$ is simply the matrix trace. However, it is also important to view the trace $\tr^*$ as the sum of the maps $\sigma_j^*$. 

To summarize, we have described a map $A_{K,B}^*$ from $\mathbb{R}^k$ into the ring $(\mathbb{R} \otimes_{\mathbb{Q}} K)_B$. We have a collection of ring homomorphisms $\sigma^*$ from the ring $(\mathbb{R} \otimes_{\mathbb{Q}} K)_B$ into $\mathbb{C}$. We also have an extension $\tr^*$ of the field trace sending $(\mathbb{R} \otimes_{\mathbb{Q}} K)_B$ into $\mathbb{R}$. By composing with $A_{K,B}^*$, this gives a family of $k$ functions $\sigma^{**} : \mathbb{R}^k \to \mathbb{C}$. Among the functions $\sigma$, $k_1 \leq k$ of them will be real-valued; these will be denoted $\rho$, the remaining $2 k_2$ of them will occur in complex conjugate pairs denoted by $\tau$ and $\overline{\tau}$. We will write $\tilde k = k_1 + k_2$.
\subsection{Oscillatory integrals associated to algebraic number fields}
Fix an algebraic number field $K$ of degree $k$ and a vector space basis $B$ for $K$ as a $\mathbb{Q}$-vector space. Let $f$ be a polynomial of degree $d$ in $n$ variables over the ring $(\mathbb{R} \otimes_{\mathbb{Q}} K)_B$. 

Throughout this section, the bold $\bm{\vec x}$ will denote an element of $\mathbb{R}^{kn}$. We will view $\bm{\vec x}$ as $\{x_{l,j}\}_{1 \leq l \leq n, 1 \leq j \leq k}$, and write $\bm{x_l}$ for the vector $\{x_{l, j}\}_{1 \leq j \leq k}$. 

Define the function $f^*(\bm{\vec x})$ by
\[f^*(\bm{\vec x}) =  f(A_{K,B}^*(\bm{x_1}), \ldots, A_{K,B}^*(\bm{x_n})).\] Hence $f^*$ is a $(\mathbb{R} \otimes_{\mathbb{Q}} K)_B$-valued function on $\mathbb{R}^{kn}$. We define the function $\phi_f: \mathbb{R}^{kn} \to \mathbb{R}$ by 
\[\phi_f(\bm{\vec x}) = \tr^*(f^*(\bm{\vec x})).\]
So $\phi_f$ is a real polynomial of $kn$ real variables.

Fix a smooth, compactly supported function $\psi : \mathbb{R}^{kn} \to \mathbb{C}$. We define the integral $I_{K, B, \psi}(f)$ to be
\begin{equation}\label{eq:Ipsi}
I_{K, B, \psi}(f) := \int_{\mathbb{R}^{kn}} e^{2 \pi i \phi_f(x)} \psi(x) \, dx
\end{equation}

We also associate to $f$ a family of complex polynomials $P_{f,\sigma}$ defined by applying the ring homomorphism $\sigma^*$ to each coefficient of $f$. The polynomials $P_{f,\sigma}$ have the key property that 
\begin{equation}\label{eq:pfjdef}
\sigma^{*}(f^*(\bm{\vec x})) = P_{f,\sigma}(\sigma^{**}(\bm{x_1}), \ldots, \sigma^{**}(\bm{x_n})).
\end{equation}

Given a vector $\bm{\vec x} \in \mathbb{R}^{kn}$ and an embedding $\sigma : K \to \mathbb{C}$, we will write $\vec{\sigma}(\bm{\vec x})$ for the vector
\[\vec{\sigma}(\bm{\vec x}) = (\sigma^{**}(\bm{x_1}), \ldots, \sigma^{**}(\bm{x_n}))\]

We will use the polynomials $P_{f,\sigma}$ to associate an $H$-functional to the polynomial $f$.

\begin{mydef}\label{def:galhfunc}
Let $K$ be a number field of degree $k$, and let $B$ be a vector-space basis for $K$ over $\mathbb{Q}$. Let $j \in \{1, \ldots, k\}$ and let $\sigma_j$ be one of the complex embeddings of $K$. Let $f$ be a polynomial over the ring $(\mathbb{R} \otimes_{\mathbb{Q}} K)_B$. We define the \textbf{pointwise $H$-functional associated to $\sigma$} on $\mathbb{R}^{kn}$ by
\begin{equation}\label{eq:singlehompointwise}
H_{f,K,B,\sigma}( \bm{\vec x}) = \max_{|\alpha| \geq 1} \left|\frac{1}{\alpha!} \partial^{\alpha} P_{f,\sigma}(\sigma^{**}(\bm{x_1}), \ldots, \sigma^{**}(\bm{x_n})) \right|^{1/|\alpha|}.
\end{equation}
Let $\psi$ be a smooth, compactly supported function on $\mathbb{R}^{nk}$. The \textbf{uniform $H$-functional associated to $\sigma$} is defined by
\begin{equation}\label{eq:singlehomunif}
H_{f,K,B,\sigma,\psi} = \inf_{\bm{\vec x} \in \supp \psi} H_{f,K,B,\sigma}(\bm{\vec x}).
\end{equation}
Let $S \subset \Sigma$ be a nonempty set with the property that if some strictly complex bedding $\tau_j \in S$, then its complex conjugate $\overline{ \tau_j} \in S$. The \textbf{uniform $H$-functional associated to the collection $S$} is defined by
\begin{equation}\label{eq:multiplehomunif}
H_{f,K,B,S,\psi}(\bm{\vec x}) = \prod_{\sigma \in S} H_{f,K,B,\sigma,\psi}.
\end{equation}
\end{mydef}

Observe that if $K = \mathbb{R}$ and $B = \{1\}$, then the $H$-functional reduces to the classical definition of Arhipov, Čubarikov, and Karacuba. If $K = \mathbb{Q}(i)$ and $B = \{1, i\}$, then the maps $\sigma_1^{**}$ and $\sigma_2^{**}$ are given by $\sigma_1^{**}(x_1, x_2) = x_1 + x_2 i$ and $\sigma_2^{**}(x_1, x_2) = x_1 - x_2 i$. Hence $P_{f,\sigma_2} (\sigma_2^{**}(x_1), \ldots, \sigma_2^{**}(x_n))$ is the complex conjugate of $P_{f,\sigma_1}(\sigma_1^{**}(x_1), \ldots, \sigma_1^{**}(x_n))$, so they have the same absolute value and the product is equal to $H_{f, \mathbb{C}}$. Hence Definition \ref{def:galhfunc} generalizes Definitions \ref{def:hrdef} and \ref{def:hcdef}.

More generally, if $\tau$ is any strictly complex embedding, then $H_{f,K,B,\tau}(\bm{ \vec x}) = H_{f,K,B,\overline{\tau}}(\bm{\vec x})$. This is clear since $P_{f, \tau}(\tau^{**}(x_1), \ldots, \tau^{**}(x_n)) = \tau^{**}(f(x_1, \ldots, x_n))$ and $P_{f, \overline{\tau}} (\overline{\tau^{**}}(x_1), \ldots, \overline{\tau^{**}}(x_N)) = \overline{\tau^{**}}(f(x_1, \ldots, x_n))$, so the quantities occurring in the maximum on the right side in \eqref{eq:singlehompointwise} are the same for $\tau$ and $\overline{\tau}$.

We point out some new features in Definition \ref{def:galhfunc}. First, it is important to consider many of the complex embeddings of the polynomial $f$ simultaneously. Hence, the behavior of the $H$-functional is not captured by the derivatives of a single function at a single point. 

Second, it is possible for the $H$-functional $H_{f,K,B,\sigma}(\bm{\vec x})$ to vanish for all $x \in \supp \psi$ even for nonzero polynomials $f$. For example, let $K = \mathbb{Q}(\sqrt{2})$ and $B = \{1, \sqrt{2}\}$, and $f$ be the polynomial $Ax^2$, where $A \in (\mathbb{R} \otimes_{\mathbb{Q}} K)$ is given by the matrix
\[\left( \begin{array}{cc}
-\sqrt{2} & 2 \\
1 & - \sqrt{2} \end{array} \right).\] 
The polynomial $P_{f,\sigma_1}(z)$ is equal to $0$ everywhere because $A$ maps to $0$ under the ring homomorphism $\sigma_1^*$. So $H_{f,K,B,\{\sigma_1, \sigma_2\}, \psi}$ vanishes for any function $\psi$ containing $0$ in its support.

This is the reason that it is useful to consider a \textit{proper subset} of the complex embeddings rather than consider \textit{all} of the complex embeddings.

We are now ready to state our main result.
\begin{mythm}\label{thm:mainthm}
Let $K$ be a number field of degree $k$, and let $B$ be a vector space basis for $K$ over $\mathbb{Q}$. Let $\psi : \mathbb{R}^{kn} \to \mathbb{C}$ be a smooth, compactly supported function. Then if $f$ is any nonzero polynomial in $n$ variables of degree at most $d$ over the ring $(\mathbb{R} \otimes_{\mathbb{Q}} K)_B$ and $S \subset \{\sigma_1, \ldots, \sigma_k\}$ has the property that $\tau \in S$ if and only if $\bar \tau \in S$, we have
\begin{equation}\label{eq:maineq}
I_{K,B,\psi}(f) \lesssim_{K,B,n, d, \psi} H_{f, K, B, S, \psi}^{-1},
\end{equation}
where the right side is interpreted as $\infty$ if $H_{f, K, B, S, \psi} = 0$. 
\end{mythm}
\subsection{Notation}
Throughout this work, normal--font variables such as $x$ or $z$ will be used to denote scalars in $\mathbb{R}$, $\mathbb{C}$, an algebraic number field $K$, or in the tensor product $\mathbb{R} \otimes_{\mathbb{Q}} K$. The parameter $k$ will be used exclusively to refer to the degree of the field extension $K/\mathbb{Q}$. Bold variables such as $\bm{x}$ will be used to refer to real or complex vectors whose components are indexed either by a subset of $\{1, 2, \ldots, k\}$, or indexed by a subset of the field embeddings $\Sigma$. A quantity with a vector symbol, such as $\vec x$, will always denote a vector whose components are indexed by a subset of $\{1, \ldots, n\}$. We use this notation even for multi-indices such as $\vec \alpha$ or $\vec \beta$. Finally, a bold variable with a vector symbol such as $\bm{\vec x}$ will denote a vector with components indexed by both a subset of $\Sigma$ (or a subset of $\{1, \ldots, k\}$) and a subset of $\{1, \ldots, n\}$. This same convention is used even for multi-indices.

We will use the notation $A \lesssim B$ to indicate that there is a constant $C$ such that $A \leq CB$. We write $B \lesssim A$ to mean the same thing as $A \lesssim B$. An expression such as $A \lesssim_{K, B, n, d, \psi} B$ means that $A \leq C_{K, B, n, d, \psi} B$ for a constant $C_{K,B,n,d,\psi}$ that may depend on $K, B, n, d,$ and $\psi$.

If $v$ is a vector and $P$ is a function, we write $\nabla_v P$ to denote the quantity $\nabla P \cdot v$, even if $v$ is not a unit vector.

We use the notation $e(x)$ to refer to the usual additive character $e^{2 \pi i x}$ on $\mathbb{R}$.
\section{Outline of proof}
The proof follows a similar outline to that of Wright \cite{Wright20}. First, we must relate the $k$ complex polynomials $P_{f,\sigma_j}$ to the real polynomial $\phi_f$ in the $n = 1$ case. We will see that the polynomials $P_{f,\sigma_j}$ will dictate the directional derivatives of $\phi_f$ in certain conveniently chosen directions. This is captured in the following result.
\begin{mylem}\label{lem:deriv}
Let $f$ be a polynomial of a single variable with coefficients in $(\mathbb{R} \otimes_{\mathbb{Q}} K)_B$. Then $\nabla \phi_f(x_1, \ldots, x_k) = 0$ if and only if $P_{f,\sigma}'(\sigma^{**}(x_1, \ldots, x_k)) = 0$ for all $1 \leq j \leq k$. Moreover, there exist constants $C_{K,B} > 0$ not depending on $f$ such that 
\begin{equation}\label{eq:deriv}
c_{K,B} \max_{1 \leq j \leq k} \left| P_{f,\sigma_j}'(\sigma_j^{**}(x_1, \ldots, x_k)) \right| \leq \left|\nabla \phi_f(x_1, \ldots, x_k) \right| \leq C_{K,B} \max_{1 \leq j \leq k} \left|P_{f,\sigma_j}'(\sigma_j^{**}(x_1, \ldots, x_k))\right|
\end{equation}
\end{mylem}
Lemma \ref{lem:deriv} allows for nonstationary phase estimates to be applied on the set in which the $H$-functional is dominated by a first derivative. In order to estimate the integral $I_{K,B, \psi}(f)$ in terms of the $H$-functional, it is also important to estimate the measure of the set of those $x$ for which  the $H$-functional $H_{f,K,B,\psi,\sigma}$ is dominated by a higher derivative. The key for this is a \textit{local sublevel set estimate}; that is, an estimate on the set of points where a function $g$ is small but has some large derivative. In fact, an estimate of this type is provided by Wright \cite{Wright20}.
\begin{mylem}[Structural sublevel set estimate(\cite{Wright20}, Proposition 9.1)]\label{lem:sublevel}
Let $Q \in \mathbb{C}[x]$ have degree $d$ and suppose $z_0 \in \mathbb{C}$ is such that $|z_0| \leq 1$. Suppose $|Q^{(k)}(z_0)| \geq \mu$ for some $k \geq 1$ and $|Q(z_0)| \leq \epsilon$ for some $0 < \epsilon \leq \epsilon_d$ is a sufficiently small positive constant, depending only on the degree of $Q$. Then there exists a zero $z_*$ of $Q^{(j)}$ for some $0 \leq j \leq d$ such that $|z_0 - z_*| \leq \left(\frac{\epsilon}{\mu} \right)^{1/k}$.
\end{mylem}
This structural sublevel set estimate allows us to obtain an upper bound on the measure of the set on which the $H$-functional is dominated by a higher--order derivative. 

\section{A decomposition of $\mathbb{C}^k$ and of $\mathbb{R}^k$}
In working with polynomials $f$ with coefficients in $(K \otimes_{\mathbb{R}} \mathbb{Q})_B$, it is convenient to decompose $\mathbb{R}^k$ into a direct sum of subspaces that are fixed by most of the maps $\sigma^{**}$.

The basis $\{\omega_1, \ldots, \omega_k\}$ of the field extension $K/\mathbb{Q}$ has nonzero discriminant. This discriminant is given by

\begin{equation}\label{eq:discriminant}
\disc_{K/\mathbb{Q}}(B) := \left| \begin{array}{cccc} 
\sigma_1(\omega_1) & \sigma_2(\omega_1) & \cdots & \sigma_k(\omega_1) \\ 
\sigma_1(\omega_2) & \sigma_2(\omega_2) & \cdots & \sigma_k(\omega_2) \\
\vdots & \vdots & \ddots & \vdots \\
\sigma_1(\omega_k) & \sigma_2(\omega_2) & \cdots & \sigma_k(\omega_k) \end{array} \right| \neq 0.
\end{equation}

For an embedding $\sigma : K \to \mathbb{C}$, we will use the notation $w_{\sigma}$ to denote the vector $(\sigma(\omega_1), \ldots, \sigma(\omega_k))$. Equation \eqref{eq:discriminant} implies that the vectors $w_{\sigma}$ are a basis of $\mathbb{C}^k$. Recall that we use $\sigma$ to denote an arbitrary embedding, $\rho$ to denote a real embedding, and $\tau$ to denote a strictly complex embedding.

\begin{mydef}\label{def:Wspacescplx}
For any embedding $\sigma$, we will denote the complex span of the vector $w_{\sigma}$ by  $W_{\sigma, \mathbb{C}}$. Notice that the vector spaces $W_{\sigma, \mathbb{C}}$ are \textbf{not} mutually orthogonal. For any $\sigma_0 \in \Sigma,$ we will let $V_{\sigma_0, \mathbb{C}}$ denote the one-dimensional complex vector space
\[V_{\sigma_0, \mathbb{C}} := \left( \bigoplus_{\sigma \in \Sigma, \sigma \neq \sigma_0}  W_{\sigma, \mathbb{C}} \right)^{\perp}.\]
\end{mydef}

We will also need a real version of this decomposition, but for complex embeddings $\tau$ we will need to consider the embeddings $\tau$ and $\overline{\tau}$ together.

\begin{mydef}\label{tildesigma}
Let $\sim$ be the equivalence relation defined on $\Sigma$ by $\sigma \sim \sigma'$ if $\sigma = \sigma'$ or $\sigma = \overline{\sigma}'$. We will use $\tilde \Sigma$ for the quotient of $\Sigma$ by this equivalence relation. So $\tilde \Sigma$ is a $\tilde k$-element set. If $\sigma \in \Sigma$, we write $\tilde \sigma$ for the equivalence class of $\sigma$.
\end{mydef}
\begin{mydef}\label{def:Wspaces}
For a real embedding $\rho : K \to \mathbb{R}$, we write $W_{\rho, \mathbb{R}}$ for the one-dimensional real vector space spanned by $w_{\rho}$. For a strictly complex embedding $\tau : K \to \mathbb{C}$, we write $W_{\tau, \mathbb{R}}$ for the two-dimensional real vector space spanned by the real vectors $w_{\tau} + w_{\overline{\tau}}$ and $\frac{1}{i} (w_{\tau} - w_{\overline{\tau}})$. Observe that $W_{\overline \tau, \mathbb{R}} = W_{\tau, \mathbb{R}}$ for each $\tau$, so it makes sense to talk about $W_{\tilde \sigma, \mathbb{R}}$ for each $\tilde \sigma \in \tilde \Sigma.$. Note that the vector spaces $W_{\tilde \sigma,\mathbb{R}}$ are \textbf{not} mutually orthogonal.
\end{mydef}

\begin{mydef}\label{def:vjdef}
For $\sigma_0 \in \tilde \Sigma$, define $V_{\tilde \sigma_0, \mathbb{R}}$ to be the orthogonal complement of the sum of $W_{\tilde \sigma,\mathbb{R}} : \tilde \sigma \neq \tilde \sigma_0$; that is,
\[V_{\sigma_0,\mathbb{R}} = \left( \bigoplus_{\tilde \sigma \in \tilde \Sigma, \tilde \sigma \tilde \neq \sigma_0} W_{\tilde \sigma,\mathbb{R}} \right)^{\perp}.\]
Then we have the decomposition of $\mathbb{R}^k$ into $\tilde k$ subspaces
\[\mathbb{R}^k = \bigoplus_{\tilde \sigma \in \tilde \Sigma} V_{\tilde \sigma,\mathbb{R}}.\]

Each subspace $V_{\tilde \rho,\mathbb{R}}$ has dimension $1$ and each subspace $V_{\tilde \tau,\mathbb{R}}$ has dimension $2$. Moreover, if $w \in W_{\tilde \sigma,\mathbb{R}}$ and $v \in V_{\tilde \sigma',\mathbb{R}}$ where $\tilde \sigma' \neq \tilde \sigma$, then $w \cdot v = 0$.
\end{mydef}

The advantage to working with these vector spaces $V_{\tilde \sigma}$ is that, for any vector $x \in \mathbb{R}^k$,  $\sigma^*(x + v) = \sigma^*(x)$ whenever $v \in V_{\sigma'}$ for $\tilde \sigma' \neq \tilde \sigma$.

For any real embedding $\rho$, we have that $V_{\tilde \rho, \mathbb{R}} = \mathbb{R}^k \cap V_{\rho, \mathbb{C}}$. For any strictly complex embedding $\tau$, we have $V_{\tilde \tau, \mathbb{R}} = \mathbb{R}^k \cap (V_{\tau, \mathbb{C}} + V_{\overline{\tau}, \mathbb{C}})$.

More generally, we will need to work with a basis of $\mathbb{R}^{kn}$. This is done by first writing $\mathbb{R}^{kn} = \bigoplus_{l = 1}^{n} \mathbb{R}^k$, and then decomposing each copy of $\mathbb{R}^k$ as in Definition \ref{def:vjdef}.

\begin{mydef}\label{def:vljdef}
We can decompose $\mathbb{C}^{kn}$ into $nk$ subspaces $W_{l, \sigma, \mathbb{C}}$ as follows: a vector $\bm{\vec x} = \{x_{l,j}\}_{1 \leq l \leq n, 1 \leq j \leq k}$ belongs to $W_{l_0,\sigma_0}$ if $x_{l,j} = 0$ whenever $l \neq l_0$ and the vector $\bm{x}_{l_0} := \{x_{l_0,j}\}_{1 \leq j \leq k}$ belongs to the vector space $W_{\sigma_0, \mathbb{C}}$. Similarly, we say that $\bm{\vec x} = \{x_{l,j}\}_{1 \leq l \leq n, 1 \leq j \leq k}$ belongs to $V_{l_0,\sigma_0, \mathbb{C}}$ if $x_{l,j} = 0$ whenever $l \neq l_0$ and the vector $\bm{x}_{l_0} := \{x_{l_0,j}\}_{1 \leq j \leq k}$ belongs to the vector space $V_{\sigma_0, \mathbb{C}}$. Observe that these vector spaces satisfy the orthogonality property that $V_{l_1, \sigma,\mathbb{C}} \perp W_{l_2, \sigma',\mathbb{C}}$ if $l_1 \neq l_2$ or if $l_1 = l_2$ and $\sigma' \neq \sigma$.

We can decompose $\mathbb{R}^{kn}$ into $n \tilde k$ subspaces $W_{l, \tilde \sigma, \mathbb{R}}$ as follows: a vector $\bm{\vec x} = \{x_{l,j}\}_{1 \leq l \leq n, 1 \leq j \leq k}$ belongs to $W_{l_0,\tilde \sigma_0}$ if $x_{l,j} = 0$ whenever $l \neq l_0$ and the vector $\bm{x}_{l_0} := \{x_{l_0,j}\}_{1 \leq j \leq k}$ belongs to the vector space $W_{\tilde \sigma_0, \mathbb{R}}$. Similarly, we say that $\bm{\vec x} = \{x_{l,j}\}_{1 \leq l \leq n, 1 \leq j \leq k}$ belongs to $V_{l_0,\tilde \sigma_0, \mathbb{R}}$ if $x_{l,j} = 0$ whenever $l \neq l_0$ and the vector $\bm{x}_{l_0} := \{x_{l_0,j}\}_{1 \leq j \leq k}$ belongs to the vector space $V_{\tilde \sigma_0, \mathbb{R}}$. Observe that these vector spaces satisfy the orthogonality property that $V_{l_1, \tilde \sigma,\mathbb{R}} \perp W_{l_2, \tilde \sigma',\mathbb{R}}$ if $l_1 \neq l_2$ or if $l_1 = l_2$ and $\tilde \sigma' \neq \tilde \sigma$.
\end{mydef}

The reason for constructing this basis is that for $v \in V_{l', \tilde \sigma', \mathbb{R}}$, we have $\sigma^{**}((\bm{ \vec v})_{l}) = 0$ if $(l,\tilde \sigma) \neq (l',\tilde \sigma')$, and $\sigma'^{**}((\bm{\vec v})_{l'}) \neq 0$ if $(l,\tilde \sigma) = (l',\tilde \sigma')$. So this decomposition of vector spaces allows us to ``nudge" exactly one of the inputs to each of the polynomials $P_{f,\sigma}$ and $P_{f, \overline{\sigma}}$ in the expression $P_{f,\sigma}(\vec \sigma (\bm{\vec x}))$ without changing the value of $\vec \sigma'(\bm{\vec x})$ for $\tilde \sigma' \neq \tilde \sigma$.

Finally, we define the spaces $V_{\tilde{\sigma}, \mathbb{R}}^n \subset \mathbb{R}^{kn}$ by
\begin{equation}\label{eq:vtildesigmarstar}
V_{\tilde \sigma, \mathbb{R}}^n := \bigoplus_{1 \leq l \leq n} V_{l, \tilde \sigma, \mathbb{R}}.
\end{equation}
Notice that $V_{\tilde \sigma, \mathbb{R}}^n$ is $n$-dimensional if $\sigma = \rho \in \Rho$ is a real embedding and $2n$-dimensional if $\sigma = \tau \in \Tau$ is a complex embedding.

We will need a technical lemma on the spaces $V_{\tilde \sigma, \mathbb{R}}^n$.

\begin{mylem}\label{lem:realpart}
Let $\sigma_0 \in \Tau$ be a strictly complex embedding. Suppose that $\vec x$ is a nonzero vector in $\mathbb{C}^n$. Then there exists a unit vector $\bm{\vec z} \in V_{\tilde \sigma_0, \mathbb{R}}^n$ such that 
\[\left| \realpart (\vec x \cdot \vec \sigma_0(\bm{\vec z})) \right| \gtrsim_{K,B} |\vec x|.\]
\end{mylem}
\begin{myrmk}
The implicit constant in Lemma \ref{lem:realpart} does not depend on $n$, but this fact will not be used anywhere.
\end{myrmk}
\begin{proof}

First we show this in the case $n = 1$. In this case, the vector $\vec x$ is simply a nonzero complex number $x$. If $\sigma_0$ is a strictly complex embedding, then the vector space $V_{\tilde \sigma_0, \mathbb{R}}$ is spanned by two real vectors $\bm{v_{\tilde \sigma_0, 1}}$ and $\bm{v_{\tilde \sigma_0, 2}}$. Observe that, since the vectors $\{\bm{w_{\sigma}}\}_{\sigma \in \Sigma}$ form a basis of $\mathbb{C}^n$, it follows that the $2$-by-$k$  matrix $M$ whose rows are given by $\bm{w_{\sigma_0}}$ and $\bm{w_{\bar \sigma_0}}$ has rank equal to $2$. Moreover, we know that, for $\bm{x} \in \bigoplus_{\tilde \sigma \neq \tilde \sigma_0} V_{\sigma, \mathbb{R}}$ that $\bm{w_{\sigma_0}} \cdot \bm{x}$ and $\bm{w_{\bar \sigma_0}} \cdot \bm{x}$ are both equal to zero. Since the range of the matrix $M$ must be two-dimensional, it follows that the $2$-dimensional vectors $(\sigma_0(\bm{v_{\tilde \sigma_0, 1}}), \sigma_0(\bm{v_{\tilde \sigma_0,2}}))$ and $(\bar \sigma_0(\bm{v_{\tilde \sigma_0, 1}}), \bar \sigma_0(\bm{v_{\tilde \sigma_0,2}}))$ are linearly independent over $\mathbb{C}$. But the second of these vectors is simply the complex conjugate of the first one; if $\sigma_0(\bm{v_{\tilde \sigma_0, 1}})$ and $\sigma_0(\bm{v_{\tilde \sigma_0, 2}})$ had the same (or opposite) complex argument, then so would $\bar \sigma_0(\bm{v_{\tilde \sigma_0, 1}})$ and $\bar \sigma_0 (\bm{v_{\tilde \sigma_0, 2}})$, and the vector $(\bar \sigma_0(\bm{v_{\tilde \sigma_0, 1}}), \bar \sigma_0(\bm{v_{\tilde \sigma_0,2}}))$ could be obtained from the vector $(\sigma_0(\bm{v_{\tilde \sigma_0, 1}}), \sigma_0(\bm{v_{\tilde \sigma_0,2}}))$ by multiplying by a unit scalar. It follows that $\sigma_0(\bm{v_{\tilde \sigma_0, 1}})$ and $\sigma_0(\bm{v_{\tilde \sigma_0, 2}})$ do not have the same (or opposite) complex argument, so the product of one of these scalars with $x$ must have nontrivial real part.

The $n$-dimensional case is similar as the above reasoning can be applied to each component of $\vec x$. Suppose the component $x_l$ of $\vec x$ is nonzero. Then the argument of the previous paragraph guarantees the existence of a vector $\bm{\vec z_l} \in V_{l, \tilde \sigma, \mathbb{R}}$ such that $|\realpart x_l \vec \sigma(\bm{\vec z_l}) | \gtrsim |x_l|$. Choosing such a vector $\bm{\vec z_l}$ for each $l$ such that $x_l$ is nonzero gives the desired vector $\bm{\vec z}$.
\end{proof}
\begin{myrmk}
Notice that the above lemma might fail if $\sigma$ is a real embedding; for example, in the $n = 1$ case, if $x$ is a purely imaginary number, then $\realpart(x  \{\vec \sigma(\bm{z})\}_l) = 0$ for any $\bm{z} = \sum_{l=1}^n \bm{\vec z}$.
\end{myrmk}
\section{Proof of Lemma \ref{lem:deriv}}

\begin{proof}[Proof of Lemma \ref{lem:deriv}]
The proof of Lemma \ref{lem:deriv} makes use of the following identity: 
\begin{equation}\label{eq:trsum}
\tr^* x = \sum_{j=1}^k \sigma_j^*(x) \quad \forall x \in (\mathbb{R} \otimes_{\mathbb{Q}} K)_B.
\end{equation}
Let $f$ be a polynomial of a single variable with coefficients in $(\mathbb{R} \otimes_{\mathbb{Q}} K )_B$ and let $\bm{x} = (x_1, \ldots, x_k) \in \mathbb{R}^k$. Then we have
\begin{IEEEeqnarray*}{Cl}
& \phi_f(x_1, \ldots, x_k) \\
= & \tr f(A^*(x_1, \ldots, x_k)) \\
= & \sum_{\sigma \in \Sigma} \sigma^*(f(A^*(x_1, \ldots, x_k))) \\
= & \sum_{\sigma \in \Sigma} P_{f,\sigma}(\sigma^*(A^*(x_1, \ldots, x_k))) \\
= & \sum_{\sigma \in \Sigma} P_{f,\sigma}(x_1 \sigma_j(\omega_1) + \cdots + x_k \sigma_j(\omega_j)) \\
= & \sum_{\sigma \in \Sigma} P_{f,\sigma} (\bm{x} \cdot \bm{w_{\sigma}}).
\end{IEEEeqnarray*}
where $\bm{w_{\sigma}} \in \mathbb{C}^k$ is the nonzero vector $(\sigma_j(\omega_1), \cdots, \sigma_j(\omega_j))$ and $\cdot$ denotes the usual dot product.

Hence, 
\[\nabla \phi_f(x_1, \ldots, x_k) = \sum_{\sigma \in \Sigma}^k P_{f,\sigma}'(\bm{x} \cdot \bm{w_{\sigma}}) \bm{w_{\sigma}}.\]
If $\nabla \phi_f(x_1, \ldots, x_k)$ is nonzero, then certainly one of the $P_{f,\sigma}$ is nonzero, and the triangle inequality guarantees that 
\[|P_{f,\sigma}'(\bm{x} \cdot \bm{w_{\sigma}})| \geq \frac{1}{k \norm{w_j} } \left|\nabla \phi_f(x_1, \ldots, x_k) \right|.\]

Conversely, assume $P_{f,\sigma}'(\bm{x} \cdot \bm{w_{\sigma_0}}) \neq 0$ for some $\sigma_0 \in \Sigma$. In order to show that $\nabla \phi_f$ is nonzero, it is enough to show that $\nabla \phi_f \cdot \bm{v}$ is nonzero for some nonzero vector $\bm{v}$. Indeed, let $\bm{v}$ be any complex unit vector in $V_{\sigma_0, \mathbb{C}}$. Then $\bm{v} \cdot \bm{w_{\sigma_0}} \neq 0$, and we have
\begin{IEEEeqnarray*}{Cl}
& \nabla \phi_f(x_1, \ldots, x_k) \cdot  \bm{v} \\
= & \left( \sum_{\sigma \in \Sigma} P_{f,\sigma}'(\bm{x} \cdot \bm{w_{\sigma}}) \right) \bm{w_{\sigma}} \cdot \bm{v} \\
 = & P_{f,\sigma_0}'(\bm{x} \cdot \bm{w_{\sigma_0}} ) \bm{w_{\sigma_0}} \cdot \bm{v} \neq 0.
\end{IEEEeqnarray*}
Hence $\nabla \phi_f(x_1, \ldots, x_k) \neq 0$ and 
\[| \nabla \phi(f_1, \ldots, x_k)| \geq |\nabla \phi_f(x_1, \ldots, x_k)| \cdot \bm{v}| \geq |P_{f,\sigma_0}'(x \cdot \bm{w_{{\sigma_0}}}) \bm{w_{\sigma_0}} \cdot \bm{v}|.\]
\end{proof}

Lemma \ref{lem:deriv} essentially reduces the study of the size of $\nabla \phi_f$ to the study of the size of the first derivatives of $P_{f,\sigma}$. 
\section{A Structural sublevel set estimate for algebraic number fields}
Wright \cite{Wright20} establishes Lemma \ref{lem:sublevel} using a Hensel--type argument. However, we will not need to repeat this argument as Wright's result can be applied directly to prove the analogue for algebraic number fields.

With a little more work, we can obtain a bound on the size of the sublevel set for polynomials of \textit{more than one variable}. Recall that given a vector $v$ and a function $P$, we will use the notation $\nabla_v P$ to denote the quantity $\nabla P \cdot v$, even if $v$ is not a unit vector.

First, we will need a lemma of Wright.
\begin{mylem}[Wright \cite{Wright20}]\label{lem:directional}
Let $r \geq 1$. Then there exists a finite collection of vectors $\mathcal{U}_r$ such that the differential operators
\[\left\{(\nabla_{u})^r : u \in \mathcal{U}_r \right\}\]
form a basis for the space of constant-coefficient differential operators of order exactly $r$. In particular, any derivative $\partial^{\alpha}$ with $|\alpha| = r$ can be written uniquely as a linear combination of the operators $(\nabla_u)^r$.
\end{mylem}
Using this lemma, we are now ready to establish the following corollary.
\begin{mycor}\label{cor:sublevelkreal}
Let $f$ be a polynomial of $n$ variables of degree at most $d$ with coefficients in $(\mathbb{R} \otimes_{\mathbb{Q}} K)_B$. Let $S \subset \Sigma$ be a collection of homomorphisms such that if $\sigma \in S$, then $\overline{\sigma} \in S$. For $\sigma \in S$, let $0 < \epsilon_{\sigma} < \mu_{\sigma}$ be real numbers, writing $\bm \epsilon$ for the $|S|$-tuple $\{\epsilon_{\sigma}\}_{\sigma \in S}$ and $\bm{\mu}$ for the $|S|$-tuple $\{\mu_{\sigma}\}_{\sigma \in S}$. For each $\sigma \in S$, Let $\vec \alpha_{\sigma}$ be a multi-index with $n$ entries. Write $\bm{\vec \alpha}$ for the $(|S|n)$-tuple $\{\vec \alpha_{\sigma}\}_{\sigma \in S}$. Let $\mathbf{B}$ denote the unit cube in $\mathbb{R}^{nk}$. Define $A_{S, \bm{\vec \alpha}, \bm{\epsilon}, \bm{\mu}}$ to be the set
\begin{equation}\label{eq:Asepsmudef}
A_{S, \bm{\vec \alpha}, \bm{\epsilon}, \bm{\mu} } := \left\{\bm{\vec x} \in \mathbf{B} : |\nabla P_{f,\sigma}(\vec \sigma(\bm{\vec x}))| \leq \epsilon_{\sigma}, \left| \partial^{\vec \alpha_{\sigma}} P_{f,\sigma}(\vec \sigma(\bm{ \vec x})) \right| \geq \mu_{\sigma} \; \forall \sigma \in S \right\}. \
\end{equation}
Then the Lebesgue measure of $A_{S, \bm{\vec \alpha}, \bm{\epsilon}, \bm{\mu}}$ satisfies the bound
\begin{equation}\label{eq:Asepsmuest}
|A_{S, \bm{\vec \alpha}, \bm{\epsilon}, \bm{\mu}}| \lesssim_{K,B,d,n} \prod_{\sigma \in S} \left(\frac{\epsilon_{\sigma}}{\mu_{\sigma}} \right)^{1/(|\alpha_{\sigma}| - 1)}.
\end{equation}
\end{mycor}
\begin{myrmk}
Even though $S$ is closed under complex conjugation, the product is extended over $\sigma \in S$, not $\tilde \sigma \in \tilde S$. So if $\sigma$ is a strictly complex embedding, then both the $\sigma$ and $\bar \sigma$ factors will appear in the product.
\end{myrmk}
\begin{proof}[Proof of Corollary \ref{cor:sublevelkreal}]
For a vector $\bm{\vec x} \in \mathbb{R}^{kn}$ written as $\bm{\vec x} = \{x_{l,j}\}_{\substack{l \in \{1, \ldots, n\} \\ j \in \{1, \ldots, k\}}}$, we write $\bm{x_l}$ for the $k$-dimensional vector $\{x_{l,j}\}_{j \in \{1, \ldots, n\}}$. 

Recall the decompositions
\[\mathbb{C}^{kn} = \bigoplus_{\sigma \in \Sigma, 1 \leq l \leq n} W_{l, \sigma, \mathbb{C}}\] 
and
\[\mathbb{C}^{kn} = \bigoplus_{\sigma \in \Sigma, 1 \leq l \leq n} V_{l, \sigma, \mathbb{C}}\] 
given in Definition \ref{def:vljdef}. As each space $W_{l, \sigma, \mathbb{C}}$ and $V_{l, \sigma, \mathbb{C}}$ is a one-dimensional complex vector space, we can select vectors $\bm{ \vec w_{l, \sigma, \mathbb{C}}} \in W_{l, \sigma, \mathbb{C}}$ and $\bm{\vec v_{l, \sigma, \mathbb{C}}} \in V_{l, \sigma, \mathbb{C}}$ that span their respective spaces.

Observe that for $l \neq l'$, we have the orthogonality relations
\begin{equation}\label{eq:vjlorthogonality}
V_{l,\sigma,\mathbb{C}} \perp V_{l', \sigma', \mathbb{C}} ; W_{l,\sigma,\mathbb{C}} \perp W_{l', \sigma', \mathbb{C}}.
\end{equation}

Let $\bm{\vec x} \in A_{S, \bm{\vec \alpha}, \bm{\epsilon}, \bm{\mu}}$. By assumption, we have for each $\sigma \in S$ that there exists some partial derivative $\partial^{\vec \alpha_{\sigma}}$ such that $\left| \partial^{\vec \alpha_{\sigma}} P_{f,\sigma}(\vec \sigma(\bm{\vec x}))\right| \gtrsim \mu_{\sigma}$. Hence, Lemma \ref{lem:directional} implies that for each $\sigma \in S$ there exists a unit vector $\vec u_{\sigma} \in \mathcal{U}_{|\alpha_{\sigma}|}$ such that $(\nabla_{u_{\sigma}})^{|\alpha_{\sigma}|} P_{f,\sigma}(\bm{\vec x}) \gtrsim \mu_{\sigma}$.

For each $\sigma$, we define the $kn$-dimensional vector
\[\bm{\vec u_{\sigma}}^* := \sum_{l=1}^n \frac{(\vec u_{\sigma})_l}{\bm {\vec w_{l,\sigma, \mathbb{C}}} \cdot \bm{\vec v_{l,\sigma, \mathbb{C}}}} \bm{\vec v_{l,\sigma, \mathbb{C}}},\]
where $(\vec u_{\sigma})_l$ denotes the $l$th component of $\vec u_{\sigma}$.
A simple calculation involving \eqref{eq:vjlorthogonality} shows that for any $\bm{\vec x} \in \mathbb{C}^{nk}$
\[\nabla_{\bm{\vec u}_{\sigma}^*} \left[ P_{f, \sigma} (\vec \sigma(\bm{\vec x})) \right] = (\nabla_{u_{\sigma}} P_{f,\sigma}) \left(\vec \sigma(\bm{\vec x}) \right).\]
Moreover, the vectors $\{\bm{ \vec u}_{\sigma}^*\}_{\sigma \in S}$ are linearly independent since each $\bm{\vec u}_{\sigma}^*$ is a nonzero linear combination of the $\mathbf{v}_{l,\sigma}^*$ for fixed $\sigma$, and the vectors $\mathbf{v}_{l,\sigma}^*$ form a basis of $\mathbb{C}^{nk}$.  Extend the collection of vectors $\mathbf{u}_{\sigma}^*$ to a basis $\{\bm{\vec u}_j^*\}_{j=1}^{kn}$ of $\mathbb{C}^{nk}$, with the understanding that $\bm{ \vec u}_{1}^*, \ldots, \bm{\vec u}_{|S|}^*$ are an enumeration of the vectors $\bm{\vec u}_{\sigma}^*$. If $\bm{\vec u}_j^* = \bm{\vec u_\sigma}^*$, write $\sigma = \sigma_j$, $\alpha_j = \alpha_{\sigma}$, $\epsilon_j = \epsilon_{\sigma}$, $\mu_j = \mu_{\sigma}$, and $P_{f,j} = P_{f, \sigma}$.

Write $\mathcal{Z}_S(f)$ for the set
\[\mathcal{Z}_S(f) = \left\{x \in \mathbb{R}^{kn} : (\nabla_{u_{\sigma}})^{\beta_{\sigma}} P_{f,\sigma}(\vec \sigma(\bm{\vec x})) = 0 \text{ for some $1 \leq \beta_{\sigma} \leq d$ and all $\sigma \in S$.} \right\}\]
For a complex vector $\bm{\vec d} = (d_{|S|+1}, \ldots, d_{nk})$, write $\mathcal{Z}_{S, \bm{\vec d}}(f)$ for the cross-section
\[\mathcal{Z}_{S,\bm{\vec d}}(f) = \left\{x \in \mathcal{Z}_S(f) : x = \sum_{j=1}^{|S|} c_j \bm{\vec u_j}^* + \sum_{j=|S|+1}^{nk} d_j \bm{\vec u_j}^*\right\}.\]
We claim that for any such $\bm{\vec d}$, the set $\mathcal{Z}_{S,\bm{\vec d}}$ consists of only finitely many points. Indeed, for a fixed $\bm{\vec d}$ and for an embedding $\sigma_j \in S$, the polynomial
\[P_{f,{\sigma_{j_0}}}\left(\vec \sigma_{j_0} \left( \sum_{j=1}^{|S|} c_j \bm{\vec u_j}^* + \sum_{j=|S|+1}^{nk} d_j \bm{\vec u_j}^* \right) \right)\]
does not depend on the variables $c_j$ for $j \neq j_0$; hence, there are only finitely many values of $c_{j_0}$ where one of the polynomials
\[(\nabla_{u_0})^{\beta}P_{f,{\sigma_{j_0}}}\left(\vec \sigma_{j_0} \left( \sum_{j=1}^{|S|} c_j \bm{\vec u_j}^* + \sum_{j=|S|+1}^{nk} d_j \bm{\vec u_j}^* \right) \right)\]
vanishes for some $\beta \leq d$. Call this collection $\mathcal{Z}_{S, \bm{\vec d}, \sigma_{j_0}}$. Observe then that
\[\mathcal{Z}_{S,\bm{\vec d}}(f) = \left\{\sum_{j=1}^{|S|} c_j \bm{\vec u_j}^* + \sum_{j=|S| + 1}^{nk} d_j \bm{\vec u_j}^* : c_j \in \mathcal{Z}_{S, d, \sigma_j} \text{for all $j$} \right\},\]
a finite set.

Now we show that for any point $\bm{\vec x} = \sum_{j=1}^{|S|} c_j \mathbf{u}_j^* + \sum_{j=|S| + 1}^{nk} d_j \mathbf{u}_{j}^* \in A_{S, \bm{\vec \alpha}, \bm{\epsilon}, \bm{\mu}}$, there exists a point $\bm{\vec y}  =  \sum_{j = 1}^{|S|} c_j' \mathbf{u}_j^*  + \sum_{j = |S| + 1}^{kn} d_j \mathbf{u}_{j^*} \in \mathcal{Z}_{S, \bm{\vec d}}(f)$ such that $|c_j - c_j'| \lesssim  \left(\frac{\epsilon_j}{\mu_j} \right)^{1/(|\alpha_j| - 1)}$ for all $j \in \{1, \ldots, |S|\}$. By applying Lemma \ref{lem:sublevel}, we observe that there exists a complex number $c_{1}'$ with $|c_{1}' - c_{1}| \lesssim \left(\frac{\epsilon_1}{\mu_1} \right)^{1/(|\alpha_1| - 1)}$ and an order $\beta_{1}$ such that $(\nabla_{u_{\sigma_1}})^{\beta_1} P_{f, 1}\left(\vec \sigma_1 (\mathbf{x_1}) \right) = 0$ where $\bm{\vec x}_1 = \bm{\vec x} + (c_{1}' - c_{1}) \mathbf{u}_{1}^*$. Observe that $\vec \sigma \left(\bm{\vec x}_1 \right) =\vec \sigma \left(\bm{\vec x} \right)$ for $\sigma \neq \sigma_1$.

 We continue in this manner---assuming we have defined $\bm{\vec x}_i$, we define $\bm{\vec x}_{i+1}$ by the equation $\bm{\vec x}_{i + 1} = \bm{\vec x}_i + (c_{i + 1}' - c_{i+1}) \bm{\vec u }_{i+1}^*$, where $c_{{i + 1}'}$ is chosen so that $(\nabla_{u_{\sigma_{i+1}}})^{\beta_{i+1}} P_{f, i+1}(\vec \sigma_{i+1} (\bm{\vec x}_{j+1})) = 0$. Observe that for $\sigma \neq \sigma_{i+1}$, we have that $\vec \sigma(\bm{\vec x_i}) =\vec \sigma(\bm{\vec x_{i + 1}})$. In particular, this means that for $i' \leq i$, we have that $(\nabla_{u_{\sigma_{i'}}})^{\beta_{i'}} P_{f, i'}\left(\vec \sigma_{i'}( \bm{\vec x}_{i + 1} ) \right)$ remains zero.

We continue this procedure until we reach $\bm{\vec x}_{|S|} =: \bm{\vec y}$. Then $P_{f, i}^{(\vec \beta_{i})} \left(\vec \sigma_{i}(\mathbf(x)) \right) = 0$ for all $i$, and $\left|c_{i} - c_{i}'\right| \lesssim \left(\frac{\epsilon_{i}}{\mu_{i}} \right)^{1/(|\alpha_{i}| - 1)}$. Hence, for any fixed choice of $\{d_i\}_{i \notin \{1, \ldots, |S|\}}$, we have that the coefficients $\{c_i\}_{1 \leq i \leq  |S|}$ must lie within a box in $\mathbb{C}^{n,k}$ with side lengths $\lesssim \left(\frac{\epsilon_i}{\mu_i} \right)^{1/(|\alpha_i|-1)}$ centered at one of the points in the finite set $\mathbf{Z}_{S,\bm{\vec d}}(f)$. Since the linear transformation sending the basis $\mathbf{u}_j^*$ to the standard basis $\mathbf{e}_j$ maps the unit box to a parallelepiped with bounded eccentricity, it follows by taking the intersection of this box with $\mathbb{R}^{nk}$ that the Lebesgue measure of the set $A_{S, \bm{\vec \alpha}, \bm{\epsilon}, \bm{\mu}}$ satisfies the inequality \eqref{eq:Asepsmuest}.
\end{proof}
\section{The $J$-functional}
Following Wright \cite{Wright20}, we define the $J$-functional to be a quantity dictating the size of the higher derivatives.

For this section, we will let $S \subset \Sigma$ be a set with the property that $\sigma \in S$ if and only if $\overline \sigma \in S$. It then makes sense to consider the reduction $\tilde S$ of $S$ under complex conjugation. We will say $\tilde \sigma \in \tilde S$ if $\sigma, \overline \sigma \in S$.  
\begin{mydef}\label{def:Jfunc}
Let $K$ be an algebraic number field of degree $k$, and let $B$ be a vector space basis for $K$. Suppose $j \in \{1, \ldots, k\}$. Let $f$ be a polynomial over the ring $(\mathbb{R} \otimes_{\mathbb{Q}} K)_B$. The \textbf{single-homomorphism pointwise $J$-functional} $J_{f,K,B,\sigma}$ is defined for $\bm{\vec x} \in \mathbb{R}^{kn}$ by
\begin{equation}\label{eq:Jfuncptwise}
J_{f,K,B,\sigma}(\bm{\vec x}) = \max_{|\alpha| \geq 2} \left| \frac{1}{\alpha!} \partial^{\alpha} P_{f,\sigma}(\vec \sigma(\bm{\vec x})) \right|^{1/|\alpha|}. 
\end{equation} 
\end{mydef} 

Given a vector $\bm{\vec x} \in \mathbb{R}^{nk}$ and a subset $S \subset \Sigma$ such that $\tau \in S$ if and only if $\overline{\tau} \in S$, we define the radius $r_{S, \sigma}(\bm{\vec x})$ for $\sigma \in \Sigma$ by 
\[r_{S,\sigma}(\bm{\vec x}) = \begin{cases}
 J_{f, K, B, \sigma}(\bm{\vec x})^{-1} & \text{if $\sigma \in S$} \\
1 & \text{if $\sigma \notin S$.}
\end{cases}.\]

\begin{myrmk}\label{rmk:Jconj}
As is the case for the $H$-functional, we have for any strictly complex embedding $\tau : K \to \mathbb{C}$:
\[J_{f,K,B,\tau}(\bm{\vec x}) = J_{f, K, B, \overline{\tau}}(\bm{\vec x}).\]
Hence, since we are imposing the condition that $\tau \in S$ if and only if $\overline{\tau} \in S$:
\[r_{S, \tau}(\bm{\vec x}) = r_{S, \overline{\tau}}(\bm{\vec x}).\]
In particular, this means that it always makes sense to consider the quantities $J_{f,K,b, \tilde \sigma}$ and $r_{S, \tilde \sigma}$ for $\tilde \sigma \in \tilde \Sigma$.
\end{myrmk}

We write $\mathbf{r}(\bm{\vec x})$ for the $\tilde k$-vector $\{r_{\tilde \sigma}(\bm{\vec x})\}_{\tilde \sigma \in \tilde \Sigma}$.

Because we wish to dilate the vector $\mathbf{r}(\bm{\vec x})$ but do not care about $\tilde \sigma \notin \tilde S$, we define the dilation of $\mathbf{r}$ in the $\tilde S$-directions as follows. For a constant $C > 0$, we will define $r_{S, \tilde \sigma,C}(\bm{\vec x})$ by
\[r_{S, \tilde \sigma,C}(\bm{\vec x}) = \begin{cases}

C \cdot J_{f, K, B, \sigma}(\bm{\vec x})^{-1} & \text{if $\tilde \sigma \in \tilde S$} \\
1 & \text{if $\tilde \sigma \notin \tilde S$.}
\end{cases}.\]
We write $\mathbf{r}_{S,C}(\bm{\vec x})$ for the vector $\{r_{S, \tilde \sigma,C}(\bm{\vec x})\}_{\tilde \sigma \in \tilde \Sigma}$.
Given a scalar $C \in \mathbb{R}$, We define the polydisc $\mathbf{P}_{S,C}(\bm{\vec x})$ to be the set
\[\left\{\bm{\vec x} + \sum_{\tilde \sigma \in \tilde \Sigma} r_{S,\tilde \sigma,C}(\bm{\vec x}) \bm{u_{\tilde \sigma}} : \bm{u_{\tilde \sigma}} \in V_{\tilde \sigma, \mathbb{R}}^n; \left\|\bm{u_{\tilde \sigma}}\right\| \leq 1. \right\}.\]
The polydiscs $\mathbf{P}_{S,C}(\bm{\vec x})$ are important because the $J$-functional does not vary much over these polydiscs. 
\begin{mylem}\label{lem:Jstability}
\leavevmode
\begin{enumerate}[(i)]
\item \label{item:twosided} There exist $C_1 > 0$ and $\epsilon_{d,n} > 0$ such that if $\bm{\vec x} \in\mathbb{R}^{kn}$, $\epsilon \leq \epsilon_{d,n}$ and $\bm{\vec x}' \in \mathbf{P}_{S,\epsilon}(\bm{\vec x})$, then for each $\sigma \in S$:
\[C_1^{-1} J_{f, K, B, \sigma}(\bm{\vec x}') \leq  J_{f, K, B, \sigma}(\bm{\vec x}) \leq C_1 J_{f, K, B, \sigma}(\bm{\vec x'}) .\]
\item \label{item:onesided} For any $A > 0$, there exists a constant $C$ depending on $A,d,n,$ and $k$ such that if $\bm{\vec x}' \in \mathbf{P}_{A}(\bm{\vec x})$, then for each $\sigma \in S$:
\[J_{f,K,B,\sigma}(\bm{\vec x}') \leq C J_{f, K, B, \sigma}(\bm{\vec x}).\]
\end{enumerate}
\end{mylem}
\begin{myrmk}\label{rmk:Jstability}
Lemma \ref{lem:Jstability} can be rephrased in terms of the vectors $\mathbf{r}(\bm{\vec x})$ and $\mathbf{r}(\bm{\vec x}')$. Part \eqref{item:twosided} can be rephrased to state that for $\bm{\vec x}' \in P_{\epsilon}(x)$, we have for each $\sigma \in \Sigma$ that 
\begin{equation}\label{eq:twosidedr}
C_1^{-1} r_{S,\sigma}(\bm{\vec x}') \leq r_{S,\sigma}(\bm{\vec x}) \leq C_1 r_{S,\sigma}(\bm{\vec x}').
\end{equation}
There is no need to restrict to $\sigma \in S$, since for $\sigma \notin S$, we have that $r_{\sigma}(\bm{\vec x}) = r_{\sigma}(\bm{\vec x}') = 1$. Similarly, if the constant $C$ is defined as in part \eqref{item:onesided}, we have for $\bm{\vec x}' \in \mathbf{P}_A(\bm{\vec x})$ and any $\sigma \in \Sigma$ that
\begin{equation}\label{eq:onesidedr}
r_{S,\sigma}(\bm{\vec x}') \geq C^{-1} r_{S,\sigma}(\bm{\vec x}).
\end{equation}
\end{myrmk}
\begin{proof}
The proof of this lemma is very similar to the proof of the corresponding lemma of Wright \cite{Wright20}. We will write $a$ for a quantity that will later be specialized to be either $\epsilon$ or $A$. 

Since $\bm{\vec x}' \in \mathbf{P}_{S,a}(\bm{\vec x})$, we can find vectors $\bm{\vec u}_{\tilde \sigma}$ with $\bm{\vec u}_{\tilde \sigma} \in V_{\tilde \sigma, \mathbb{R}}^n$ and $\norm{\bm{\vec u}_{\tilde \sigma}} \leq 1$ for $\tilde \sigma \in \tilde \Sigma$ such that
\[\bm{\vec x}' = \bm{\vec x} + \sum_{\tilde \sigma \in \tilde \Sigma}^k r_{S, \tilde \sigma, a}(\bm{\vec x}) \bm{\vec u}_{\tilde \sigma}.\]
Fix $\sigma \in S$. Since $\vec \sigma'(\bf{\vec u}_{\tilde \sigma}) = 0$ if $\tilde \sigma' \neq \tilde \sigma$, we have
\[\vec \sigma(\bm{\vec x}') - \vec \sigma (\bm{\vec x}) = a r_{S, \tilde \sigma} (\bm{\vec x}) \vec \sigma(\bm{\vec u}_{\tilde \sigma})\]
Observe that $ar_{S, \sigma}(\bm{\vec x}) \vec \sigma(u_{\tilde \sigma})$ is a vector in $\mathbb{C}^n$ with norm at most $C_2 a r_{S, \sigma}(\bm{\vec x})$ for some constant $C_2(K,B,n)$ that does not depend on $\bm{\vec x}$ or $\bm{\vec x}'$. Write $\vec z$ for the $n$-dimensional vector $\vec \sigma(u_{\tilde \sigma})$. 

Choose a multi-index $\vec \alpha_0$ with $|\vec \alpha_0| \geq 2$ such that $J_{f,K,B,\sigma}(\bm{\vec x}) = \left|\frac{1}{\vec \alpha_0!} \partial^{\vec \alpha_0} P_{f,\sigma}(\bm{\vec x}) \right|^{1/{|\alpha_0|}}.$ Write $k_0$ for $|\vec \alpha_0|$. Let $\vec \alpha$ be an arbitrary multi-index  with $2 \leq |\alpha| \leq d$, and write $k$ for $|\vec \alpha|$. We will write $Q =\partial^{\vec \alpha}P$. Then, by expanding the Taylor series of $Q$ about $\bm{\vec x}$, we have
\[Q(\bm{\vec x}') = Q(\bm{\vec x}) + \sum_{1 \leq |\vec \beta| \leq d - k} \frac{1}{\vec \beta !} \partial^{\vec \beta} Q(\bm{\vec x}) (a \mathbf{r}_{S, \tilde \sigma}(\bm{\vec x}))^{|\vec \beta|} (\vec z)^{\vec \beta}.\]
Observe that $|(\vec z)^{\vec \beta}| \leq C_{K,B,d,n}$ for some constant $C_{K,B,d,n}$ that does not depend on the specific vectors $\bm{\vec x}$ and $\bm{\vec x}'$. For each $\vec \beta$, we have
\[\left|\partial^{\vec \beta} Q(\bm{\vec x}) \right| (a r_{S,\tilde \sigma}(\bm{\vec x}))^{|\vec \beta|} = \left| \partial^{(\vec \alpha + \vec \beta)} P_{f,\sigma}(\bm{\vec x}) \right| a^{|\vec \beta|} \cdot \frac{1}{\left|\frac{1}{\vec \alpha_0!} \partial^{\vec \alpha_0} P_{f,\sigma}(\bm{\vec x}) \right|^{|\vec \beta|/k_0}}.\] 
However, it is also easily seen that
\begin{IEEEeqnarray*}{rCl}
\left|\partial^{\vec \alpha + \vec \beta} P_{f,\sigma}(\bm{\vec x}) \right| & = & (
\vec \alpha + \vec \beta)! \left| \frac{1}{(\vec \alpha + \vec \beta)!}  \partial^{\vec \alpha + \vec \beta} P_{f,\sigma} (\bm{\vec x}) \right|^{\frac{|\vec \alpha| + |\vec \beta|}{|\vec \alpha| +  |\vec \beta|}} \\
& \leq & (\vec \alpha + \vec \beta)! J_{f,K,B,\sigma} (\bm{\vec x})^{k + |\vec \beta|} \\
& = & (\vec \alpha + \vec \beta)! \left| \frac{1}{\vec \alpha_0!} \partial^{\vec \alpha_0} P_{f,\sigma}(\bm{\vec x}) \right|^{\frac{k + |\vec \beta|}{k_0}}.
\end{IEEEeqnarray*}
Hence, we observe that
\[\left| \partial^{\vec \beta} Q (\bm{\vec x}) \right| (a r_{S,\tilde \sigma}(\bm{\vec x}))^{|\vec \beta|} \leq (\vec \alpha + \vec \beta)! a^{|\vec \beta|} \left| \frac{1}{\vec \alpha_0!} \partial^{\vec \alpha_0} P_{f,\sigma}(\bm{\vec x}) \right|^{k/k_0} \]
Hence
\begin{IEEEeqnarray*}{rCl}
\left| \sum_{1 \leq |\vec \beta| \leq d - k} \frac{1}{\vec \beta!} \partial^{\vec \beta} Q(\bm{\vec x}) (a r_{S,\sigma}(\bm{\vec x}))^{|\vec \beta|} z^{\vec \beta} \right| \leq C_{K,B,d,n} \left|\frac{1}{\vec \alpha_0!} \partial^{\vec \alpha_0} P_{f,\sigma}(\bm{\vec x}) \right|^{k/k_0} \sum_{1 \leq t \leq d - k} a^{t} \sum_{|\vec \beta| = t} \frac{(\vec \alpha + \vec \beta)!}{\vec \beta!}
\end{IEEEeqnarray*}
Write 
\[C_{K,B,d,n}(a) := C_{K,B,d,n} \max_{2 \leq |\vec \alpha| \leq d} \sum_{1 \leq t \leq d - |\alpha|} a^t \sum_{|\vec \beta| = t} \frac{(\vec \alpha + \vec \beta)!}{\vec \beta!}.\]
Observe that $C_{K,B,d,n}(a) \to 0$ as $a \to 0$. Then, the above calculations show that for each $\alpha$, we have
\begin{equation}\label{eq:alphadiff}
\left|\partial^{\vec \alpha} P_{f,\sigma}(\bm{\vec x}') - \partial^{\vec \alpha} P_{f,\sigma}(\bm{\vec x}) \right| \leq C_{K,B,d,n}(a) \left|\frac{1}{\vec \alpha_0!} \partial^{\vec \alpha_0} P_{f,\sigma}(\bm{\vec x}) \right|^{|\vec \alpha|/k_0} = C_{K,B,d,n}(a) J_{f,K,B,\sigma}(\bm{\vec x})^{|\vec \alpha|}.
\end{equation}

Hence we have the upper bound
\begin{IEEEeqnarray*}{rCl}
\left|\frac{1}{\vec \alpha!} \partial^{\vec \alpha} P_{f,\sigma}(\bm{\vec x}') \right| & \leq & \left|\frac{1}{\vec \alpha!} \partial^{\vec \alpha} P_{f,\sigma}(\bm{\vec x})\right| + \frac{1}{\vec \alpha!} C_{K,B,d,n}(a) J_{f,K,B,\sigma}(\bm{\vec x})^{|\vec \alpha|} \\
\left|\frac{1}{\vec \alpha!} \partial^{\vec \alpha} P_{f,\sigma}(\bm{\vec x}') \right| & \leq & J_{f,K,B,\sigma}(\bm{\vec x})^{|\vec \alpha|} + \frac{1}{\vec \alpha!} C_{K,B,d,n}(a) J_{f,K,B,\sigma}(\bm{\vec x})^{|\vec \alpha|} \\
\left|\frac{1}{\vec \alpha!} \partial^{\vec \alpha} P_{f,\sigma}(\bm{\vec x}') \right|^{1/|\vec \alpha|} & \leq & \left(1 + \frac{1}{\vec \alpha!} C_{K,B,d,n}(a)\right)^{1/|\vec \alpha|} J_{f,K,B,\sigma}(\bm{\vec x}). \IEEEyesnumber \label{eq:xprimeupper}
\end{IEEEeqnarray*}
and the lower bound
\begin{IEEEeqnarray*}{rCl}
\left|\frac{1}{\vec \alpha!} \partial^{\vec \alpha} P_{f,\sigma}(\bm{\vec x}') \right| & \geq & \left|\frac{1}{\vec \alpha!} \partial^{\vec \alpha} P_{f,\sigma}(\bm{\vec x})\right| - \frac{1}{\vec \alpha!} C_{K,B,d,n}(a) J_{f,K,B,\sigma}(\bm{\vec x})^{|\vec \alpha|} \\
\left|\frac{1}{\vec \alpha!} \partial^{\vec \alpha} P_{f,\sigma}(\bm{\vec x}') \right| & \geq & J_{f,K,B,\sigma}(\bm{\vec x})^{|\vec \alpha|} - \frac{1}{\vec \alpha!} C_{K,B,d,n}(a) J_{f,K,B,\sigma}(\bm{\vec x})^{|\vec \alpha|} \\
\left|\frac{1}{\vec \alpha!} \partial^{\vec \alpha} P_{f,\sigma}(\bm{\vec x}') \right|^{1/|\vec \alpha|} & \geq & (1 - \frac{1}{\vec \alpha!} C_{K,B,d,n}(a))^{1/|\vec \alpha|} J_{f,K,B,\sigma}(\bm{\vec x}). \IEEEyesnumber \label{eq:xprimelower}
\end{IEEEeqnarray*}

To prove \eqref{item:twosided}, choose $\epsilon_{K,B,d,n}$ so that $\ C_{K,B,d,n}(\epsilon) < \frac{1}{2}$ for all $\epsilon< \epsilon_{K,B,d,n}$. Then for any $\vec \alpha,$ we have that 
\[\frac{1}{2} \leq \left(1 - \frac{1}{\vec \alpha!} C_{K,B,d,n}(\epsilon) \right)^{1/|\vec \alpha|}  \leq \left(1 + \frac{1}{\vec \alpha!} C_{K,B,d,n}(\epsilon) \right)^{1/|\vec \alpha|} \leq \frac{3}{2}.\]
Taking the maximum over all $\vec \alpha$ in inequalities \eqref{eq:xprimeupper} and \eqref{eq:xprimelower} gives the desired result \eqref{item:twosided}.

To prove \eqref{item:onesided}, we use the upper bound \eqref{eq:xprimeupper}. Taking $a = A$, we can deduce \eqref{item:onesided} by taking the constant $C = \max_{|\vec \alpha| \geq 2} (1 + \frac{1}{\vec \alpha!} C_{K,B,d,n}(A))^{1/|\vec \alpha|}$.
\end{proof}
Fix $\epsilon < \epsilon_{K,B,d,n}$. We now prove a basic covering property of these polydiscs. Consider the polydiscs $\mathbf{P}_{S, 3 \epsilon}(\bm{\vec x})$ for $\bm{\vec x} \in \supp \psi$. As these polydiscs cover the set $\supp \psi$, the Vitali covering lemma implies that there exists a finite collection $\{\bm{\vec x}_e\}_{e \in E}$ of points such that the polydiscs $\mathbf{P}_{S, 3 \epsilon}(\bm{\vec x}_e)$ cover $\supp \psi$ while the smaller polydiscs $\mathbf{P}_{S, \epsilon}(\bm{\vec x}_e)$ are disjoint.

The next lemma states that dilations of these polydiscs have bounded overlap.
\begin{mylem}\label{lem:overlap}
For any $C > 0$, there exists an integer $N(C) \in \mathbb{N}$ such that each point $\bm{\vec x}^*$ lies in at most $N$ of the polydiscs $\mathbf{P}_{S,C}(\bm{\vec x}_e)$. Here $N$ depends only on $C, d, n, \epsilon, K$, and $B$, but does not otherwise depend on the polynomial $f$ or the point $\bm{\vec x}^*$.
\end{mylem}
\begin{proof}
For simplicity of notation, any implicit constants appearing in this argument will be allowed to depend on the parameters $C, d, n, \epsilon, K$, and $B$. Let $\bm{ \vec x_e}$ be any point such that $\bm{\vec x^*} \in \mathbf{P}_{S,C}(\bm{\vec x_e})$. For any point $\bm{\vec x}$ in the polydisc $\mathbf{P}_{S,C}(\bm{\vec x}_e)$, we have from \eqref{eq:onesidedr} of Remark \ref{rmk:Jstability} that for some constant $C'$ and for each $\tilde \sigma \in \tilde \Sigma$:
\begin{equation}\label{eq:rxelowerbound}
r_{S,\tilde \sigma}(\bm{\vec x}_e) \leq C' r_{S,\tilde \sigma}(\bm{\vec x}).
\end{equation}
Let $\{\tilde \sigma_1, \ldots, \tilde \sigma_{\tilde k}\}$ be an enumeration of the elements of $\tilde \Sigma$. Then we can write

\[\bm {\vec x_e} = \sum_{j=1}^{\tilde k} \bm{\vec x}_{\tilde \sigma_j}; \quad \bm{\vec x_{\tilde \sigma_j}} \in V_{\tilde \sigma_j, \mathbb{R}}^n\]
and we can write
\[\bm{\vec x^*} = \sum_{j=1}^{\tilde k} \left(\bm{\vec{x}_{\tilde \sigma_j}} + r_{S, \tilde \sigma_{j}, C}(\bm{ \vec x_e}) \bm{ \vec u_{\tilde \sigma_{j}}} \right); \quad \bm{\vec u_{\tilde \sigma_j}} \in V_{\tilde \sigma_j, \mathbb{R}}^n; \quad \norm{\bm{\vec u_{\tilde \sigma_j}}} \leq 1.\]
For $j_0 \in \{1, \ldots, \tilde k\}$, define $\bm{\vec{y_{j_0}}}$ to be the point 
\[\bm{\vec y_{j_0}} = \sum_{j=1}^{j_0} \left( \bm{\vec x_{\tilde \sigma_j}} + r_{S, \tilde \sigma_j, C}(\bm{\vec x_e}) \bm{\vec u_{\tilde \sigma_j}} \right)  + \sum_{j = j_0 + 1}^{\tilde k} \bm{\vec x_{\tilde \sigma_j}},\]
taking $\bm{\vec y}_0 = \bm{\vec x}_e$ and $\bm{\vec y}_{\tilde k} = \bm{\vec x}^*$. In particular, this implies that $r_{\tilde \sigma_j}(\bm{\vec y}_j) = r_{\tilde \sigma_j}(\bm{\vec x}^*)$ for each $1 \leq j \leq \tilde k$ and that $r_{\tilde \sigma_{j+1}}(\bm{\vec y}_j) = r_{\tilde \sigma_{j+1}}(\bm{\vec x}_e)$ for every $0 \leq j \leq \tilde k - 1$. We define $L_j$ and vectors $\bm{\vec z}_{j,i}$ for $1 \leq j \leq \tilde k$ and $0 \leq i \leq L_j$ by 
\[\bm{\vec z}_{j,0} = \bm{\vec y}_{j-1}\]
and such that, for $i \geq 0$, the vector $\bm{\vec z}_{j,i} = \bm{\vec y}_{j-1} + t_l (\bm{\vec y}_j - \bm{\vec y}_{j-1})$, where $1 \leq i \leq L_j$ and the scalars $t_i$ are chosen so that for each $1 \leq i \leq L_j$:
\[ \left|\bm{\vec z}_{j,i} - \bm{\vec z}_{j,i-1} \right| = \epsilon r_{\tilde \sigma_j}(\bm{\vec z}_{j, i - 1})\]
and so that
\[\left| \bm{\vec y}_j - \bm{\vec z}_{j, L_j} \right| \leq \epsilon r_{\tilde \sigma_{j+1}}(\bm{\vec z}_{j, L_j}).\]
For simplicity, we will also write $\bm{\vec z}_{j, L_j + 1} := \bm{\vec y}_j$. We can obtain an upper bound on $L_j$ using the lower bound on $r_{\tilde \sigma_j}(\bm{\vec z}_{j,i})$. Since $\bm{\vec z}_{j,i} \in \mathbf{P}_{S, C}(\bm{\vec x}_e)$ for every $j$ and $i$, it follows from \eqref{eq:rxelowerbound} that 
\[r_{\tilde \sigma_j}(\bm{\vec x}_e) \leq C' r_{\tilde \sigma_j}(\bm{\vec z}_{j,i}).\]
Hence, we conclude
\[(L_j + 1) \epsilon r_{\tilde \sigma j} (\bm{\vec x}_e) \leq C' \epsilon \sum_{i=0}^{L_j} r_{\tilde \sigma_{j}}(\bm{\vec z}_{j,i}) = C' \sum_{i=0}^{L_j} \left| \bm{\vec z}_{j,i + 1} - \bm{\vec z}_{j,i} \right| \leq C' \left|\bm{\vec y}_j - \bm{\vec y}_{j-1} \right| \leq C C' r_{\tilde \sigma_j}(\bm{\vec x}_e).\]
Thus $L_j \lesssim 1$ for each $j$.

In order to complete the argument, we need a \textit{lower} bound on each component of the vector $\mathbf{r}(\bm{\vec x}_e)$. In order to do this, we need to apply the bound \eqref{item:twosided} of Lemma \ref{lem:Jstability} successively to each $\bm{\vec z}_{j,i}$. Fix $j$. Since $\left| \bm{\vec z}_{j,i + 1} - \bm{\vec z}_{j, i} \right| \leq \epsilon r_{\tilde \sigma_j}(\bm{\vec z}_{j,i}) $ for each $i$, we have for each $0 \leq i \leq L_j$ that

\[r_{\tilde \sigma_j}(\bm{\vec z_{j,i + 1}}) \leq C_1 r_{\tilde \sigma_j}(\bm{\vec z_{j,i}}) \]
Hence, by applying this $L_{j+1}$ times, we conclude
\[r_{\tilde \sigma_j}(\bm{\vec x}^*) = r_{\tilde \sigma_j}(\bm{\vec y}_j) \leq C_1^{L_j + 1} r_{\tilde \sigma_j} (\bm{\vec y}_{j-1}) = C_1^{L_j + 1} r_{\tilde \sigma_j}(\bm{\vec x}_e).\]
Observe that $C_1^{L_j + 1} \lesssim 1$ for each $j$. Define $C_2 = \max_{1 \leq j \leq k} C_1^{L_j + 1}$. Then for each $j$, we have the inequality
\begin{equation}\label{eq:rxeupperbound}
r_{\tilde \sigma_j}(\bm{\vec x}^*) \leq C_2 r_{\tilde \sigma_j}(\bm{\vec x}_e).
\end{equation}
This means that a polydisc of the desired form that contains $\bm{\vec x}^*$ cannot be too small. Now define 
\[\mathcal{P}(\bm{\vec x}^*) = \{e : \bm{\vec x}^* \in \mathbf{P}_{S,C}(\bm{\vec x}_e)\}.\]
We then have by \eqref{eq:rxelowerbound} that
\[\bigcup_{e \in \mathcal{P}(\bm{\vec x}^*)} \mathbf{P}_{S,C} (\bm{\vec x}_e) \subset 2 \mathbf{P}_{S,CC'}(\bm{\vec x}^*)\]
Since the polydiscs on the left side of this inclusion are disjoint, we conclude from \eqref{eq:rxeupperbound} that
\[\# \mathcal{P} \cdot \left|\mathbf{P}_{S,C'}(\bm{\vec x}^*)\right| \lesssim \sum_e \left| \mathbf{P}_{S,\epsilon}(\bm{\vec x}_e) \right| \lesssim 2^{nk} \left|\mathbf{P}_{S,C'}(\bm{\vec x}^*) \right|.\]
Hence $\# \mathcal{P} \lesssim 1$ as desired.
\end{proof}

We can also obtain a good bound on the first derivatives of $P_{f,\sigma_j}$ at $\bm{\vec x}' \in \mathbf{P}_{S,\epsilon}(\bm{\vec x})$.
\begin{mylem}\label{lem:gradstability}
Suppose $\bm{\vec x}' \in \mathbf{P}_{S,\epsilon}(\bm{\vec x})$ where $\epsilon < \epsilon_{K,B,d,n}$ from Lemma \ref{lem:Jstability}. Then we have
\[|\nabla P_{f, \sigma}(\bm{\vec x}') - \nabla P_{f, \sigma}(\bm{\vec x})| \lesssim_{K,B,d,n} \epsilon J_{f,K,B,\sigma}(\bm{\vec x}').\]
\end{mylem}
The implicit constant does not depend on the specific points $\bm{\vec x}$ and $\bm{\vec x}'$ or on the polynomial $f$. 
\begin{proof}[Proof of Claim]
We will prove this using part \eqref{item:twosided} of Lemma \ref{lem:Jstability}. Observe that for any $l$, we have
\begin{IEEEeqnarray*}{Cl}
& \left| \partial_l P_{f,\sigma}(\bm{\vec x}') - \partial_l P_{f,\sigma}(\bm{\vec x}) \right|\\
= & \left| \int_{\bm{\vec x}}^{\bm{\vec x}'} \nabla \partial_l P_{f,\sigma}(\bm{\vec y}) \cdot \frac{\bm{\vec x}' - \bm{\vec x}}{\norm{\bm{\vec x}' - \bm{\vec x}}} \, d \bm{\vec y} \right| \\
\lesssim_{K,B,d,n} & \int_{\bm{\vec x}}^{\bm{\vec x}'} J_{f,K,B,\sigma}(\bm{\vec y})^2 \, d \bm{\vec y}.\\
\lesssim_{K,B,d,n} & J_{f,K,B,\sigma}(\bm{\vec x})^2 \cdot \epsilon r_{\sigma}(\bm{\vec x})  \\
= & \epsilon J_{f,K,B,\sigma}(\bm{\vec x}),
\end{IEEEeqnarray*}
as desired.
\end{proof}
\section{Proof of the main theorem}\label{sec:mainthm}
We are now ready to prove Theorem \ref{thm:mainthm}. The proof will follow the proof of Lemma 11.2 of Wright \cite{Wright20}. 
\begin{proof}
Let $\epsilon > 0$ be smaller than $\frac{\epsilon_{K,B,d,n}}{6}$, where $\epsilon_{K,B,d,n}$ is as in Lemma \ref{lem:Jstability}. Choose a smooth, nonnegative function $\gamma : \mathbb{R}^{k} \to \mathbb{R}_{\geq 0}$ with $\gamma(x) = 1$ for $|x| \leq 1$ and $\gamma(x) = 0$ for $|x| \geq 2$. Define the function
\[\Gamma(\bm{\vec x}) = \sum_e \gamma_e(\bm{\vec x})\]
where
\[\gamma_e(\bm{\vec x}) = \prod_{\tilde \sigma \in \tilde \Sigma} \gamma((3 \epsilon r_{\tilde \sigma} (\bm{\vec x}_e) )^{-1} (\bm{\vec x_{\tilde \sigma}} - \bm{\vec x}_{e, \tilde \sigma})). \]
Here, the vector $\bm{\vec x}$ is written as 
\begin{equation}\label{eq:xdecomp}
\bm{\vec x} = \sum_{\tilde \sigma \in \tilde \Sigma} \bm{\vec x_{\tilde \sigma}} : \bm{\vec x_{\tilde \sigma}} \in V_{\tilde \sigma, \mathbb{R}}^n,
\end{equation}
and the vector $\bm{\vec x_e}$ is written as
\[\bm{\vec x_e} = \sum_{\tilde \sigma \in \tilde \Sigma} \bm{\vec x_{e, \tilde \sigma}} : \bm{\vec x_{e, \tilde \sigma}} \in V_{\tilde \sigma, \mathbb{R}}^n.\]
For later notational convenience, we also introduce the function 
\[\gamma_0(\bm{\vec y}) = \prod_{\sigma} \gamma(\bm{\vec y_{\sigma}}),\]
where $\bm{\vec y}$ is decomposed into $\bm{\vec y_{\tilde \sigma}}$ as in equation \eqref{eq:xdecomp}.

The function $\Gamma$ is bounded by the constant $N(6 \epsilon)$ from Lemma \ref{lem:overlap}. Since the polydiscs $\mathbf{P}_{S,3 \epsilon}(\bm{\vec x_e})$ cover $\supp \psi$, the function $\Gamma$ satisfies $\Gamma \geq 1$ on $\supp \psi$. Moreover,
\[\supp \Gamma \subset \bigcup_{e = 1}^M \mathbf{P}_{S,6 \epsilon} (\bm{\vec x}_e).\]
Write $I$ for the integral $I_{K,B, \psi}(f)$. We decompose $I$ as 
\begin{equation}\label{eq:idecomp}
I = \sum_{e \in E} \int_{\mathbb{R}^{k n}} e(\phi_f(\bm{\vec x})) \psi_e(\bm{\vec x}) d\bm{\vec x},
\end{equation}
where the function $\psi_e$ is defined by
\[\psi_e(\bm{\vec x}) := \frac{1}{\Gamma(\bm{\vec x})} \psi(\bm{\vec x}) \gamma_e(\bm{\vec x}). \]
Write $I_e$ for the summand in \eqref{eq:idecomp}. 

We will partition the set $E$ according to various data related to the values $H_{f,K,B,\sigma}(\bm{\vec x}_e)$ and $J_{f,K,B,\sigma}(\bm{\vec x}_e)$ for $\sigma \in S$.

First, we define the $n|S|$-multi-index $\bm{\vec \alpha}(e) = \{\alpha_{\sigma, l}(e)\}_{\sigma \in S, 1 \leq l \leq n}$ so that the $n$-multi-index $\vec \alpha_{\sigma}(e)$ is defined so that
\[J_{f, K, B, \sigma}(\bm{\vec x_e}) = \frac{1}{\vec \alpha_{\sigma}!} \partial^{\vec \alpha_{\sigma}} P_{f, \sigma} (\vec \sigma(\bm{\vec x_e})) .\]
If more than one value of $\vec \alpha_{\sigma}$ achieves the maximum in \eqref{eq:Jfuncptwise}, we define $\vec \alpha_{\sigma}$ to be the first one in the lexicographic order. Note that $\vec \alpha_{\sigma} = \vec \alpha_{\bar \sigma}$ for any $\sigma \in \Sigma$, so it makes sense to talk about $\vec \alpha_{\tilde \sigma}$ as well.

For an $n |S|$-multi-index $\bm{\vec \alpha}$, define
\[T_{\bm{\vec \alpha}} = \{e \in T : \bm{\vec \alpha}(e) = \bm{\vec{\alpha}}\}.\]

Fix a large constant $C$ and a large number $N$. Suppose $S_1 \cup S_2 \cup S_3$ is a partition of $S$ into three disjoint sets closed under taking complex conjugates.  Then we define the set $E_{\bm{\vec \alpha}, S_1, S_2, S_3}$ to be the set of those $e \in E_{\bm{\vec \alpha}}$ such that
\begin{IEEEeqnarray*}{rClL}
J_{f,K,B,\sigma}(\bm{\vec x}_e) & \geq & C^{-1} |\nabla P_{f, \sigma}(\vec \sigma(\bm{\vec x}_e))| & \text{for $\sigma \in S_1$} \\ 
1  \leq J_{f,K,B,\sigma}(\bm{\vec x}_e) & \leq & C^{-1}| \nabla P_{f, \sigma}(\vec \sigma(\bm{\vec x}_e))| & \text{for $\sigma \in S_2$} \\
J_{f, K, B, \sigma} (\bm{\vec x}_e) \leq 1 &  \leq & C^{-1} |\nabla P_{f, \sigma}(\vec \sigma(\bm{\vec x}_e))| & \text{for $\sigma \in S_3$}.
\end{IEEEeqnarray*}
For $\sigma \in S_1$ we will need additional information about the size of $J_{f, K, B, \sigma}(\bm{\vec x_e})$ relative to the uniform $H$-functional $H_{f,K,B,\sigma, \psi}$ ; for $\sigma \in S_2$ we will need information about the size of $\nabla P_{f, \sigma}(\vec \sigma(\bm{\vec x}_e))$ relative to both $H_{f, K, B, \sigma, \psi}$ and to $J_{f, K, B, \sigma}(\bm{\vec x}_e)$. In order to keep track of this information, we will need additional data in the form of the vectors $\bm{r}_1 = \{r_{1, \sigma}\}_{\sigma \in S_1}$, $\bm{r_2} = \{r_{2, \sigma}\}_{\sigma \in S_2}$, and $\bm{l} = \{l_\sigma\}_{\sigma \in S_2}$.

We define $E_{\bm{\vec \alpha}, S_1, S_2, S_3,\bm{r_1}, \bm{r_2}, \bm{l}}$ to be the set of those $e \in E_{\bm{\vec \alpha}, S_1, S_2, S_3}$ such that
\begin{IEEEeqnarray*}{rClL}
J_{\sigma}(\bm{\vec x_e}) & \sim & C^{-1} 2^{r_{1, \sigma}} H_{f,K, B, \psi, \sigma} \quad & \text{for $\sigma \in S_1$} \\
|\nabla P_{\sigma}(\bm{\vec x_e})| & \sim & C^{-1} 2^{r_{2,\sigma}} H_{f, K, B, \psi, \sigma} \quad & \text{for $\sigma \in S_2$.} \\
|\nabla P_{\sigma}(\bm{\vec x_e})| & \sim & C 2^{l_{\sigma}} J_{\sigma}(\bm{\vec x_e}) \quad & \text{for $\sigma \in S_2$.} \\
\end{IEEEeqnarray*}
Observe that the definition of $E_{\bm{\vec \alpha}, S_1, S_2, S_3}$ guarantees that each component of $\bm{r_1}$, $\bm{r_2}$, and $\bm{l}$ is nonnegative.

We will estimate the contribution of $I_{e}$ for $e \in T_{\bm{\vec \alpha}, S_1, S_2, S_3, \bm{r_1}, \bm{r_2}, \bm{l}}$. Recall that $I_e$ is defined by
\[I_e = \int_{\bm{\vec x} \in \mathbb{R}^{kn}} \frac{1}{\Gamma}(\bm{\vec x}) \psi(\bm{\vec x}) \gamma_{e} (\bm{\vec x}) e( \phi_f(\bm{\vec x})).\]
We expand this integral, making the change of variable $\bm{\vec x_{\tilde \sigma}} = \bm{\vec x_{e, \tilde \sigma}} + r_{\tilde \sigma}(\bm{\vec x_e}) \bm{ \vec y_{\tilde \sigma}}$ for each $\tilde \sigma \in \tilde \Sigma$. In making this change-of-variables we accrue a constant Jacobian factor of $|\mathbf{P}_{S, 1} (\bm{\vec x_e})|$. Hence $I_e$ is given by
\begin{IEEEeqnarray*}{Cl}
& |\mathbf{P}_{S,1} (\bm{\vec x_e}) | \int_{\bm{\vec y}} \gamma_e \left(\bm{\vec x_e} + \sum_{j=1}^{\tilde k} r_{\tilde \sigma_j} (\bm{\vec x_e}) \bm{\vec y_{\tilde \sigma_j}} \right) \frac{\psi}{\Gamma} \left(\bm{\vec x_e} + \sum_{j=1}^{\tilde k}  r_{\tilde \sigma_j}(\bm{\vec x_e}) \bm{\vec y_{\tilde \sigma_j}} \right)e \left(\phi_f\left(\bm{\vec x_e} + \sum_{j=1}^{\tilde k} r_{\tilde \sigma_j}(\bm{\vec x_e}) \bm{\vec y_{\tilde \sigma_j}} \right) \right) d \bm{\vec y} . \\
= & |\mathbf{P}_{S,1} (\bm{\vec x_e}) | \int_{\bm{\vec y}} \gamma_0 \left((3 \epsilon)^{-1} \bm{\vec y}\right) \frac{\psi}{\Gamma} \left(\bm{\vec x_e} + \sum_{j=1}^{\tilde k}  r_{\tilde \sigma_j}(\bm{\vec x_e}) \bm{\vec y_{\tilde \sigma_j}} \right) \prod_{\sigma \in \Sigma} e \left(P_{f, \sigma} (\vec \sigma (\bm{\vec x_e} + r_{\tilde \sigma}(\bm{\vec x_e}) \bm{\vec y_{\tilde \sigma}})) \right) \, d \bm{\vec y}. \\
\end{IEEEeqnarray*}
By combining the $\sigma$ and $\bar \sigma$ factors in the product, we have the equation
\[I_e = |\mathbf{P}_{S,1} (\bm{\vec x_e}) | \int_{\bm{\vec y}} \gamma_0 \left((3 \epsilon)^{-1} \bm{\vec y}\right) \frac{\psi}{\Gamma} \left(\bm{\vec x_e} + \sum_{j=1}^{\tilde k}  r_{\tilde \sigma_j}(\bm{\vec x_e}) \bm{\vec y_{\tilde \sigma_j}} \right) \prod_{\tilde \sigma \in \tilde \Sigma} e \left(C_{\tilde \sigma} \realpart(P_{f, \sigma} (\vec \sigma (\bm{\vec x_e} + r_{\tilde \sigma}(\bm{\vec x_e}) \bm{\vec y_{\tilde \sigma}}))) \right) \, d \bm{\vec y}, \]
where $C_{\tilde \sigma} = 1$ if $\sigma$ is a real embedding (in which case the expression inside of $\realpart$ is already real) and $C_{\tilde \sigma} = 2$ if $\sigma$ is a strictly complex embedding.

Using integration by parts in the form of Corollary \ref{cor:intbyparts}, we have the estimate
\begin{IEEEeqnarray*}{rCl}
|I_e| & \lesssim_{\psi} & |\bm{P_{S,1}}(\bm{\vec x_e})| \prod_{\tilde \sigma \in \tilde S_2 \cup \tilde S_3} |\nabla P_{f, \sigma}(\bm{\vec x_e})|^{-N} \frac{(1 + J_{f, K, B, \sigma}(\bm{\vec x_e})^{-N})}{J_{f, K, B, \sigma}(\bm{\vec x_e})^{-N}} \\
& \lesssim_{\psi} &  |\bm{P_{S,1}}(\bm{\vec x_e})| \cdot \left( \prod_{\tilde \sigma \in \tilde S_2} \frac{J_{f, K, B, \sigma}(\bm{\vec x_e})^N}{|\nabla P_{f, \sigma}(\bm{ \vec x_e})|^N} \right) \cdot \left(\prod_{\tilde \sigma \in \tilde S_3}  \left|\nabla P_{f, \tilde \sigma}(\bm{\vec x_e}) \right|^{-N} \right).
\end{IEEEeqnarray*}
Because $e$ belongs to $E_{\bm{\vec \alpha}, S_1, S_2, S_3, \bm{r_1}, \bm{r_2}, \bm{l}}$, we have for $\sigma \in S_2$ that $\frac{J_{f, K, B, \sigma}(\bm{\vec x_e})^N}{|\nabla P_{f \sigma}(\bm{\vec x_e})|^N} \lesssim 2^{l_{\tilde \sigma} N}$. Hence the total contribution of $|I_e| : e \in E_{\bm{\vec \alpha}, S_1, S_2, S_3, \bm{r_1}, \bm{r_2}, \bm{l}}$ is 

\[\sum_{e \in E_{\bm{\vec \alpha}, S_1, S_2, S_3, \bm{r_1}, \bm{r_2}, \bm{l}}} |I_e| \lesssim_{C, N \psi} \sum_{e \in T_{\bm{\vec \alpha}, S_1, S_2, S_3, \bm{r_1}, \bm{r_2}, \bm{l}}} |\mathbf{P}_{S,1}(\bm{\vec x_e})| \prod_{\tilde \sigma \in \tilde S_2} 2^{-(l_{2, \tilde \sigma})N} \prod_{\tilde \sigma \in S_3} H_{\tilde \sigma}^{-N}.\]

Now, for an appropriate constant $C_{K,B,d,n,\epsilon}$, we define the set\footnote{Note that these sets are not disjoint.} $\mathfrak{S}_{\bm{\vec \alpha}, S_1, S_2, S_3, \bm{r_1}, \bm{r_2}, \bm{l}}$ to be the set of those $x \in \supp \psi$ such that $J_{f, K, B, \sigma}(\bm{\vec x}) \sim |\partial^{\vec \alpha_{\sigma}} P_{f, \sigma}(\bm{\vec x})|^{1/|\vec \alpha_{\sigma}|}$ for all $\sigma \in S_1 \cup S_2$ and 
\begin{IEEEeqnarray}{rCll}
J_{f, K, B, \sigma}(\bm{\vec x}) & \geq & C_{K,B,d,n,\epsilon}^{-1} |\nabla P_{f, \sigma}(\vec \sigma(\bm{\vec x}))| \quad & \text{for $\sigma \in S_1$} \label{eq:S1eqn} \\
C_{K,B,d,n,\epsilon}^{-1} \leq J_{f,K,B,\sigma} (\bm{\vec x}) & \leq & C_{K,B,d,n,\epsilon} |\nabla P_{f, \sigma} (\vec \sigma(\bm{\vec x}))| \quad & \text{for $\sigma \in S_2$} \label{eq:S2eqn}\\
J_{f, K, B, \sigma}(\bm{\vec x}) & \leq & C_{K,B,D,n, \epsilon}^{-1} \quad & \text{for $\sigma \in S_3$}\label{eq:S3eqn} \\
J_{f,K,B, \sigma}(\bm{\vec x}) & \sim & 2^{r_{1, \sigma}} H_{f, K, B,  \psi, \sigma} \quad & \text{for $\sigma \in S_1$}\label{eq:S1r} \\
|\nabla P_{f,\sigma}(\bm{\vec x})| & \sim & 2^{r_{2, \sigma}} H_{f, K, B, \psi, \sigma} \quad & \text{for $\sigma \in S_2$}\label{eq:S2r} \\
|\nabla P_{f, \sigma}(\bm{\vec x})| & \sim &  2^{l_{\sigma}} J_{f, K, B, \sigma}(\bm{x}) \quad & \text{for  $\sigma \in S_2$}\label{eq:S2l}.
\end{IEEEeqnarray}
Suppose $C$ is sufficiently large that $C \gg C_{K,B,d,n}$, where $C_{K,B,d,n}$ is the implicit constant from Lemma \ref{lem:gradstability}. For $\bm{\vec x} \in \bm{P}_{S,6 \epsilon}(\bm{\vec x_e})$, we can deduce that \eqref{eq:S1eqn} and \eqref{eq:S2eqn} hold by Lemma \ref{lem:gradstability}. We can deduce that \eqref{eq:S3eqn} and \eqref{eq:S1r} hold by applying Lemma \ref{lem:Jstability}. We deduce that \eqref{eq:S2r} holds by applying Lemma \ref{lem:gradstability} and that \eqref{eq:S2l} holds by applying both Lemmas \ref{lem:Jstability} and \ref{lem:gradstability}. Hence $\bm{P}_{S,6 \epsilon}(\bm{\vec x_e}) \subset \mathfrak{S}_{\bm{\vec \alpha}, S_1, S_2, S_3, \bm{r_1}, \bm{r_2}, \bm{l}}$. It remains to estimate the size of the set $\mathfrak{S}_{\bm{\vec \alpha}, S_1, S_2, S_3, \bm{r_1}, \bm{r_2}, \bm{l}}$. This size will be estimated using Corollary \ref{cor:sublevelkreal} with $S$ replaced by $S_1 \cup S_2$. For $\sigma \in S_1$, we have from \eqref{eq:S1eqn} and \eqref{eq:S1r} that the estimates in equation \eqref{eq:Asepsmudef} are met for $\epsilon_{\sigma} \lesssim \mu_{\sigma} ^{1/|\vec \alpha_{\sigma}|}$ and $\mu_{\sigma} \gtrsim 2^{r_{1, \sigma} |\vec \alpha_{\sigma}|} H_{f, K, B, \psi, \sigma}^{|\alpha|}.$ Hence, for $\sigma \in S_1$, we have
\begin{IEEEeqnarray*}{Cl}
& \left( \frac{\epsilon_{\sigma}}{\mu_{\sigma}} \right)^{1/(|\vec \alpha_{\sigma}| - 1)} \\
\lesssim & (\mu_{\sigma})^{-1/|\vec \alpha_{\sigma}|} \\
\lesssim & 2^{-r_{1,\sigma}} H_{f,K,B,\psi,\sigma}^{-1}.
\end{IEEEeqnarray*}
On the other hand, for $\sigma \in S_2$, we have from \eqref{eq:S2r} and \eqref{eq:S2l} that the estimates in equation \eqref{eq:Asepsmudef} are met with $\epsilon \sim 2^{r_{2, \sigma}} H_{f, K, B, \psi, \sigma}$ and $\mu \sim 2^{-|\vec \alpha_{\sigma}| l_{\sigma}} \epsilon^{|\vec \alpha_{\sigma}|}$. Hence, we have the estimate
\begin{IEEEeqnarray*}{Cl}
& \left( \frac{\epsilon_{\sigma}}{\mu_{\sigma}} \right)^{1/(|\vec \alpha_{\sigma}| - 1)} \\
\lesssim & \left(2^{|\vec \alpha_{\sigma}| l_{\sigma}} \epsilon^{1 - |\vec \alpha_{\sigma}|} \right)^{1/(|\vec \alpha_{\sigma}| - 1)} \\
\lesssim & 2^{\frac{|\vec \alpha_{\sigma}|}{|\vec \alpha_{\sigma}| - 1} l_{\sigma}} 2^{-r_{2, \sigma}} H_{f,K,B,\psi,\sigma}^{-1}.
\end{IEEEeqnarray*}
Hence, Corollary \ref{cor:sublevelkreal} provides the estimate
\[|\mathfrak{S}_{\bm{\vec \alpha}, S_1, S_2, S_3, \bm{r_1}, \bm{r_2}, \bm{l}}| \lesssim_{K, B, n, d} \prod_{\sigma \in S_1} 2^{-r_{1,\sigma}} H_{f, K, B, \psi, \sigma}^{-1} \cdot \prod_{\sigma \in S_2} 2^{\frac{|\vec \alpha_{\sigma}|}{|\vec \alpha_{\sigma}| - 1} l_{\sigma}} 2^{-r_{2, \sigma}} H_{f,K,B,\psi,\sigma}^{-1}.\]

So
\begin{IEEEeqnarray*}{Cl}
& \sum_{e \in E_{\bm{\vec \alpha}, S_1, S_2, S_3, \bm{r_1}, \bm{r_2}, \bm{l}}} |I_e|  \\
\lesssim_{K, B, n, d, \psi} & \prod_{\sigma \in S_1} 2^{-r_{1, \sigma}} H_{f, K, B, \psi, \sigma}^{-1} \cdot \prod_{\sigma \in S_2} 2^{\left(-N + \frac{|\vec \alpha_{\sigma}|}{|\vec \alpha_{\sigma}| - 1} \right) l_{\sigma}} 2^{-r_{2, \sigma}} H_{f,K,B, \psi, \sigma}^{-1} \cdot \prod_{\sigma \in S_3} H_{f, K, B, \psi, \sigma}^{-N} \\
\leq & H_{f, K, B, \psi}^{-1} \prod_{\sigma \in S_1} 2^{-r_{2,\sigma}} \prod_{\sigma \in S_2} 2^{\left(-N + \frac{|\vec \alpha_{\sigma}|}{|\vec \alpha_{\sigma}| - 1} \right) l_{\sigma}} 2^{-r_{1,\sigma}}.
\end{IEEEeqnarray*}
If $N \geq 4 >  \frac{|\vec \alpha_{\sigma}|}{|\vec \alpha_{\sigma}| - 1}$, then this expression is summable in $\bm{l}, \bm{r_1}$, and $\bm{r_2}$ to a value that is $\lesssim H_{f, K, B, \psi}^{-1}$. Moreover, there are only finitely many choices for $\bm{\vec \alpha}, S_1, S_2,$ and $S_3$. Hence 
\[\sum_{e} I_e \lesssim_{K, B, n, d, \psi}  H_{f,K,B,\psi}^{-1},\]
as desired.
\end{proof}
By keeping track of the dependence on $\psi$ above, we observe the following.
\begin{mycor}\label{cor:mainthmcor}
If $\psi$ is a function supported on the unit ball in $\mathbb{R}^{kn}$ and if $f$ is such that for some $N \geq 4$, we have the inequality  $J_{K,B,f,\sigma}(\bm{\vec x})^t \geq R \norm{\psi}_{C^t}$ for each $0 \leq t \leq kN$ and each $\sigma \in S$ and each $\bm{\vec x} \in \supp \psi$, then 
\begin{equation}\label{eq:corresult}
\int_{\mathbb{R}^{kn}} e(\phi_f(\bm{\vec x})) \psi(\bm{\vec x}) \lesssim_{K,B,n, d, R} H_{K,B,f,S, \psi}^{-1}
\end{equation}
where the implicit constant is independent of $\psi$.
\end{mycor}
\begin{proof}
Corollary \ref{cor:intbyparts} implies that the constant arising from the use of integration by parts in the above proof is independent of $\psi$.
\end{proof}

\section{A Fourier transform on $(\mathbb{R} \otimes_{\mathbb{Q}} K)_B^n$}\label{sec:FT}
In Sections \ref{sec:FT} to \ref{sec:Tarrysharpness} of this paper, we present an application of Theorem \ref{thm:mainthm} and Corollary \ref{cor:mainthmcor} to the singular integral in a version of Tarry's problem for algebraic number fields. We want to relate the Fourier transform on the algebra $(\mathbb{R} \otimes_{\mathbb{Q}} K)_B$ to the usual Fourier transform on $\mathbb{R}^k$.

If $\psi : \mathbb{R}^k \to \mathbb{C}$ is a function, we define the $(K,B)$-Fourier transform of $\psi$ by 
\begin{equation}\label{eq:KFourier}
\mathcal{F}_{K,B} \psi(\bm{x}) = \int \exp \left(2 \pi i \tr^*(A^*(\bm{x}) A^*(\bm{y})) \right) \psi(\bm{y}) \, d \bm{y}
\end{equation}
Recall that $A^*(\bm{x})$ and $A^*(\bm{y})$ are $k$-by-$k$ matrices with real entries, and thus the trace is real. Hence the phase $\tr^*(A^*(\bm{x}) A^*(\bm{y}))$ is real-valued and the oscillatory integral \eqref{eq:KFourier} does not require any special interpretation.

We will focus on the phase $\tr^*(A^*(\bm{x}) A^*(\bm{y}))$. This is a non-degenerate symmetric bilinear form, so there is an invertible, symmetric matrix $T$ such that
\begin{equation}\label{eq:Tequation}
\tr(A^*(\bm{x}) A^*(\bm{y})) = T \bm{x} \cdot \bm{y},
\end{equation}
where $\cdot$ denotes the usual dot product on $\mathbb{R}^k$. Thus
\begin{IEEEeqnarray*}{rCl}
\mathcal{F}_{K,B} \psi(\bm{x}) & = & \int \exp \left(2 \pi i \tr^*(A^*(\bm{x}) A^*(\bm{y})) \right) \psi(\bm{y}) \, d \bm{y} \\
& = & \int \exp(2 \pi i T \bm{x} \cdot \bm{y}) \psi(y) \, dy \\
& = & \hat \psi(T \bm{x}).
\end{IEEEeqnarray*}
So the $(K,B)$-Fourier transform is just a composition of the usual Fourier transform and an invertible, symmetric, linear change-of-variables.

In higher dimensions, the $(K,B)$-Fourier transform is defined in a similar manner. If $\psi : \mathbb{R}^{kn} \to \mathbb{C}$ is a function, then
\begin{equation}\label{eq:KFouriern}
\mathcal{F}_{K,B} \psi(\bm{\vec x}) = \exp \left( 2 \pi i \left(\tr^*(A^*(\bm{x}_1) A^*(\bm{y}_1)) + \cdots + \tr^*(A^*(\bm{x}_n) A^*(\bm{y}_n)) \right) \right)
\end{equation}
and, writing $\vec{T}(\mathbf{x})$ for $(T \bm{x}_1, \ldots, T \bm{x}_n)$, we have that the $(K,B)$-Fourier transform is related to the usual Fourier transform by
\begin{equation}\label{eq:Fourierlinear}
\mathcal{F}_{K,B} \psi(\bm{x}) = \hat \psi(\vec T \bm{\vec x}).
\end{equation}

\section{An application to Tarry's problem for algebraic number fields}\label{sec:Tarryintro}
Let $K$ be a number field of degree $k$ with basis $B = \{\omega_1, \ldots, \omega_k\}$. Then there are uniquely defined rational polynomials $Q_{l,j}$ such that for any $q_1, \ldots, q_k \in \mathbb{Q}$,
\begin{equation}\label{eq:qjl}
(q_1 \omega_1 + \cdots + q_k \omega_k)^l = \sum_{j = 1}^k Q_{l,j}(q_1, \ldots, q_k) \omega_j.
\end{equation}
The polynomials $Q_{l,j}$ can be viewed as real-valued functions of $k$ real variables. We define the \textbf{$(K,B)$-moment curve of degree $n$} to be the $k$-dimensional manifold in $\mathbb{R}^{nk}$ parameterized by the polynomials $\{Q_{l,j}\}_{1 \leq l \leq n, 1 \leq j \leq k}$. For simplicity of notation, we write $\bm{\vec Q}(\bm{x})$ for the vector consisting of the polynomials $Q_{l,j}$. Fix a smooth cutoff $\psi$ with the property that $\psi(\bm{x}) = 1$ when $|x| \leq 2$. To the $k$--moment curve, we associate an extension operator $E_{K,B,\psi}$ defined by 
\begin{equation}\label{eq:extension}
E_{K,B, \psi}g(\bm{\vec \xi}) = \int_{\mathbb{R}^{k}} \exp \left(2 \pi i \bm{\vec \xi} \cdot \bm{\vec Q(x)} \right) \, f(\bm{x}) \psi(\bm{x}) d \bm{x}
\end{equation}
In order to apply Theorem \ref{thm:mainthm} to the integral \eqref{eq:extension}, we need to recognize the phase as a trace of a polynomial in $(\mathbb{R} \otimes_{\mathbb{Q}} K)_B.$. By equations \eqref{eq:KFouriern} and \eqref{eq:Fourierlinear}, we have that 
\[E_{K,B,\psi} g(\bm{\vec \xi}) = \int_{\mathbb{R}^k} \exp \left(2 \pi i \tr^* \left( \sum_{l=1}^n A^*(T^{-1} \bm{\xi}_l) A^*(\bm{Q}_l (\bm{x})) \right) \right) g(\bm{x}) \psi(\bm{x}) \, d \bm{x}.\]
Writing $\bm{\vec \eta} = (\vec T)^{-1} \bm{\vec \xi}$, we have
\[E_{K,B, \psi} g(\xi) = \int_{\mathbb{R}^k} \exp(2 \pi i \phi_{f_{\bm{\vec \eta}(\bm{\vec \xi})}}(\bm{x}) )\psi(\bm{x}) g(\bm{x}) \, dx\]
where 
\[f_{\bm{\vec \eta}}(x) = A_{K,B}^*(\bm{\eta_1}) x + \cdots + A_{K,B}^*(\bm{\eta_n}) x^n.\]
Hence the oscillatory integral \eqref{eq:extension} is written in a form for which we can apply Theorem \ref{thm:mainthm}.
\section{Estimating the $H$-functional for Tarry's problem}\label{sec:HfuncTarry}

It remains to estimate the $H$-functional $H_{f_{\eta},K,B,\psi,S}$ for varying values of $\bm{\vec \eta}$. Observe that the polynomial $P_{f_{\bm{\vec \eta}}, \sigma}$ is given by
\[P_{f_{\bm{ \vec \eta}}, \sigma}(x) = \sigma^{**}(\bm{\eta_1}) x + \cdots + \sigma^{**}(\bm{\eta_n}) x^n.\]
We will use the following basic lemma of Arhipov, Čubarikov, and Karacuba (in the real case) and Wright (in the complex case).
\begin{mylem}\label{lem:basicpolynomial}
\begin{enumerate}[(1)]
\item (Real Case) Let $Q > 1$. For a vector $\vec a = (a_1, \ldots, a_n) \in \mathbb{R}^n$ , define the real-variable polynomial
\[P_{\vec a}(x) = a_1 x + \cdots + a_n x^n.\]
Let $\psi$ be a function supported on $[-1,1]$. Let $\mathfrak{A}(Q)$ denote the set of those vectors $\vec a$ such that $H_{P_{\vec a}, \mathbb{R}, \psi } \leq Q$ . The measure in $\mathbb{R}^n$  of $\mathfrak{A}(Q)$ is $\lesssim Q^{\frac{1}{2} (n^2 + n) + 1}$ .
Moreover, the set $\mathfrak{A}(Q)$ can be contained in a union of $\lesssim Q$ parallelepipeds $\mathfrak{P}_{b}(Q): b \in [-Q, Q] \cap \mathbb{Z}$ of measure $Q^{\frac{1}{2}n(n+1)}$ .

The parallelepipeds $\mathfrak{P}_{b}$ are of the following shape. There is an interval $J_n$ of length $\sim Q^n$ such that $a_n \in J_n$. For $1 \leq l \leq n-1$, $a_l$ is contained in an interval $J(b, a_{l+1}, \ldots, a_n)$ of length $\sim Q^{l}$.
\item (Complex Case) Let $Q > 1$. For a vector $\vec a = (a_1, \ldots, a_n) \in \mathbb{C}^n$ , define the complex-variable polynomial
\[P_{\vec a}(z) = a_1 z + \cdots + a_n z^n.\]
Let $\psi$ be a function supported on the unit disc. Let $\mathfrak{A}(Q)$ denote the set of those vectors $\vec a$ such that $H_{P_{\vec a}, \mathbb{C}, \psi } \leq Q$ . The measure in $\mathbb{R}^n$  of $\mathfrak{A}(Q)$ is $\lesssim Q^{(n^2 + n) + 2}$ .
Moreover, the set $\mathfrak{A}(Q)$ can be contained in a union of $\lesssim Q^2$ parallelepipeds $\mathfrak{P}_{b_1, b_2}(Q): b_1, b_2 \in [-Q, Q] \cap \mathbb{Z}$ of measure $Q^{n(n+1)}$.

The parallelepipeds $\mathfrak{P}_{b_1, b_2}$ are of the following shape. There is a ball$J_n$ of diameter $\sim Q^n$ such that $a_n \in J_n$. For $1 \leq l \leq n-1$, $a_l$ is contained in a ball $J(b_1, b_2, a_{l+1}, \ldots, a_n)$ of diameter $\sim Q^{l}$. Moreover, $\overline{\mathfrak{P}_{b_1, b_2}} = \mathfrak{P}_{b_1, -b_2}$ for any $b_1, b_2$. Hence $\vec a \in \mathfrak{A}(Q)$ if and only if $\overline{\vec a} \in \mathfrak{A}(Q)$.
\end{enumerate}

\end{mylem}
\begin{proof}
We present the proof appearing in Wright \cite{Wright20}. We will only present the complex case, as the real case is similar (and simpler).

For integers $b_1, b_2 \in [-Q,Q] \cap \mathbb{Z}$, define the point $z_{b_1, b_2} = \frac{b_1}{Q} + \frac{b_2}{Q} i$.  We claim that if $\vec a$ is such that $H_{P_{\vec a}, \mathbb{C}, \psi} \leq Q$, then there exists a point $z_{b_1, b_2}$ such that $H_{P_{\vec a}, \mathbb{C}}(z_{b_1, b_2}) \lesssim_n Q$, where the implicit constant does not depend on $Q$.

Since $H_{P_{\vec a}, \mathbb{C}, \psi} \leq Q$, there exists $z \in [-1,1] \times [-1,1] \subset \mathbb{C} $ such that $H_{P_{\vec a}, \mathbb{C}}(z) \leq Q$. Choose $b_1 = \lfloor Q \text{Re}(z) \rfloor$ and $b_2 = \lfloor Q \text{Im}(z) \rfloor$, so that $|z - z_{b_1, b_2}| \lesssim Q^{-1}$. We claim that $H_{P_{\vec a}, \mathbb{C}}(z_{b_1, b_2}) \lesssim_n Q$.

To see this, observe that for each $1 \leq l \leq n$, we have the Taylor expansion
\[\partial^l P_{\vec a}(z_{b_1, b_2}) = \sum_{l' = 0}^{n - l} \frac{1}{l'!} \partial^{l + l'} P_{\vec a}(z) (z_{b_1, b_2} - z)^{l'}.\]
Since $H_{P_{\vec a, \mathbb{C}}}(z) \leq Q$, it follows that each derivative $|\partial^{l + l'} P_{\vec a}(z)| \lesssim Q^{l + l'}$; since $|z_{b_1, b_2} - z| \lesssim Q^{-1}$, it follows that each summand above is bounded above by an $n$-dependent constant times $Q^l$. Hence $|\partial^l P_{\vec a}(z_{b_1, b_2})| \lesssim Q^l$ for each $1 \leq l \leq n$ and thus $H_{P_{\vec a}, \mathbb{C}}(z) \lesssim Q$.

Define $\mathfrak{A}_{b_1, b_2}(Q)$ to be the set of points $\vec a$ for which $H_{P_{\vec a}, \mathbb{C}}(z_{b_1, b_2}) \lesssim Q$. Then the above argument shows that the set of points $\mathfrak{A}$ such that $H_{P_{\vec a}, \mathbb{C}, \psi} \leq Q$ is contained in $\bigcup_{b_1, b_2} \mathfrak{A}_{b_1, b_2}$. Hence it only remains to check the Lebesgue measure of $\mathfrak{A}_{b_1, b_2}$.

The set $\mathfrak{A}_{b_1, b_2}(Q)$ consists of those vectors in $\vec a \in \mathbb{C}^n$ satisfying the following system of inequalities.
\begin{IEEEeqnarray*}{rCl}
|a_1 + 2 a_2 z_{b_1, b_2} + \cdots + n a_n z_{b_1, b_2}^{n-1}| & \lesssim &  Q \\
\vdots \qquad \vdots \qquad \vdots & \vdots & \vdots \\
\left|\frac{l!}{0!} a_l + \frac{(l + 1)!}{1!} a_{l+1} z_{b_1, b_2}^1 + \cdots + \frac{n!}{(n - l)!} a_n z_{b_1, b_2}^{n - l} \right| & \lesssim & Q^l \\
\vdots \qquad \vdots \qquad \vdots & \vdots & \vdots \\
|(n-1)! a_{n-1} + n! a_n z_{b_1, b_2}| & \lesssim & Q^{n-1} \\
|n! a_n| & \lesssim & Q^n.
\end{IEEEeqnarray*}
The last of these inequalities implies that $a_n$ lies in a disc  $J_n$ of radius $\sim Q^n$. The second--last of these inequalities implies that $a_{n-1}$ lies in a certain disc $J_{n-1}(b_1, b_2, a_n)$ of radius $\sim Q^{n-1}$. In fact, solving for $a_l$ in the $l$th inequality gives for $1 \leq l \leq n$ that $a_l$ lies in disc $J_{l}(b_1, b_2, a_{l+1}, \ldots, a_n)$ of radius $\sim Q^l$. Hence $\vec a$ must lie in a parallelepiped of measure $Q^{2 \left(\frac{n(n + 1)}{2} \right)}$. So $|\mathfrak{A}_{b_1, b_2}| \sim Q^{n(n+1)}$. Since there are $\sim Q^2$ choices for the pair $(b_1, b_2)$, it follows that $|\mathfrak{A}(Q)| \lesssim Q^{n(n+1) + 2}$.
\end{proof}
This lemma implies the following.
\begin{mylem}\label{lem:embeddingsublevel}
Let $\sigma \in \Sigma$ be an embedding, let $Q > 1$, and let $\bm{w_{\sigma}}$ denote the $k$-vector 
\[\bm{w_{\sigma}} = (\sigma(\omega_1), \sigma(\omega_2), \ldots, \sigma(\omega_k)).\]
Then $H_{f_{\eta}, K, B, \psi, \sigma} \geq Q$ provided that the $n$-vector $|\{\bm{\eta_l} \cdot \mathbf{w}_{\sigma}\}_{1 \leq l \leq n} | \notin \mathfrak{A}(Q)$.
\end{mylem}
It remains to characterize the set of $\bm{\eta_l}$ such that $\bm{\eta_l} \cdot \bm{w_{\sigma}} \in \mathfrak{A}(Q)$. Recall that even if $\sigma$ is a complex embedding, we are still interested in considering only real vectors $\bm{\vec \eta}$.
\begin{mylem}\label{lem:intervalfiber}
If $\sigma \in \Sigma$ is a real embedding and $J$ is any interval, then the set $\mathcal{N}_{l, \sigma}(J)$ of those $\bm{\eta_l}$ such that $\bm{\eta_l} \cdot \bm{w_{\sigma}} \in J$ is a $\lesssim |J|$-neighborhood of a translate of $\bm{w_{\sigma}}^{\perp}$.

If $\sigma \in \Sigma$ is a complex embedding and $J$ is any disc, then the set $\mathcal{N}_{l, \sigma}(J)$ of $\bm{\eta_l}$ such that $\bm{\eta_l} \cdot \bm{w_{\sigma}} \in J$ is a $\lesssim \diam(J)$-neighborhood of a translate of $(\bm{w_{\sigma}} + \bm{w_{\bar \sigma}})^{\perp} \cap (\frac{1}{i} (\bm{w_{\sigma}} - \bm{w_{\bar \sigma}}))^{\perp}$. Furthermore $\mathcal{N}_{l, \sigma}(J) = \mathcal{N}_{l, \bar \sigma}(\bar J)$.
\end{mylem}

Given a set $S \subset \Sigma$ with the property that $\bar \sigma \in S$ for any complex embedding $\sigma \in S$, and any $|S|$ vector $\bm{\alpha} = \{\alpha_{\tilde \sigma}\}_{\tilde \sigma \in \tilde S}$ of positive integers,  we define the set $\mathfrak{S}_{S, \bm{\alpha}}$ to be the set of those $\bm{\vec \eta} \in \mathbb{R}^{nk}$ such that
\begin{IEEEeqnarray*}{rCll}
2^{\alpha_{\sigma}} \leq H_{f_{\eta},K,B,\psi, \sigma}  & < & 2^{\alpha_{\sigma} + 1} \quad & \text{if $\sigma \in S$} \\ 
H_{f_{\eta},K,B,\psi, \sigma} & < & 1 \quad &  \text{if $\sigma \notin S$.}
\end{IEEEeqnarray*}
This is well--defined because $H_{f_{\eta}, K, B, \psi \sigma} = H_{f_{\eta}, K, B, \psi, \bar \sigma}.$
We will use the convention that $\mathfrak{S}_{\emptyset}$ is the set of points for which $H_{f_{\eta}, K, B, \psi \sigma} \geq 1$ for each $\sigma \in \Sigma$.
\begin{mycor}
The set $\mathfrak{S}_{S, \bm{\alpha}}$ has measure $\lesssim \prod_{\sigma \in S}2^{\left(\frac{1}{2} n(n + 1) + 1 \right) \alpha_{\sigma }}.$
\end{mycor}
\begin{proof}
Suppose $\bm{\vec \eta} \in \mathfrak{S}_{S, \bm{\alpha}}$. Then for each real embedding $\sigma \in S$, there exists a $b(\sigma)$ such that $\{\bm{\eta_l} \cdot \bm{w_{\sigma}}\}_{1 \leq l \leq n} \in \mathfrak{P}_{b(\sigma)}(2^{\alpha_{\sigma}})$, and for each conjugate pair $\tilde \sigma$ of strictly complex embeddings there exists $b_1(\tilde \sigma), b_2(\tilde \sigma)$ such that $\{\bm{\eta_l} \cdot \bm{w_{\sigma}}\}_{1 \leq l \leq n}$ is contained in $\mathfrak{P}_{b_1(\sigma), b_2(\sigma)}(2^{\alpha_{\sigma}})$. For simplicity, we will just write $b(\sigma)$ for the pair $(b_1(\sigma), b_2(\sigma))$ in the complex case. Hence for each $\sigma$, we have that $\bm{\eta_n}$ lies in the intersection
\[\bm{\eta_n} \in \bigcap_{\tilde \sigma \in \tilde S} \mathcal{N}_{n, \tilde \sigma}(J_n(2^{\alpha_{\sigma}})) \cap \bigcap_{\tilde \sigma \notin \tilde S} \mathcal{N}_{n, \tilde \sigma}(J_n(1)) .\] 
Because the sets $\mathcal{N}_{n,\tilde \sigma}(Q)$ are $\sim Q^n$-neighborhoods of linearly independent affine subspaces of codimension $1$ for real embeddings $\sigma$ and codimension $2$ for pairs of complex embeddings, it follows that the above set has measure $\lesssim_{K,B} \prod_{\sigma \in S} 2^{n \alpha_{\sigma}}$. Observe that for strictly complex embeddings both $\sigma$ and $\overline \sigma$ are counted in this product because $\tilde \sigma$ corresponds to a codimension $2$ affine subspace.
Then, reverse inductively, for $1 \leq l \leq n - 1$, we determine that 
\[\bm{\eta_l} \in \bigcap_{\tilde \sigma \in \tilde S} \mathcal{N}_{l, \sigma}(J_n(2^{\alpha_{\sigma}}, b(\sigma) , \sigma(\bm{\eta_{l+1}}), \ldots, \sigma(\bm{\eta_n}))) \cap \bigcap_{\tilde \sigma \notin \tilde S} \mathcal{N}_{l, \sigma}(J_n(1, b(\sigma) , \sigma(\bm{\eta_{l+1}}), \ldots, \sigma(\bm{\eta_n}))) \]
Each set in the intersection is a neighborhood of a linearly independent affine subspace of codimension $1$ or $2$, so $\bm{\vec \eta}_l$ is contained in a set of measure $\lesssim \prod_{\sigma \in S} 2^{l \alpha_{\sigma}}$. Hence the set of $\bm{\vec \eta}$ such that $\{\bm{\eta_l} \cdot \bm{w_{\sigma}}\}_{1 \leq l \leq n}$ is contained in $\mathfrak{P}_{b_1(\sigma), b_2(\sigma)}$ has measure at most $\prod_{\sigma \in S} 2^{\frac{1}{2} \alpha_{\sigma}n(n+1)}$. Since there are $\sim 2^{\alpha_{\tilde \sigma}}$ choices of $b(\sigma)$ for each real embedding $\sigma$ and $\sim 2^{2 \alpha_{\tilde \sigma}}$ choices of $b(\sigma)$ for each pair of complex embeddings $\sigma, \bar \sigma$, it follows that the total measure of $\mathfrak{S}_{S, \bm{\alpha}}$ is no more than $\prod_{\sigma \in S} 2^{\left(\frac{1}{2} n(n + 1) + 1 \right) \alpha_{\sigma}}$, as desired.
\end{proof}

We are finally ready to show the integrability of $E_{K,B,\psi} 1$.
\begin{mythm}
The function $E_{K,B, \psi} 1$ lies in $L^q$ for $q > \frac{n(n+1)}{2} + 1$.
\end{mythm}
\begin{proof}
We need to estimate the integral
\[\int_{\bm{\vec \xi} \in \mathbb{R}^{kn}}\left| \int_{\mathbb{R}^k} \exp \left(2 \pi i \phi_{f_{\bm{\vec \eta}(\bm{\vec \xi})}}(\bm{x}) \right) \psi(\bm{x}) \, d\bm{x} \right|^{q} \, d \bm{\vec \xi}.\]
Since $\bm{\vec \eta}(\bm{\vec \xi})= (\vec T)^{-1} \bm{\vec \xi}$, where $\vec T$ is an invertible linear transformation, we have that this integral is equal to a constant times
\[\int_{\bm{\vec \eta} \in \mathbb{R}^{kn}} \left| \int_{\mathbb{R}^k} \exp \left(2 \pi i \phi_{f_{\bm{\vec \eta}}}(\bm x) \right) \psi(\bm{x}) \, d \bm{x} \right|^q \, d \bm{\vec \eta}.\]
By Theorem \ref{thm:mainthm}, the inner integral is bounded above by a constant times $H_{f_{\bm{\vec \eta}}, K, B, \psi, S}^{-1}$ for any $S \subset \Sigma$. We define $H_{\bm{\vec{\eta}}, K, B, \psi, \emptyset} = 1$ for convenience. The choice of $S$ is allowed to depend on $\eta$, and we choose $S(\eta)$ depending on which set $\mathfrak{S}_{S, \bm{\alpha}}$ contains $\eta$. In other words, we estimate the integral above by
\begin{IEEEeqnarray*}{Cl}
& \sum_{\substack{S \subset \Sigma \\ S = \bar S}} \sum_{\bm{\alpha}} \int_{\bm{\vec \eta} \in \mathfrak{S}_{S, \bm{\alpha}}} H_{f_{\bm{\vec \eta}}, K, B, \psi, S}^{-q} \\
 \lesssim & \sum_{\substack{S \subset \Sigma \\ S = \bar S}} \sum_{\bm{\alpha}} \int_{\bm{\vec \eta} \in \mathfrak{S}_{S, \bm{\alpha}}} \left(\prod_{\sigma \in S}2^{-q\alpha_{\sigma}} \right) \\
\lesssim & \sum_{\substack{S \in \Sigma \\ S = \bar S}} \sum_{\bm{\alpha}} \left(\prod_{\sigma \in S} 2^{\left(\frac{1}{2} n (n + 1) + 1 \right) \alpha_{\sigma}} \cdot 2^{- q \alpha_{\sigma}} \right),
\end{IEEEeqnarray*}
which is finite provided that $q > \frac{1}{2} n(n + 1) + 1$.
\end{proof}
\begin{myrmk}
Wright \cite{Wright20} is able to obtain a pointwise bound $|E 1(\xi)| \lesssim (1 + |\xi|)^{\frac{-4n}{n^2 + n + 2}}$ for the extension operator associated to the complex moment curve. This shows in $n = 2$ complex dimensions that the complex moment curve (i.e. the ``parabola") is a Salem set of Hausdorff dimension $2$ in $\mathbb{R}^4$.

However, a similar bound is not available in our setting as the pointwise bound on $|E_{K,B, \psi} 1(\xi)|$ where $\bm{\vec \eta}(\bm{\vec \xi}) \in \mathfrak{S}_{S} := \bigcup_{\bm{\alpha}} \mathfrak{S}_{S, \bm{\alpha}}$ is only on the order of $(1 + |\xi|)^{\frac{-2 |S|n}{n^2 + n + 2}}$. Because there are arbitrarily large $\bm{\vec{\eta}} \in \mathfrak{S}_S$ for any nonempty $S$ with $S = \overline{S}$, this bound only shows that the $n$-dimensional $K$-moment curve has Fourier dimension at least $\frac{2}{n^2 + n + 2}$ if $K$ has at least one real embedding, or $\frac{4}{n^2 + n + 2}$ if $K$ has no real embeddings. Even in the $n = 2$ case, this does not yield a $k$-dimensional Salem set in any case for which $k \geq 3$. However, the sets $\mathfrak{S}_S$ for $S \neq \Sigma$ are sufficiently small that we have the same $L^q$-integrability for $E_{K,B,\psi} 1$ as in the real or complex cases.
\end{myrmk}
\section{Sharpness of the integrability exponent for Tarry's problem}\label{sec:Tarrysharpness}
We are able to show that the integrability exponent of Tarry's problem is exactly $\frac{n(n+1)}{2} + 1$ as in the real and complex cases. This proof follows the proof of Wright \cite{Wright20} for the complex case, but there are more more technicalities because we must consider each embedding of $K$ into $\mathbb{C}$.

Following Wright, we choose $A > 2$ and fix a lacunary sequence $Q_m = A^m$. We define the set $E_m$ to be the set of points $\bm{x} \in \mathbb{R}^d$ such that $Q_m^n \leq |\bm{x}| \leq (2 Q_m)^n $ and $|\sigma^{**}(\bm{x})| \geq C_{K,B} Q_m^{n}$ for some constant $C_{K,B}$. A straightforward argument in fact shows that $|E_m| \gtrsim Q_m^{kn}$. 

For $1 \leq r,s \leq Q_m/10$, we define 
\[\mathcal{P}_m := \{\bm{r} \in \mathbb{N}^k : 1 \leq r_j \leq Q_m/10 \quad \text{for $1 \leq j \leq k$.}\}\]

For $\bm{r} \in \mathcal{P}_m$, define the point $\bm{x_r} = \frac{1}{Q_m} \bm{r}$. Recall that $A_{K,B}^*$ is an invertible linear map from $\mathbb{R}^k$ into $(K \otimes_{\mathbb{Q}} \mathbb{R})_B$. Given an $nk$-vector $\bm{\vec \eta} = (\bm{\eta_1}, \ldots, \bm{\eta_n})$, write $\bm{\eta}'$ for the $(n-1)k$-vector $(\bm{\eta_1}, \ldots, \bm{\eta_{n-1}})$. For $\bm{r} \in \mathcal{P}_m$ and $\bm{\eta_n} \in E_m$, We define the affine-linear transformation $T_{\bm{r}, \bm{\eta_n}}: \mathbb{R}^{(n-1)k} \to \mathbb{R}^{(n-1)k}$  by $T_{\bm{r}, \bm{\eta_n}} \bm{\vec \theta}' = \bm{\vec \eta}'$, where
\[A^*(\bm{\eta_n}) z^n + \cdots + A^*(\bm{\eta_1}) z = A^*(\bm{\eta_n}) (z - A^*(\bm{x_r}))^n + \sum_{l=1}^{n-1} A^*(\bm{\theta_{l}}) (z - A^*(\bm{x_r}))^l.\]
For simplicity of notation, we will write $\bm{\theta_n} = \bm{\eta_n}$. As in the complex case, the linear part of the map $T_{\bm{r}, \bm{\eta_n}}$ is upper-triangular with $1$'s along the diagonal, and thus $T_{\bm{r}, \bm{\eta_n}}$ is an invertible transformation with determinant $1$.

Fix a small, positive constant $c_1$. Given the map $T_{\bm{r}, \bm{\eta_n}}$, we define the set $R_{m, \bm{r}}'(\bm{\eta_n}) = T_{\bm{r}, \bm{\eta_n}} \mathcal{B}$, where
\[\mathcal{B} := \{\bm{\vec \theta}' =  (\bm{\theta_1}, \ldots, \bm{\theta_{n-1}}) \in \mathbb{R}^{(n-1)k} : |\bm{\theta_l}| \leq (c_1 Q_m)^l; 1 \leq l \leq n-1\}.\]

The point of the set $R_{m, \bm{r}}'(\bm{\eta_n})$ is that the $H$-functionals $H_{f_{\bm{\vec \eta}},K,B,\sigma}(\bm{\vec x_r})$ are small when $\bm{\vec \eta'} \in R_{m, \bm{r}}'(\bm{\eta_n})$.

\begin{mylem}\label{lem:Rmrhfunc}
If $\bm{\eta_n} \in E_m$ and $\bm{\vec \eta'} \in R_{m, \bm{r}}'(\bm{\eta_n})$, then $H_{f,K,B,\sigma}(\bm{x_r}) \lesssim Q$ for all $\sigma \in \Sigma$.
\end{mylem}
\begin{proof}
Let $\sigma \in \Sigma$. By the definition of $R_{m,\bm{r}}'(\bm{\eta_n})$, there exists $\bm{\vec \theta} \in \mathcal{B}$ such that $\bm{\vec \eta}' = T_{\bm{r}, \bm{\eta_n}}(\bm{\vec \theta'})$. Thus
\[P_{f_{\bm{\vec \eta}}, \sigma}(\sigma^{**}(\bm{x})) = \sigma^{**} (\bm{\theta_1}) (\sigma^{**}(\bm{x}) - \sigma^{**}(\bm{x_r})) + \cdots + \sigma^{**}(\bm{\theta_n}) (\sigma^{**}(\bm{x}) - \sigma^{**}(\bm{x_r}))^n.\]
We use the convention that $\bm{\theta_n} = \bm{\eta_n}$. The $l$th derivative is
\[\partial^l P_{f_{\bm{\vec \eta}}, \sigma} (\sigma^{**}(\bm{x})) = \sum_{m=0}^{n - l} \frac{(m + l)!}{m!} \sigma^{**}(\bm{\theta_m}) (\sigma^{**}(\bm{x}) - \sigma^{**}(\bm{x_r}))^m.\]
At the point $\bm{x} = \bm{x_r}$, all of the terms in this sum except for the $m = 0$ term vanish. Hence
\[\partial^l P_{f_{\bm{\vec \eta}}, \sigma} (\sigma^{**}(\bm{x_r})) = l! \sigma^{**}(\bm{\theta_m}).\]
Since $\bm{\theta_n} = \bm{\eta_n} \in E_m$, we have that $|\sigma^{**}(\bm{\theta_n})| \sim Q^n$. Moreover, since $\bm{\vec \theta}' \in \mathcal{B}$, it follows that $|\bm{\theta}_l| \leq (c_1 Q_m)^l$ for each $1 \leq l \leq n$. But
\[|\sigma^{**}(\bm{\theta_l})| = \left| \sum_{j=1}^k \theta_{l,j} \sigma(\omega_j) \right| \lesssim_{K,B} |\bm{\theta_l}| \leq c_1 Q^l\]
as desired.
\end{proof}

The set $R_{m}'(\bm{\eta_n})$ is defined to be the union
\[R_m'(\bm{\eta_n}) = \bigcup_{\bm{r} \in \mathcal{P}_m} R_{m, \bm{r}}'(\bm{\eta_n})\]
Finally, we define the set $R_m$ by 
\[R_m = \{(\bm{\vec \eta}', \bm{\eta_n}) : \bm{\eta_n} \in E_m; \bm{\vec \eta}' \in R_m(\bm{\eta_n})\}.\]

For simplicity of notation, write
\[I(\bm{\vec \eta}) = \int_{\bm{x} \in \mathbb{R}^k} e(\phi_{f_{\bm{\vec \eta}}}(\bm{x})) \psi(\bm{x}) \, d \bm{x}\]
Recall that $\psi(\bm{x}) = 1$ when $|x| \leq 2$. Our goal is to show the following result.
\begin{mythm}\label{thm:tarrysharpness}
Let $q_n = \frac{n(n + 1)}{2} + 1$. Then $E_{K,B,\psi}(1)$ does not lie in $L^q$ for any $q \leq q_n$. In other words, if $q \leq q_n$, then 
\[\int_{\bm{\vec \eta}} |I(\bm{\vec \eta})|^q = \infty.\]
\end{mythm}

We will show this estimate by obtaining a lower bound on $|I(\bm{\vec \eta})|$ on $R_m$. First, we verify that for a fixed $\bm{\eta_n}$ that the sets $R_{m,\bm{r}}'(\bm{\eta_n})$ are disjoint. 
\begin{mylem}
Suppose the constant $c_1$ is chosen to be sufficiently small. Let $\bm{\eta_n} \in E_m$, and suppose there exists a vector $\bm{\vec \eta}' \in \mathbb{R}^{k (n-1)}$ and vectors $\bm{r_1}, \bm{r_2} \in \mathcal{P}_m$  such that $\bm{\vec \eta}' \in R_{m, \bm{r_1}}'(\bm{\eta_n}) \cap R_{m, \bm{r_2}}'(\bm{\eta_n})$. Then $\bm{r_1} = \bm{r_2}$. Hence for a fixed $\bm{\eta_n}$, the sets $R_{m, \bm{r}}' : r \in \mathcal{P}_m$ are disjoint.
\end{mylem}
\begin{proof}
Suppose $\bm{r_1} \neq \bm{r_2}$. Then $|\bm{x_{r_1}} - \bm{x_{r_2}}| \geq Q_m^{-1}$. Hence there exists an embedding $\sigma \in \Sigma$ such that
\begin{equation}\label{eq:sigmadiff}
|\sigma^{**}(\bm{x_{r_1}}) - \bm{x_{r_2}}| \geq c_{K,B}' Q_m
\end{equation}
where $c_{K,B}'$ is some constant depending only on the field $K$ and the basis $B$. 

If $\bm{\eta'} \in R_{m, \bm{r_1}}' (\bm{\eta_n}) \cap R_{m, \bm{r_2}}' (\bm{\eta_n})$, then we have that
\[\left|(n-1)! \sigma^{**}(\bm{\eta}_{n-1}) + n! \sigma^{**}(\bm{\eta_n}) \sigma^{**}(\bm{x_{r_1}})\right| \leq (c_1 Q_m)^{n-1}\]
and
\[\left|(n-1)! \sigma^{**}(\bm{\eta}_{n-1}) + n! \sigma^{**}(\bm{\eta_n}) \sigma^{**}(\bm{x_{r_2}})\right|\leq (c_1 Q_m)^{n-1}\]
Hence
\[\left|\sigma^{**}(\bm{x_{r_1}} - \bm{x_{r_2}}) \right| \leq \frac{c_1^{n-1} Q_m^{n-1}}{n! |\sigma^{**}(\bm{\eta_n})|}.\]
Since $\bm{\eta_n} \in E_m$, we know that
\[|\sigma^{**}(\bm{\eta_n})| \geq C_{K,B} Q_m^n.\]
Hence
\[\left|\sigma^{**}(\bm{x_{r_1}} - \bm{x_{r_2}}) \right| \leq \frac{c_1^{n-1} Q_m^{-1}}{n! C_{K,B}}.\]
This contradicts the inequality \eqref{eq:sigmadiff} if $c_1$ is chosen to be sufficiently small.
\end{proof}
Therefore,
\begin{IEEEeqnarray*}{rCl}
\int_{\bm{\vec \eta}} |I(\bm{\vec \eta})|^q \, d \bm{\vec \eta} & \geq & \sum_{m \geq m_0} \int_{R_m} |I(\bm{\vec \eta})|^q \, d \bm{\vec \eta} \\
& = & \sum_{m \geq m_0} \sum_{\bm{r} \in \mathcal{P}_m} \int_{\bm{\eta}_n \in E_m} \int_{\bm{\vec \eta}' \in R_{m, \bm{r}}'(\bm{\eta_n})} |I(\bm{\vec \eta})|^q d \bm{\vec \eta}
\end{IEEEeqnarray*}
We will focus on the innermost integral
\[\int_{\bm{\vec \eta'} \in R_{m,\bm{r}}'(\bm{\eta_n})} |I(\bm{\vec \eta})|^q d \bm{\vec \eta}.\]
We make the change of variables $\bm{\vec \eta} = T_{\bm{r}, \bm{\eta_n}}(\bm{\vec \theta})$ (again taking $\bm{\theta_n}= \bm{\eta_n}$) to arrive at the integral
\[\int_{\bm{\vec \eta'} \in R_{m,\bm{r}}(\bm{\eta_n})} |I(\bm{\vec \eta})|^q d \bm{\vec \eta} = \int_{|\bm{\theta_{n-1}}| \leq (c_1 Q_m)^{n-1}} \cdots \int_{|\bm{\theta_1}| \leq c_1 Q_m} |\mathrm{II}_{\bm{r}}(\bm{\vec \theta})|^q d \bm{\vec \theta},\]
where
\[\mathrm{II}_{\bm{r}}(\bm{\vec \theta}) = \int_{\bm{x} \in \mathbb{R}^k} e(\phi_{f_{\bm{\vec \theta}}}(\bm{x} - \bm{x_{r}})) \psi(\bm{x}) \, d\bm{x}\] 
After translating by $\bm{x_{r}}$, we have
\[\mathrm{II}_{\bm{r}}(\bm{\vec \theta}) = \int_{\bm{x} \in \mathbb{R}^k} e(\phi_{f_{\bm{\vec \theta}}}(\bm{x})) \psi(\bm{x} + \bm{x_r}) \, d\bm{x}.\] 
We claim the following lemma. 
\begin{mylem}\label{lem:IIlem}
Suppose $\bm{\vec \theta}$ is such that $\bm{\theta}_n \in E_m$ and $\bm{\theta_j} \leq (c_1 Q_m)^j$ for every $j \leq n - 1$. Then $\mathrm{II}_{\bm{r}}(\bm{\vec \theta}) \gtrsim Q_m^{-k}$, where the implicit constant is independent of $m$ and $\bm{\vec \theta}$.
\end{mylem}
\begin{proof}[Proof of Theorem \ref{thm:tarrysharpness} assuming Lemma \ref{lem:IIlem}]
Suppose that $|\mathrm{II}_{\bm{r}}(\bm{\vec \theta})| \gtrsim Q_m^{-k}$ when $\bm{\theta_n} \in E_m$ and $\bm{\theta_j} \leq (c_1 Q_m)^j$ for each $1 \leq j \leq n-1$. Thus, for $\bm{\theta_n} \in E_m$, we have the estimate
\[\int_{|\bm{\theta_{n-1}|} \leq (c_1 Q_m)^{n-1}} \cdots \int_{|\bm{\theta_1}| \leq c_1 Q_m} |\mathrm{II}_{\bm{r}}(\bm{\vec \theta})|^q \, d \bm{\vec \theta} \gtrsim (c_1 Q_m)^{\frac{kn(n-1)}{2}} \cdot Q_m^{-kq}.\]
Integrating this over $\bm{\theta_n} \in E_m$, and summing over $\bm{r} \in \mathcal{P}_m$, we get
\[\int_{\bm{\vec \eta} \in R_m} |I(\bm{\vec \eta})|^q \, d \bm{\vec \eta} \gtrsim Q_m^k Q_m^{kn} Q_m^{\frac{kn(n-1)}{2}} Q_m^{-kq} = Q_m^{k \left(\frac{n(n+1)}{2} + 1 \right) - k q}.\]
Summing over $m \geq m_0$, we see that the sum diverges whenever 
\[q \leq \frac{n(n+1)}{2} + 1.\]
This proves the sharpness of the exponent in Tarry's problem.
\end{proof}
It only remains to prove Lemma \ref{lem:IIlem}. 
Fix a nonnegative function $\chi \in C_c^{\infty}(\mathbb{R}^k)$ supported on $\mathbf{B}(\bm{0};2)$ such that $\chi \equiv 1$ on $\mathbf{B}(\bm{0}; 1)$. Write $\chi_m(\bm{x}) = \chi \left(\frac{Q_m}{a} \bm{x} \right)$, where $a$ is a large constant. Here, we will assume that $a$ and $c_1$ satisfy the condition that
\begin{equation}\label{eq:c1acondition}
c_1 a^{k+1} \ll 1.
\end{equation}

We write
\[\mathrm{II}_{\bm{r}}^1(\bm{\vec \theta}) = \int_{\mathbb{C}} e( \phi_{f_{\bm{\vec \theta}}}(\bm{x})) \psi(\bm{x} + \bm{x_r}) \chi_m(\bm{x}) \, d \bm{x} \]
and
\[\mathrm{II}_{\bm{r}}^2(\bm{\vec \theta}) = \int_{\mathbb{C}} e( \phi_{f_{\bm{\vec \theta}}}(\bm{x})) \psi(\bm{x} + \bm{x_r}) (1 - \chi_m(\bm{x})) \, d \bm{x}. \]
We will first obtain an upper bound on $\mathrm{II}_{\bm{r}}^2(\bm{\vec \theta})$. We expect that $\mathrm{II}_{\bm{r}}^2(\bm{\theta})$ will be small since the support of $1 - \chi_m(\bm{\vec x})$ is far from $0$, which is the only place where every derivative of $P_{f,\sigma}$ is small for every $\sigma \in \Sigma$. Hence, we expect that the $H$-functional will be large for this integral.
\begin{mylem}\label{lem:II2lem}
Let $\mathrm{II}_{\bm{r}}^2(\bm{\vec \theta})$ be as above. Then 
\begin{equation}\label{eq:II2lem}
|\mathrm{II}_{\bm{r}}^2(\bm{\vec \theta})| \lesssim_{K,B,n, \psi} a^{-\frac{1}{n-1}} Q_m^{-k}.
\end{equation}
\end{mylem}
\begin{proof}[Proof of Lemma \ref{lem:II2lem}]
Observe that there exists a small constant $c_{K,B}$ such that for each $\bm{x} \in \mathbb{R}^k$, there exists $\sigma \in \Sigma$ such that $|\sigma^{**}(\bm{x})| \geq c_{K,B} |\bm{x}|$. Hence $\mathbb{R}^k \setminus \{0\}$ is covered by open sets $U_{\sigma} = \{\bm{x} \in \mathbb{R}^k : \sigma^{**}(\bm{x}) >\frac{c_{K,B}}{2} |\bm{x}|\}$. Let $\{\Delta_{\sigma} : \sigma \in \Sigma\}$ be a smooth partition of unity subordinate to this cover and write
\[\Psi_{m,\bm{r}, \sigma}(\bm{x}) = \psi(\bm{x} + \bm{x_r}) (1 - \chi_m(\bm{x})) \Delta_{\sigma}(x).\]
It is easy to see for each $N \geq 1$ that $\norm{\Psi_{m, \bm{r}, \sigma}}_{C^N} \leq \left(\frac{Q_m}{a} \right)^N$. 

Fix $\sigma_0 \in \Sigma$. We will obtain a bound on $H_{f_{\bm{\vec \theta}} ,K,B,\sigma, \Psi_{m, \bm{r}, \sigma_0}}$ when $\sigma = \sigma_0$ and a different bound when $\sigma \neq \sigma_0$.  For $\sigma \neq \sigma_0$, we observe that the derivative
\[\partial^n P_{f_{\bm{\vec \theta}}, \sigma}(\sigma^{**}(\bm{x})) = |n! \sigma^{**}(\bm{\theta_n})| \gtrsim Q_m^n\]
since $\bm{\theta}_n \in E_m$. Hence for $\sigma \neq \sigma_0$, we have the bound
\begin{equation}\label{eq:Hfunctionalsigmanot0}
H_{f_{\bm{\vec \theta}}, K, B, \sigma, \Psi_{m, \bm{r}, \sigma_0}} \gtrsim Q_m,
\end{equation}
with a similar bound for the corresponding $J$--functional.

When $\sigma = \sigma_0$, we obtain a better bound by looking at the $(n-1)$st derivative instead of the $n$th derivative. On the support of $\Psi_{m, \bm{r}, \sigma}$, the $(n - 1)$st derivative is given by
\begin{IEEEeqnarray*}{rCl}
\left| \partial^{(n-1)} P_{f_{\bm{\vec \theta}, \sigma_0}}(\sigma^{**}(\bm{x})) \right| & = & \left| n! \sigma_0^{**} (\bm{\theta_n}) \sigma_0^{**}(\bm{x}) + (n - 1)! \sigma_0^{**}(\bm{\theta_{n-1}})\right| \\
& \geq & \left| n! C_{K,B}^{-1} c_{K,B} a  Q_m^{n-1} + (n-1)! (c_1 Q_m)^{n-1} \right| \\
& \gtrsim_{K,B} & a Q_m^{n-1}
\end{IEEEeqnarray*}
since $a$ is large and $c_1$ is small. 

Hence for $\sigma = \sigma_0$, we have the bound
\begin{equation}\label{eq:Hfunctionalsigma0}
H_{f_{\bm{\vec \theta}}, K, B, \sigma_0, \Psi_{m, \bm{r}, \sigma_0}} \gtrsim_{K,B} a^{\frac{1}{n-1}} Q_m,
\end{equation}
with a similar bound for the $J$-functional. Note that this is significantly larger than the bound we obtain for $\sigma \neq \sigma_0$.

Combining \eqref{eq:Hfunctionalsigmanot0} and \eqref{eq:Hfunctionalsigma0}, we have that
\[H_{f_{\bm{\vec \theta}, K, B, \Sigma, \Psi_{m, \bm{r}, \sigma_0}}} \gtrsim_{K,B} a^{\frac{1}{n-1}} Q_m^k.\]
Moreover, since \eqref{eq:Hfunctionalsigmanot0} and \eqref{eq:Hfunctionalsigma0} also hold for the $J$-functional, we can apply Corollary \ref{cor:mainthmcor} to obtain the estimate
\[\left|\mathrm{II}_{\bm{r}}^2 \right| \lesssim_{K,B,n} a^{-\frac{1}{n-1}} Q_m^{-k}.\]
\end{proof}

Now, we obtain a \textit{lower} bound on $\mathrm{II}_{\bm{r}}^1(\bm{\vec \theta})$. Since $\chi_m$ is supported on $\mathbf{B}(\bm{0}; 2a Q_m^{-1})$, we have that $|\bm{x} + \bm{x_r}| \leq 1/2$ on the support of $\chi_m$. Hence $\psi(\bm{x} + \bm{x_r}) \equiv 1$ on the support of $\chi_m$. Thus, for all $\bm{x}$,
\[\psi(\bm{x} + \bm{x_r}) \chi_m(\bm{x}) = \chi_m(\bm{x})\]
Hence
\[\mathrm{II}_{\bm{r}}^1 = \int_{\mathbb{R}^k} e(\phi_{f_{\bm{\vec \theta}}}(\bm{x})) \chi_m(\bm{x}) \, d \bm{x}.\]
Define $g_{\bm{\vec \theta}}(\bm{x}) = g_{\bm{\theta_n}}(\bm{x}) = \bm{\theta_n}x^n$ for $x \in (K \otimes_{\mathbb{Q}} \mathbb{R})_B$. 

Write
\[\mathrm{II}_{\bm{r}}^1(\bm{\vec \theta}) = \int e \left(\phi_{g_{\bm{\theta_n}}}(\bm{x}) \right) \chi_m(\bm{x}) \, d \bm{x} + E,\]
where $E$ is an error term given by
\[|E| = \left| \int_{\mathbb{R}^k} \left[ e(\phi_{f_{\bm{\vec \theta}}}(\bm{x})) - e(\phi_{g_{\bm{\theta_n}}}(\bm{x})) \right] \chi_m(\bm{x}) \, d \bm{x} \right|. \]
\begin{mylem}\label{lem:Elem}
Let $E$ be as above. Then $E$ satisfies the bound
\[|E| \lesssim_{K,B} a^{k+1} c_1 Q_m^{-k}.\]
Note that this is significantly smaller than $Q_m^{-k}$ in light of \eqref{eq:c1acondition}.
\end{mylem}
\begin{proof}[Proof of Lemma \ref{lem:Elem}]
The error term $E$ satisfies the estimate
\begin{IEEEeqnarray*}{rCl}
|E| & = & \left| \int_{\mathbb{R}^k} \left[ e(\phi_{f_{\bm{\vec \theta}}}(\bm{x})) - e(\phi_{g_{\bm{\theta_n}}}(\bm{x})) \right] \chi_m(\bm{x}) \, d \bm{x} \right| \\
& \lesssim & \int_{\mathbf{B}(\bm{0}; 2a Q_m^{-1})} \left| \phi_{f_{\bm{\vec \theta}}}(\bm{x}) - \phi_{g_{\bm{\theta_n}}}(\bm{x}) \right| \, d \bm{x} \\
& \lesssim & \int_{\mathbf{B}(\bm{0}; 2 a Q_m^{-1})} \sum_{\sigma \in \Sigma} \left| P_{f_{\bm{\vec \theta}}, \sigma}(\sigma^{**}(\bm{x})) - P_{g_{\bm{\theta_n}}, \sigma}(\sigma^{**}(\bm{x})) \right| \, dx
\end{IEEEeqnarray*}
Recall that $c_1$ is a small constant and $a$ is a large constant. By \eqref{eq:c1acondition}, $c_1 a \ll 1$. Thus
\begin{IEEEeqnarray*}{Cl}
& \left| P_{f_{\bm{\vec \theta}}}(\sigma^{**}(\bm{x})) - P_{g_{\theta_n}, \sigma}(\sigma^{**}(\bm{x})) \right| \, d \bm{x} \\
= & \left| \sigma^{**}(\bm{\theta_1}) \sigma^{**}(\bm{x}) + \cdots + \sigma^{**}(\bm{\theta_{n-1}}) \sigma^{**}(\bm{x}^{n-1}) \right| \\
\lesssim_{K,B} & \left| (c_1 Q_m) (2a Q_m)^{-1} + \cdots + (c_1 Q_m)^{n-1} (2a Q_m)^{n-1} \right| \\
\lesssim_{K,B} & \left| 2a c_1 + \cdots + (2a c_1)^{n-1} \right| \\
\lesssim_{K,B} & 2a c_1
\end{IEEEeqnarray*}
Hence 
\begin{IEEEeqnarray*}{rCl}
|E| & \lesssim_{K,B} & (2a Q_m^{-1})^k \cdot 2a c_1 \\
& \lesssim_{K,B} & a^{k+1} c_1 Q_m^{-k}.
\end{IEEEeqnarray*}
\end{proof}

Using Lemma \ref{lem:Elem} and the triangle inequality, we have the estimate
\[| \mathrm{II}_{\bm{r}}^1(\bm{x})| \geq \left| \int_{\mathbb{R}^k} e(\phi_{g_{\bm{\theta_n}}}(\bm{x})) \chi_m(\bm{x}) \right| \, d \bm{x} - C_{K,B,m} a^{k-1} c_1 Q_m^{-k}.\]
So we are left to obtain a lower bound on the main term. Note that this term does not depend on $\bm{\vec{\theta}'}$ and depends only on $\bm{\theta_n}$. We write
\[\mathrm{II}_{\bm{r}}^{\text{main}}(\bm{\theta_n}) = \int_{\mathbb{R}^k} e(\phi_{g_{\bm{\theta_n}}}(\bm{x})) \chi_m(\bm{x}) \, d \bm{x}.\]
Let $\Psi : \mathbb{R}^k \to \mathbb{R}_{\geq 0}$ be a smooth function supported on $\mathbb{B}(\bm{0}; 1)$ with $\Psi(\bm{x}) = 1$ for $|\bm{x}| \leq \frac{1}{2}$. In order to obtain a lower bound on $\mathrm{II}_{\bm{r}}^{\text{main}}(\bm{\theta_n})$, we begin by writing it as $A - B$, where
\[A :=  := \lim_{R \to \infty} A_R := \lim_{R \to \infty} \int_{\mathbf{B}(\bm{0}, R)} e(\phi_{g_{\bm{\theta_n}}}(\bm{x})) \Psi(R^{-1} \bm{x}) \, d \bm{x} \]
and
\[B := \lim_{R \to \infty} B_R := \lim_{R \to \infty} \int_{\mathbf{B}(\bm{0}, R)} e( \phi_{g_{\bm{\theta_n}}}(\bm{x})) \left[1 - \chi\left(\frac{Q_m}{a} \bm{x}\right)\right] \Psi(R^{-1} \bm{x}) \, d \bm{x}\]
We will first obtain an upper bound on $B_R$. 
\begin{mylem}\label{lem:BRlem}
Suppose $B_R$ is defined as above. Then 
\[|B_R| \lesssim_{K,B,n} a^{-(n-1)} Q_m^{-k}.\]
Here, the implicit constant is independent of $R$.
\end{mylem}
\begin{proof}[Proof of Lemma \ref{lem:BRlem}]
We make the substitution $\bm{w} = R^{-1} \bm{x}$.  Then
\[B_R = R^k \int e(\phi_{g_{\bm{\theta_n}}}(R \bm{w}))\left[1 - \chi (\frac{Q_m}{a} R \bm{w}) \right] \Psi(\bm{w}) \, d \bm{w}.\]
Write $\Upsilon(\bm{w}) = \left[1 - \chi \left(\frac{Q_m R}{a} \bm{w} \right) \right] \Psi(\bm{w})$. Defining $\Delta_{\sigma}$ as before, we define $\Upsilon_{\sigma}(\bm{w}) = \Upsilon(\bm{w}) \Delta_{\sigma}(\bm{w})$.
Then, for every integer $N \geq 1$ we have 
\[\norm{\Upsilon_{\sigma}}_{C^N} \lesssim \left(\frac{Q_m R}{a} \right)^N\]
and we have for $\sigma = \sigma_0$ that 
\begin{IEEEeqnarray*}{Cl}
& H_{g_{\bm{\theta_n}}, K, B, \sigma_0, \Upsilon_{\sigma_0}} \\
\geq & \inf_{\bm{w} \in \supp \Upsilon_{\sigma_0}} \partial \left[\sigma^{**}(\bm{\theta_n}) \sigma^{**}(R \bm{w})^n \right] \\
\gtrsim & \inf_{\bm{w} \in \supp \Upsilon_{\sigma_0}} R^n Q_m^N |\bm{w}|^{n-1} \\
\gtrsim & a^{n-1} Q_m R
\end{IEEEeqnarray*}
and
\begin{IEEEeqnarray*}{Cl}
& J_{g_{\bm{\theta_n}}, K, B, \sigma_0, \Upsilon_{\sigma_0}} \\
\geq & \left[ \inf_{\bm{w} \in \supp \Upsilon_{\sigma_0}} \partial^2 \left[ \sigma^{**}(\bm{\theta_n}) \sigma^{**}(R \bm{w})^n \right] \right]^{1/2}  \\
\gtrsim & \left[ \inf_{\bm{w} \supp \Upsilon_{\sigma_0}} R^n Q_m^n |\bm{w}|^{n-2} \right]^{1/2} \\
\gtrsim & \left(a^{n-2} Q_m^2 R^2 \right)^{1/2} \\
= & a^{\frac{n-2}{2}} Q_m R
\end{IEEEeqnarray*}
For $\sigma \neq \sigma_0$, we use the $n$th derivative to obtain a bound on the $H$- and $J$- functionals.
\begin{IEEEeqnarray*}{rCl}
H_{g_{\bm{\theta_n}}, K, B, \sigma, \Upsilon_{\sigma_0}} & \gtrsim & \left[\partial^n \sigma^{**}(\bm{\theta_n}) (R \bm{w})^n \right]^{1/n} \\
& \gtrsim & \left[R^n Q_m^n \right]^{1/n} \\
& = & R Q_m
\end{IEEEeqnarray*}
with a similar estimate for the $J$-functional. Hence
\[H_{g_{\bm{\theta_n}, K, B, \Sigma, \Upsilon_{\sigma_0}}} \gtrsim a^{n-1} (Q_m R)^k. \]
Since $J_{g_{\bm{\theta_n}, K, B, \sigma, \Upsilon_{\sigma_0}}} \gg \norm{\Upsilon_{\sigma_0}}_{C^N}^{1/N}$ for every $N$ and every $\sigma$, we have from Corollary \ref{cor:mainthmcor} that 
\[|B_R| \lesssim_{K, B, n} a^{-(n-1)} Q_m^{-k}.\] 
This estimate is independent of $R$.
\end{proof}

Finally, we need to obtain a lower bound on the integral
\[A_R = \int_{\mathbb{R}^k} e(\phi_{g_{\bm{\theta_n}}}(x)) \Psi(R^{-1}\bm{x}) \, d \bm{x}.\]
\begin{mylem}\label{lem:ARlem}
Let $A_R$ be defined as above. Then 
\[|A_R| \gtrsim_{K,B,n} Q_m^{-k}.\]
\end{mylem}
\begin{proof}[Proof of Lemma \ref{lem:ARlem}]
Recall that the phase $g_{\bm{\theta_n}}$ is given by
\[\tr^*(g_{\bm{\theta_n}}(A^*\bm{x})) =  T^{-1}(\bm{\theta_n}) \cdot (A^*)^{-1}(A^*(\bm{x}))^n.\]
For any $\bm{\theta_n}$, this is a homogeneous polynomial in $\bm{x}$ of degree $n$ whose largest coefficient is on the order of $Q_m^n$. Write $\mathcal{Q}_n(Q)$ for the space of homogeneous polynomials of degree $n$ whose largest coefficient is on the order of $Q$. We will show that
\[\int_{\mathbb{R}^k} e(P(\bm{x})) \Psi(R^{-1} \bm{x}) \, d \bm{x} \gtrsim Q_m^{-k}\]
uniformly for $R \geq 1$ and for $P \in \mathcal{Q}_n(Q_m^n)$. In fact, by a rescaling argument, it is enough to show that
\[\int_{\mathbb{R}^k} e(P(\bm{x})) \Psi(R^{-1} \bm{x}) \, d \bm{x} \gtrsim 1\]
for $R \geq 1$ and $P \in \mathcal{Q}_n(1)$.
In fact, we will show the stronger statement that for any unit vector $\bm{u}$ that 
\[\int_{t >0} \text{Re}(e(P(t \bm{u}))) \Psi(R^{-1}(t \bm{u})) \, dt \gtrsim 1.\]
It is to be understood that the left side of the inequality is positive. Here, the implicit constant is uniform over sufficiently large $R$ and for $P \in \mathcal{Q}_n(1)$.

Indeed, this integral is of the form
\[\int_{t > 0} \text{Re}(e(C_{P, \bm{u}} t^n)) \Psi(R^{-1} t \bm{u}) \, dt.\]
Here, $C_{P, \bm{u}}$ is a constant such that $C_{P, \bm{u}} \lesssim 1$.

But this integral is 
\[\int_{t > 0} \cos(C_{P, \bm{u}} t^n) \Psi(R^{-1} t \bm{u}) \, dt.\]
We may assume $C_{P,\bm{u}} \geq 0$ since the cosine function is even. If $C_{P,\bm{u}} = 0$, this integral is clearly positive (and in fact approaches $\infty$ in the limit as $R \to \infty$). If $C_{P, \bm{u}} \neq 0$, then we make the change-of-variables $s = C_{P, \bm{u}}^{1/n} t$. This gives
\[C_{P,\bm{u}}^{-1/n} \int_{s > 0} \cos(s^d) \Psi(R^{-1} C_{P,\bm{u}} s \bm{u}) \, ds\]
But the series
\[\sum_{j=0}^{\infty} \int_{s \in [(j \pi)^{1/n}, ((j + 1) \pi)^{1/n}]} \cos(s^d)\]
is seen to satisfy the conditions of the alternating series test and therefore converges to a strictly positive number. Since $C_{P,\bm{u}}^{-1/n} \gtrsim 1$, the desired estimate holds uniformly for $P \in \mathcal{Q}_n(1)$ and for large $R$.
\end{proof}
\begin{proof}[Proof of Lemma \ref{lem:IIlem}]
Combining Lemmas \ref{lem:II2lem}, \ref{lem:Elem}, \ref{lem:BRlem}, and \ref{lem:ARlem}, we see that 
\begin{IEEEeqnarray*}{rCl}
\left|\mathrm{II}_{\bm{r}}(\bm{\theta}) \right| & \gtrsim_{K,B,n, \psi} &  Q_m^{-k} - a^{-\frac{1} {n - 1}} Q_m^{-k} - c_1 a^{k+1} Q_m^{-k} - a^{-(n-1)} Q_m^{-k} \\
& \gtrsim & Q_m^{-k},
\end{IEEEeqnarray*}
where the last step follows from the fact that $a$ is chosen to be large, and the fact that $c_1 a^{k+1}$ is small by \eqref{eq:c1acondition}.
\end{proof}

\section{Appendix: Integration by parts}
We carry out the integration by parts described in Section \ref{sec:mainthm}. Throughout this section, we will assume that $e \in E_{\bm{\vec \alpha}, S_1, S_2, S_3, \bm{r_1}, \bm{r_2}, \bm{l}}$, where $E_{\bm{\vec \alpha}, S_1, S_2, S_3, \bm{r_1}, \bm{r_2}, \bm{l}}$ is as in Section \ref{sec:mainthm}.

\begin{mylem}[Integration by parts]\label{lem:intbyparts}
Let $\chi$ and $\Psi$ be smooth functions with $\chi$ supported on the unit box. Let $I_{\chi, \Psi}$ denote the integral
\[\int_{\bm{\vec y}} \chi(3 \epsilon \bm{\vec y}) \Psi \left(r_{\sigma_j} (\bm{\vec x_e}) + \sum_{j=1}^{\tilde k} r_{\tilde \sigma_j}(\bm \vec{x_e}) \bm{\vec y_{\tilde \sigma_j}}\right) \prod_{j=1}^{\tilde k}  \prod_{\sigma \in \Sigma} e \left(P_{f, \sigma} (\vec \sigma (\bm{\vec x_e} + r_{\tilde \sigma}(\bm{\vec x_e}) \vec y_{\tilde \sigma})) \right) \, d \bm{\vec y}.\]
Suppose $\tilde \sigma \in \tilde S_2 \cup \tilde S_3$, and let $N \geq 1$. Suppose $\Lambda_{\tilde \sigma,e}$ is defined by
\[\Lambda_{\tilde \sigma, e} = |\nabla P_{f,\tilde \sigma}(\bm{\vec x_e})|.\]
Then there exist finite collections of smooth functions $\chi_1, \ldots, \chi_{C(N)}$ and $\Psi_1, \ldots, \Psi_{C(N)}$ such that for each $b > 0$, and $1 \leq i \leq C(N)$ there exists an exponent $0 \leq N_i \leq N$ such that
\[\norm{\Psi_i}_{C^b} \lesssim_{\chi, N, b, \epsilon} r_{\tilde \sigma}(\bm{\vec x_e})^{N_i} \norm{\Psi}_{C^{b + N_i}}\]
such that
\[|I_{\chi, \Psi}| \lesssim \Lambda^{-N} \sum_{i=1}^{C(N)} |I_{\chi_i, \Psi_i}|.\]
In particular, the implicit constants do not depend on $f, r(\bm{\vec x_e}),$ or on $\Lambda_{\tilde \sigma,e}$, or on $\Psi$.
\end{mylem}
\begin{proof} Suppose $|\nabla P_{f, \tilde \sigma}(\vec \sigma(\bm{\vec x_e}))| \geq C J_{f, K, B, \sigma}(\bm{\vec x_e})$ for some large constant $C$. Then we have by Lemma \ref{lem:gradstability} for $\bm{\vec x} \in \mathbf{P}_{S,1}(\bm{\vec x_e})$ that 
\[|\nabla P_{f, \tilde \sigma}(\vec \sigma(\bm{\vec x})) - \nabla P_{f, \tilde \sigma}(\vec \sigma(\bm{\vec x_e}))| \leq \frac{1}{100} |\nabla P_{f, \tilde \sigma}(\vec \sigma(\bm{\vec x_e}))|.\]

Fix $\tilde \sigma \in \tilde S_2 \cup \tilde S_3$. 

By Lemma \ref{lem:realpart}, we can select for each $\tilde \sigma \in \tilde S_2 \cup \tilde S_3$ a unit vector $\bm{\vec z_{\tilde \sigma}} \in V_{\tilde \sigma, \mathbb{R}}^n$ such that
\begin{equation}\label{eq:intbypartsrealpart}
\realpart (\nabla P_{f, \tilde \sigma}(\vec \sigma(\bm{\vec x_e})) \cdot \vec \sigma(\bm{\vec z_{\tilde \sigma}})) \geq c |\nabla P_{f, \vec \sigma}(\vec \sigma(\bm{\vec x_e}))|.
\end{equation}
Then, for a different constant $c$, we have the following estimate for any vector $\bm{\vec z_{\tilde \sigma}}^{\perp}$ of length at most $1$ orthogonal to $\bm{\vec z_{\tilde \sigma}}$ lying in $V_{\tilde \sigma, \mathbb{R}}^n$:
\begin{equation}\label{eq:intbypartsrealparteverywhere}
\realpart (\nabla P_{f, \tilde \sigma}(\vec \sigma(\bm{\vec x_e} + r_{\tilde \sigma}(\bm{\vec x_e}) \bm{\vec z_{\tilde \sigma}^{\perp}})) \cdot \vec \sigma(\bm{\vec z_{\tilde \sigma}})) \geq c |\nabla P_{f, \vec \sigma}(\vec \sigma(\bm{\vec x_e} + r_{\tilde \sigma}(\bm{\vec x_e}) \bm{\vec z_{\tilde \sigma}^{\perp}}))|.
\end{equation}

We will apply integration by parts to the integral in $\vec y_{\sigma}$ only. We can pull the factors of $\realpart P_{f, \tilde \sigma_j} (\vec \sigma_j(\bm{\vec x_e} + r_{\tilde \sigma_j} (\bm{\vec x_e}) \bm{\vec y_{\tilde \sigma_j}}))$ for $\tilde \sigma_j \neq \tilde \sigma$ out of the integral and focus only on the integral
\[\int_{\vec y_{\tilde \sigma}} \chi \left(3 \epsilon \bm{\vec y}\right) \Psi \left(\bm{\vec x_e} + \sum_{j=1}^{\tilde k}  r_{\tilde \sigma_j}(\bm{\vec x_e}) \vec y_{\tilde \sigma_j} \right) e \left(C_{\sigma} \realpart(P_{f, \sigma} (\vec \sigma (\bm{\vec x_e} + r_{\tilde \sigma}(\bm{\vec x_e}) \vec y_{\tilde \sigma})))) \right) \, d \vec y_{\sigma}. \]
We only want to integrate by parts in the $\bm{\vec z_{\tilde \sigma}}$ direction, so will rewrite the phase to highlight the dependence on this variable. Write $\vec y_{\tilde \sigma} = a \bm{\vec z_{\tilde \sigma}} + \bm{\vec z_{\tilde \sigma}}^{\perp}$, where $a \in [-1,1]$ and $|\bm{\vec z_{\tilde \sigma}}^{\perp}| \leq 3 \epsilon$. Then the integral becomes
\begin{IEEEeqnarray*}{cl}
\int_{\bm{\vec z_{\tilde \sigma}}^{\perp}} \int_{a} & \chi(3 \epsilon (\bm{\vec z_{\tilde \sigma}}^{\perp} + a \bm{\vec z_{\tilde \sigma}})) \Psi \left(\bm{\vec x_e} + \sum_{\tilde \sigma_j \neq \tilde \sigma} r_{\tilde \sigma_j} (\bm{\vec x_e}) \bm{\vec y_{\tilde \sigma_j}} + r_{\tilde \sigma}(\bm{\vec x_e}) \bm{\vec z_{\tilde \sigma}}^{\perp} + r_{\tilde \sigma}(\bm{\vec x_e}) a \bm_{\vec z_{\tilde \sigma}}\right) \\
& \cdot e\left (C_{\sigma} \realpart(P_{f, \sigma}(\vec \sigma(\bm{\vec x_e} + r_{\tilde \sigma}(\bm{\vec x_e}) \bm{\vec z_{\tilde \sigma}}^{\perp} + a r_{\tilde \sigma}(\bm{\vec x_e}) \bm{\vec z_{\tilde \sigma}}))) \right)d a d \bm{\vec z_{\tilde \sigma}}^{\perp}.
\end{IEEEeqnarray*}
We will write the Taylor expansion of the phase in the $a$ variable, centered at the point $\vec \sigma(\bm{\vec x_e} + r_{\tilde \sigma}(\bm{\vec x_e}) \bm{\vec z_{\tilde \sigma}}^{\perp}).$ We see that
\begin{IEEEeqnarray*}{CCl}
& & \realpart(P_{f, \sigma}(\vec \sigma(\bm{\vec x_e} + r_{\tilde \sigma}(\bm{\vec x_e}) \bm{\vec z_{\tilde \sigma}}^{\perp} + a r_{\tilde \sigma}(\bm{\vec x_e}) \bm{\vec z_{\tilde \sigma}}))) \\
= & & \realpart(P_{f, \sigma}(\vec \sigma(\bm{\vec x_e} + r_{\tilde \sigma}(\bm{\vec x_e}) \bm{\vec z_{\tilde \sigma}}^{\perp}))) + a r_{\tilde \sigma}(\bm{\vec x_e}) \realpart \nabla_{\vec \sigma(\bm{\vec z_{\tilde \sigma})}} P_{f, \sigma}(\vec \sigma(\bm{\vec x_e} + r_{\tilde \sigma}(\bm{\vec x_e}) \bm{\vec z_{\tilde \sigma}}^{\perp})) \\
& + & \sum_{\alpha \geq 2} \frac{1}{\alpha!} (a r_{\tilde \sigma}(\bm{ \vec x_e}))^{\alpha} \realpart \left( (\nabla_{\vec \sigma(\bm{\vec z_{\tilde \sigma}})})^{\alpha} P_{f, \sigma}(\vec \sigma(\bm{\vec x_e} + r_{\tilde \sigma}(\bm{\vec x_e}) \bm{\vec z_{\tilde \sigma}}^{\perp}))\right).
\end{IEEEeqnarray*}
For clarity, we emphasize that $\alpha$ above is a one-dimensional index, not a multi-index.

Observe that 
\[\left| \frac{1}{\alpha!} (\nabla_{\vec \sigma(\bm{\vec z_{\tilde \sigma}})})^{\alpha} P_{f, \sigma}(\vec \sigma (\bm{\vec x_e} + r_{\tilde \sigma}(\bm{\vec x_e)} \bm{\vec z_{\tilde \sigma}}^{\perp})) \right| \lesssim_n J_{f,K,B,\sigma}(\bm{\vec x_e} + r_{\tilde \sigma}(\bm{\vec x_e}) \bm{\vec z_{\tilde \sigma}}^{\perp})^{\alpha} \lesssim_{K,B,n} J_{f,K,B, \sigma}(\bm{\vec x_e})\]
where the last inequality follows since $|\bm{\vec z_{\tilde \sigma}}^{\perp}| \leq 3 \epsilon$ on the support of $\gamma$. Recall that we chose $r_{\tilde \sigma}(\bm{\vec x_e}) = \frac{1}{J_{f,K,B, \sigma}(\bm{\vec x_e})}$; therefore, the coefficients on $a^{\alpha}$ are bounded above by a constant. Moreover, writing $\Lambda_{e, \sigma}(\bm{\vec z_{\tilde \sigma}}^{\perp})$ for the coefficient on $a$; that is, $r_{\tilde \sigma}(\bm{\vec x_e}) \realpart \nabla_{\vec \sigma(\bm{\vec z_{\tilde \sigma}})} P_{f, \sigma} (\vec \sigma(\bm{\vec x_e} + r_{\tilde \sigma}(\bm{\vec x_e}))\bm{\vec z_{\tilde \sigma}}^{\perp})$, we observe from the fact that $\tilde \sigma \in \tilde S_2 \cup \tilde S_3$ and the estimate \eqref{eq:intbypartsrealparteverywhere} that $\Lambda_{e, \tilde \sigma}(\mathbf{z}_{\tilde \sigma}^{\perp}) \gtrsim 1$. Hence the phase takes the form 

\[\realpart(P_{f, \sigma}(\vec \sigma(\bm{\vec x_e} + r_{\tilde \sigma}(\bm{\vec x_e}) \bm{\vec z_{\tilde \sigma}}^{\perp}))) + \Lambda_{e, \tilde \sigma}(\bm{\vec z_{\tilde \sigma}}^{\perp}) \cdot \left( a + \sum_{2 \leq \alpha \leq d}C_{e, \tilde \sigma, \alpha}(\bm{\vec z_{\tilde \sigma}}) a^{\alpha} \right),\]
where each coefficient $|C_{e, \tilde \sigma, \alpha}| \lesssim 1$ and $\Lambda_{e, \tilde \sigma}(\bm{\vec z_{\tilde \sigma}}) \gtrsim 1.$ Write $Q(a)$ for the polynomial $\left( a + \sum_{2 \leq \alpha \leq d}C_{e, \tilde \sigma, \alpha}(\bm{\vec z_{\tilde \sigma}}) a^{\alpha} \right)$. Then $Q'(a) \gtrsim 1$ and $Q^{\beta}(a) \lesssim 1$ for all $a \in [-6 \epsilon, 6 \epsilon]$.
Thus the integral becomes
\begin{IEEEeqnarray*}{ccl}
\int_{\bm{\vec z_{\tilde \sigma}}^{\perp}} & & e \left(\realpart(P_{f, \sigma}(\vec \sigma(\bm{\vec x_e} + r_{\tilde \sigma}(\bm{\vec x_e}) \bm{\vec z_{\tilde \sigma}^{\perp}}))) \right) \cdot \\
& \int_{a} &\chi \left(3 \epsilon (\sum_{\tilde \sigma_j \neq \tilde \sigma} \bm{\vec y_{\tilde \sigma}} + \bm{\vec z_{\tilde \sigma}}^{\perp} + a \bm{\vec z_{\tilde \sigma}}) \right) \Psi  \left(\bm{\vec x_e} + \sum_{\tilde \sigma_j \neq \tilde \sigma} r_{\tilde \sigma_j} (\bm{\vec x_e}) \bm{\vec y_{\tilde \sigma_j}} + r_{\tilde \sigma}(\bm{\vec x_e}) \bm{\vec z_{\tilde \sigma}}^{\perp} + r_{\tilde \sigma}(\bm{\vec x_e}) a \bm{\vec z_{\tilde \sigma}}\right) \cdot \\
& & \cdot e\left (C_{\sigma} \Lambda_{e, \tilde \sigma}(\bm{\vec z_{\tilde \sigma}}^{\perp}) \cdot Q(a) \right)\,da \, d \bm{\vec z_{\tilde \sigma}}^{\perp}.
\end{IEEEeqnarray*}
We are now ready to integrate by parts $N$ times in the $a$-variable. Defining the differential operator $D$ by 
\[Df(a) := \frac{1}{2 \pi i Q'(a)} \frac{d}{da} f,\]
we see that the adjoint $D^t$ of $D$ is given by
\[D^t f = \frac{1}{2 \pi i}\frac{d}{da} \left[ \frac{1}{Q'(a)} f \right].\]
The integral in $a$ becomes
\begin{IEEEeqnarray*}{ccl}
& & C_N (\Lambda_{\tilde \sigma, e}(\bm{\vec z_{\tilde \sigma}}^{\perp}))^{-N} \cdot  \\
\cdot & \int_{a} & \chi \left(3 \epsilon (\sum_{\tilde \sigma_j \neq \tilde \sigma} \bm{\vec y_{\tilde \sigma}} + \bm{\vec z_{\tilde \sigma}}^{\perp} + a \bm{\vec z_{\tilde \sigma}}) \right) \Psi  \left(\bm{\vec x_e} + \sum_{\tilde \sigma_j \neq \tilde \sigma} r_{\tilde \sigma_j} (\bm{\vec x_e}) \bm{\vec y_{\tilde \sigma_j}} + r_{\tilde \sigma}(\bm{\vec x_e}) \bm{\vec z_{\tilde \sigma}}^{\perp} + r_{\tilde \sigma}(\bm{\vec x_e}) a \bm{\vec z_{\tilde \sigma}}\right) \cdot \\
& & \cdot D^n  e\left (C_{\sigma} \Lambda_{e, \tilde \sigma}(\bm{\vec z_{\tilde \sigma}}^{\perp}) \cdot Q(a) \right)\,da \, d \bm{\vec z_{\tilde \sigma}}^{\perp}. \\
= & & C_N (\Lambda_{\tilde \sigma, e} (\bm{\vec z_{\tilde \sigma}}^{\perp}))^{-N} \cdot \\
\cdot  &  \int_{a} & (D^t)^N \left[ \chi \left(3 \epsilon (\sum_{\tilde \sigma_j \neq \tilde \sigma} \bm{\vec y_{\tilde \sigma}} + \bm{\vec z_{\tilde \sigma}}^{\perp} + a \bm{\vec z_{\tilde \sigma}}) \right) \Psi  \left(\bm{\vec x_e} + \sum_{\tilde \sigma_j \neq \tilde \sigma} r_{\tilde \sigma_j} (\bm{\vec x_e}) \bm{\vec y_{\tilde \sigma_j}} + r_{\tilde \sigma}(\bm{\vec x_e}) \bm{\vec z_{\tilde \sigma}}^{\perp} + r_{\tilde \sigma}(\bm{\vec x_e}) a \bm{\vec z_{\tilde \sigma}}\right) \right] \cdot \\
& & \cdot e\left (C_{\sigma} \Lambda_{e, \tilde \sigma}(\bm{\vec z_{\tilde \sigma}}^{\perp}) \cdot Q(a) \right)\,da \, d \bm{\vec z_{\tilde \sigma}}^{\perp}. \\
\end{IEEEeqnarray*}
It only remains to compute
\[(D^t)^N \left[ \chi \left(3 \epsilon (\sum_{\tilde \sigma_j \neq \tilde \sigma} \bm{\vec y_{\tilde \sigma_j}} + \bm{\vec z_{\tilde \sigma}}^{\perp} + a \bm{\vec z_{\tilde \sigma}}) \right) \Psi  \left(\bm{\vec x_e} + \sum_{\tilde \sigma_j \neq \tilde \sigma} r_{\tilde \sigma_j} (\bm{\vec x_e}) \bm{\vec y_{\tilde \sigma_j}} + r_{\tilde \sigma}(\bm{\vec x_e}) \bm{\vec z_{\tilde \sigma}}^{\perp} + r_{\tilde \sigma}(\bm{\vec x_e}) a \bm{\vec z_{\tilde \sigma}}\right) \right] \cdot \]
But this derivative is of the form
\begin{IEEEeqnarray*}{CCl}
\sum_{\beta_1 + \beta_2 + \beta_3 = N} & & \left(\frac{d}{da} \right)^{\beta_1} \frac{1}{Q'}(a) \\
& \cdot & (3 \epsilon)^{\beta_2} (\nabla_{\bm{\vec z_{\tilde \sigma}}})^{\beta_2} \chi\left(3 \epsilon \left( \sum_{\tilde \sigma_j \neq \sigma} \bm{\vec y_{\tilde \sigma_j}} \bm{\vec z_{\tilde \sigma}}^{\perp} + a \bm{\vec z_{\tilde \sigma}}\right) \right)\\
& \cdot & r_{\tilde \sigma}(\bm{\vec x_e})^{\beta_3} (\nabla_{\bm{z_{\tilde \sigma}}})^{\beta_3} \Psi \left(\bm{\vec x_e} + \sum_{\tilde \sigma_j \neq \tilde \sigma} r_{\tilde \sigma_j} (\bm{\vec x_e}) \bm{\vec y_{\tilde \sigma_j}} + r_{\tilde \sigma}(\bm{\vec x_e}) \bm{\vec z_{\tilde \sigma}}^{\perp} + r_{\tilde \sigma}(\bm{\vec x_e}) a \bm{\vec z_{\tilde \sigma}}\right).
\end{IEEEeqnarray*}
Observe that each summand is of the desired form, since $\left(\frac{d}{da} \right)^{N} \frac{1}{Q'}(\alpha)$ is bounded above by a constant depending only on $N$ by our observation that $Q'(a) \sim 1$ and $Q^{(i)}(a) \lesssim 1$ on the relevant domain $[-6 \epsilon, 6 \epsilon]$ for $i \geq 1$. Observe that the coefficient $r_{\tilde \sigma}(\bm \vec{x_e})^{\beta_3}$ on each term is always bounded above by $1 + r_{\tilde \sigma}(\bm \vec{x_e})^N$ and that the quantity $\Lambda_{e, \tilde \sigma}(\bm{\vec z_{\tilde \sigma}}) \gtrsim \Lambda_{e, \tilde \sigma}$ for $\norm{z_{\tilde \sigma}} \leq 3 \epsilon$.
\end{proof}
\begin{mycor}[Integration by parts, Several times]\label{cor:intbyparts}
Let $\chi$ and $\Psi$ be smooth functions with $\chi$ supported on the box $[-2,2]^{kn}$. Let $I_{\chi, \Psi}$ denote the integral
\[\int_{\bm{\vec y}} \chi((3 \epsilon)^{-1} \bm{\vec y}) \Psi \left(\bm{\vec x_e}  + \sum_{j=1}^{\tilde k} r_{\tilde \sigma_j}(\bm{\vec x_e}) \bm{\vec y_{\sigma_j}} \right) \prod_{j=1}^{\tilde k} \prod_{\sigma \in \Sigma} e \left(P_{f, \sigma}(\vec \sigma(\bm{\vec x_e} + r_{\tilde \sigma}(\bm{\vec x_e}) \bm{y}_{\tilde \sigma})) \right) \, d \bm{\vec y}.\]
Suppose $N \geq 1$. For each $\tilde \sigma \in \tilde S_2 \cup S_3$, suppose $\Lambda_{\tilde \sigma, e}$ is defined by 
\[\Lambda_{\tilde \sigma, e} = | \nabla P_{f, \tilde \sigma}(\bm{\vec x_e}) |.\]
Then we have the bound
\begin{equation}\label{eq:intbypartsgeneral}
|I_{\chi, \Psi}| \lesssim_{\chi, \Psi, n, N, d, \epsilon} \prod_{\tilde \sigma \in \tilde S_2 \cup \tilde S_3} \left( \Lambda_{\tilde \sigma, e}^{-N} \left(1 + r_{\tilde \sigma}(\bm{\vec x_e})^N \right) \right).
\end{equation}
If, in addition, there exists a positive constant $Z$ such that $\Psi$ satisfies the condition that for all $\alpha \geq 0$, we have the bound
\begin{equation}\label{eq:intbypartscond}
\norm{\psi}_{C^{\alpha}} \leq Z \min_{\tilde \sigma \in \tilde S_2 \cup \tilde S_3} \left( r_{\tilde \sigma}(\bm{\vec x_e})^{-\alpha}  \right)
\end{equation}
for all $1 \leq \alpha \leq k N$, then $|I_{\chi, \Psi}|$ satisfies the stronger bound 
\begin{equation}\label{eq:intbypartsspecial}
I_{\chi, \Psi} \lesssim_{\chi, n, N, d, \epsilon, Z} \prod_{\tilde \sigma \in \tilde S_2 \cup \tilde S_3} \Lambda_{\tilde \sigma, e}^{-N}.
\end{equation}
In particular, the bound \eqref{eq:intbypartsspecial} does not depend on the function $\Psi$ beyond the dependence on $R$. 
\end{mycor}
\begin{proof}
By iteratively applying Lemma \ref{lem:intbyparts} to each $\tilde \sigma \in \tilde S_2 \cup \tilde S_3$, we arrive at a family of integrals $I_{\chi_i, \Psi_i}$, where each $\Psi_i$ satisfies a bound of the form
\begin{equation}\label{eq:psiibound}
\norm{\Psi_i}_{C^b} \lesssim_{\chi, N, b, \epsilon} \norm{\Psi}_{C^{b + N_i}} \prod_{\tilde \sigma \cup S_2 \cup S_3} r_{\tilde \sigma}(\bm{\vec x_e})^{N_{i, \tilde \sigma}},
\end{equation}
for some $0 \leq N_1 \leq |\tilde S_2 \cup \tilde S_3| N$, and some $0 \leq N_{i, \tilde \sigma} \leq N$ that sum to $N_i$.

Applying the trivial bound $|I_{\chi_i, \Psi_i}| \lesssim_{\chi} \norm{\Psi}_{C^0}$, we get that
\begin{equation}\label{eq:Ichipsibound}
|I_{\chi_i, \Psi_i}| \lesssim_{\chi, N, \epsilon, n ,d} \norm{\Psi}_{C^{N_i}} \prod_{\tilde \sigma \in \tilde S_2 \cup \tilde S_3} \Lambda_{\tilde \sigma}^{-N} r_{\tilde \sigma}(\bm{\vec x_e})^{N_{i, \tilde \sigma}}.
\end{equation}
In any case, we have that 
\[|I_{\chi, \Psi_i}| \lesssim_{\chi, N, \epsilon, n , d, \Psi} \prod_{\tilde \sigma \in \tilde S_2 \cup \tilde S_3} \Lambda_{\tilde \sigma}^{-N} (1 + r_{\tilde \sigma}(\bm{\vec x_e})^N).\]
If we also know that the condition \eqref{eq:intbypartscond} holds, then the bound \eqref{eq:Ichipsibound} implies 
\[|I_{\chi_i, \Psi_i}| \lesssim_{\chi, n, N, d, \epsilon, R} \prod_{\tilde \sigma \in \tilde S_2 \cup \tilde S_3} \Lambda_{\tilde \sigma}^{-N},\]
as desired.
\end{proof}
\bibliographystyle{plain}
\bibliography{Oscillatory_Integrals_ANT}
\end{document}